\documentclass[12pt, leqno]{amsart}

\usepackage{shuffle} 

\usepackage{graphicx}
\usepackage{amssymb,amsmath,amsthm,amscd}
\usepackage{mathrsfs}
\usepackage{enumerate}
\usepackage[usenames,dvipsnames]{color}
\usepackage[colorlinks=true, pdfstartview=FitV, linkcolor=blue,citecolor=blue,urlcolor=blue]{hyperref}

\usepackage[all]{xy}
\usepackage[normalem]{ulem}  
\usepackage{enumitem}

\allowdisplaybreaks[2]
\usepackage{oldgerm}
\usepackage{xspace}
\usepackage{marginnote}

\usepackage{array}
\usepackage{tikz}
\usepackage{enumitem}

\newcommand{\nc}{\newcommand}
\numberwithin{equation}{section}

\newenvironment{red}
{\relax\color{red}}
{\hspace*{.5ex}\relax}

\newcommand{\ber}{\begin{red}}
\newcommand{\er}{\end{red}}

\nc{\berE}{\ber{}\marginnote{\fbox{\scshape\lowercase{E}}}}
\newenvironment{blue}
{\relax\color{Blue}}
{\hspace*{.5ex}\relax}

\newcommand{\beb}{\begin{blue}}
\newcommand{\eb}{\end{blue}}

\newenvironment{yellow}
{\relax\color{Dandelion}}
{\hspace*{.5ex}\relax}

\newcommand{\bey}{\begin{yellow}}
\newcommand{\ey}{\end{yellow}}

\renewcommand{\le}{\leqslant}
\renewcommand{\ge}{\geqslant}

\theoremstyle{plain}
\newtheorem{lemma}{Lemma}[section]
\newtheorem{prop}[lemma]{Proposition}
\newtheorem{theorem}[lemma]{Theorem}

\newcommand{\Prop}{\begin{prop}}
\newcommand{\enprop}{\end{prop}}
\newcommand{\Lemma}{\begin{lemma}}
\newcommand{\enlemma}{\end{lemma}}
\newcommand{\Th}{\begin{theorem}}
\newcommand{\enth}{\end{theorem}}
\newtheorem{corollary}[lemma]{Corollary}
\newcommand{\Cor}{\begin{corollary}}
\newcommand{\encor}{\end{corollary}}

\newtheorem{sublemma}[lemma]{Sublemma}
\newcommand{\Sublemma}{\begin{sublemma}}
\newcommand{\ensub}{\end{sublemma}}

\theoremstyle{definition}

\newtheorem{remark}[lemma]{Remark}
\newtheorem{example}[lemma]{Example}
\newtheorem{Convention}[lemma]{Convention}
\newtheorem{definition}[lemma]{Definition}
\newcommand{\Def}{\begin{definition}}
\newcommand{\edf}{\end{definition}}
\newcommand{\Conv}{\begin{Convention}}
\newcommand{\enconv}{\end{Convention}}

\newcommand{\C}{{\mathbb C}}
\newcommand{\Q}{\mathbb {Q}}

\newcommand{\Z}{{\mathbb Z}}
\newcommand{\A}{{\mathbb A}}

\newcommand{\D}{\mathscr{D}}

\newcommand{\one}{{\bf{1}}}

\newcommand{\seteq}{\mathbin{:=}}

\newcommand{\hd}{{\operatorname{hd}}}
\newcommand{\soc}{{\operatorname{soc}}}

\newcommand{\g}{{\mathfrak{g}}}

\newcommand{\Hom}{\operatorname{Hom}}

\newcommand{\M}{{\mathscr M}}

\newcommand{\hs}{\hspace*}
\newcommand{\ms}{\mspace}

\newcommand{\To}[1][{\hs{2ex}}]{\xrightarrow{\,#1\,}}

\newcommand{\on}{\operatorname}

\newcommand{\bni}{\begin{enumerate}[label=\rm(\roman*)]}
\newcommand{\bnum}{\bni}
\newcommand{\bna}{\begin{enumerate}[label=\rm(\alph*)]}
\nc{\ee}{\end{enumerate}}

\newcommand{\id}{\on{id}}
\newcommand{\ba}{\begin{array}}
\newcommand{\ea}{\end{array}}

\newcommand{\supp}{\operatorname{supp}}

\newcommand{\eqsub}{\begin{subequations}\begin{eqnarray}}
\newcommand{\eneqsub}{\end{eqnarray}\end{subequations}}

\newcommand{\ol}{\overline}

\nc{\la}{\lambda}
\nc{\U}[1][\g]{U_q(#1)}
\nc{\te}{\tilde{e}}
\nc{\tem}{\tilde{e}^{\mathrm{max}}}
\nc{\tei}{\tilde{e}_i}
\nc{\tf}{\tilde{f}}
\nc{\tfm}{\tilde{f}^{\mathrm{max}}}
\nc{\tfi}{\tilde{f}_i}
\nc{\tU}{\widetilde U_q(\g)}
\nc{\tE}{\widetilde{E}}
\nc{\tF}{\widetilde{F}}
\nc{\tK}{\widetilde{K}}
\nc{\tEs}{\widetilde{E}^*}
\nc{\tFs}{\widetilde{F}^*}

\nc{\ttE}{\widetilde{\mathcal{E}}}
\nc{\ttF}{\widetilde{\mathcal{F}}}
\nc{\ttEs}{\ttE^*}
\nc{\ttFs}{\ttF^*}
\nc{\tfs}{\tf^*}
\nc{\tfss}[1]{\tf^{* \hskip 0.05em #1}}
\nc{\tess}[1]{\te^{* \hskip 0.05em #1}}
\nc{\tes}{\te^*}
\nc{\tesm}{\tilde{e}^{* \hskip 0.05em \mathrm{max}}}
\nc{\tfsm}{\tilde{f}^{* \hskip 0.05em \mathrm{max}}}

\nc{\tk}{\tilde{k}}
\nc{\tkone}{\tk_{\ol{1}}}
\nc{\teone}{\tilde{e}_{\ol{1}}}
\nc{\tfone}{\tilde{f}_{\ol{1}}}

\nc{\teibar}{\tilde{e}_{\ol{i}}} \nc{\tfibar}{\tilde{f}_{\ol{i}}}
\nc{\tki}{{\tk}_{\ol {i}}}

\nc{\BZ}{{\mathbb{Z}}}
\nc{\al}{\alpha}
\nc{\tal}{\widetilde{\al}}
\nc{\tch}{\widetilde{h}}

\nc{\lan}{\langle}
\nc{\ran}{\rangle}
\nc{\re}{{\mathrm{re}}}
\nc{\wt}{\operatorname{wt}}
\nc{\hwt}{\widehat{\wt}}
\nc{\Ht}{\mathrm{ht}}
\nc{\hHt}{\widehat{\Ht}}
\nc{\ch}{\operatorname{ch}}
\nc{\Um}[1][\g]{U^-_q(#1)}
\nc{\Ue}{U^+_q(\g)}
\nc{\ep}{\varepsilon}
\nc{\hep}{\widehat{\ep}}
\nc{\vphi}{\varphi}
\nc{\sphi}{\varphi^*}
\nc{\eps}{\ep^*}
\nc{\heps}{\hep^{ \hskip 0.2em  *}}

\nc{\nn}{\nonumber}

\nc{\vp}{\varpi}
\nc{\cls}{{\operatorname{cl}}}
\nc{\Wt}{{\operatorname{Wt}}}
\nc{\Us}{U'_q(\g)}
\nc{\La}{\Lambda}
\nc{\tLa}{\widetilde\Lambda}
\nc{\ro}{{\rm(}}
\nc{\rf}{{\rm)}}
\nc{\norm}{{\mathrm{norm}}}
\nc{\qbox}{\quad\mbox}
\nc{\braid}{{\mathfrak{B}}}
\nc{\Ad}{\operatorname{Ad}}
\nc{\Aut}{\operatorname{Aut}}
\nc{\dt}[1]{\tilde{\tilde #1}}
\nc{\Sn}{S^{{\mathrm{norm}}}}
\nc{\aff}{{\mathrm{aff}}}
\nc{\rk}{{\mathrm{rk}}}
\nc{\tP}{\widetilde{P}}
\nc{\tW}{\widetilde{W}}
\nc{\Dyn}{\mathrm{Dyn}}
\nc{\tD}{\widetilde{\Delta}}
\nc{\height}[1]{{\operatorname{ht}}(#1)}
\nc{\bl}{\bigl(}
\nc{\br}{\bigr)}
\nc{\Hecke}{\mathrm{H}}
\nc{\HA}{\Hecke^{\mathrm{A}}}
\nc{\HB}{\Hecke^{\mathrm{B}}}
\newcommand{\scbul}{{\,\raise1pt\hbox{$\scriptscriptstyle\bullet$}\,}}
\nc{\vac}{{\phi}}
\nc{\Bt}{\B_\theta(\g)}
\nc{\low}{{\mathrm{low}}}
\nc{\upper}{{\mathrm{up}}}
\nc{\Zodd}{\Z_{\mathrm{odd}}}
\nc{\Ft}[1][n]{\mathbb{P}\mathrm{ol}_{#1}}
\nc{\Ftf}[1][n]{\widetilde{\mathbb{P}\mathrm{ol}}_{#1}}
\nc{\KA}{\on{K}^{\mathrm{A}}}
\nc{\KB}{\on{K}^{\mathrm{B}}}
\nc{\Res}{\on{Res}}
\nc{\Fc}[1][{n,m}]{\mathbf{F}_{#1}}
\nc{\tphi}{\tilde{\varphi}}
\nc{\CO}{\mathscr{O}}
\nc{\inte}{\mathrm{int}}
\nc{\Oint}{\mathcal{O}^{\ge0}_{\inte}}
\nc{\vs}{\vspace*}
\nc{\tLt}{\widetilde{L}}
\nc{\tL}{\widetilde{\Lambda}}
\nc{\tu}{\tilde{u}}
\nc{\noi}{\noindent}
\nc{\heigh}{\mathfrak{t}}
\nc{\lowest}{\mathfrak{l}}
\nc{\rootl}{\mathsf{Q}}
\nc{\rlQ}{\rootl}

\nc{\KLR}{KLR algebra}
\nc{\KLRs}{KLR algebras}
\nc{\cor}{\mathbf{k}}
\nc{\cora}{{\cor(A)}}
\nc{\haut}{\mathrm{ht}}
\nc{\tens}{\mathop\otimes}
\nc{\gmod}{\mbox{-$\mathrm{gmod}$}}
\nc{\gMod}{\mbox{-$\mathrm{gMod}$}}
\nc{\proj}{\mbox{-$\mathrm{proj}$}}
\nc{\gproj}{\mbox{-$\mathrm{gproj}$}}
\nc{\smod}{\mbox{-$\mathrm{mod}$}}
\nc{\Mod}{\mbox{-$\mathrm{Mod}$}}
\nc{\h}{\mathfrak h}
\nc{\Rnorm}{R^{\mathrm{norm}}}
\nc{\Runiv}{R^{\mathrm{univ}}}
\nc{\Rren}{R^{\mathrm{ren}}}

\nc{\Vhat}{\widehat{V}}
\nc{\F}{\mathcal{F}}

\def\T{{\mathcal T}}

\def\Gup{\mathrm{G}^{\up}}
\def\tGup{\widetilde{\mathrm{G}}^{\up}}

\def\Binf{\mathrm{B}(\infty)}

\nc\Bc{\mathrm{B}}

\nc{\fd}[1][A]{\on{\mathrm{flat.dim}_{#1}}}
\nc{\bP}{{\mathbb{P}}}
\nc{\bPh}{\widehat{\mathbb{P}}}
\nc{\bK}[1][{n}]{\widehat{\mathbb{K}}_{#1}}
\nc{\bV}[1][{n}]{\widehat{V}^{\otimes{#1}}}
\nc{\bVK}[1][{n}]{\widehat{V}^{\otimes{#1}}_{\widehat{\mathbb{K}}}}
\nc{\hV}{\widehat{V}}
\nc{\opp}{\mathrm{opp}}
\nc{\col}{\colon}
\nc{\oep}{\epsilon}
\nc{\qtext}{\quad\text}
\nc{\qtextq}[1]{\quad\text{#1}\quad}
\nc{\longtwoheadrightarrow}[1][]{\xymatrix{\ar@{->>}[r]^-{{#1}}&}}
\nc{\epiTo}[1][]{\longtwoheadrightarrow[{#1}]}
\nc{\epito}{\twoheadrightarrow}
\nc{\monoTo}[1][]{\xymatrix{\ar@{>->}[r]^-{{#1}}&}}
\nc{\sym}{\mathfrak{S}}
\nc{\inp}[1]{{({#1})_{\mathrm{n}}}}
\nc{\rtl}{\rootl}
\nc{\wtd}{\widetilde}
\nc{\etens}{\boxtimes}
\nc{\ds}[1]{\mathrm{d}(#1)}
\nc{\rmat}[1]{{\mathbf{r}}_%
{\mspace{-2mu}\raisebox{-.6ex}{${\scriptstyle{#1}}$}}}
\nc{\rmats}[1]{{\mathbf{r}}_%
{\mspace{-2mu}\raisebox{-.6ex}{${\scriptscriptstyle{#1}}$}}}
\nc{\shc}{\mathcal{C}}
\nc{\shs}{\mathcal{S}}
\nc{\Fct}{{\on{Fct}}}
\nc{\tC}{\widetilde{\shc}}
\nc{\Zp}{\Z_{\ge0}}
\nc{\tPhi}{\widetilde{\Phi}}
\nc{\tT}{{\widetilde{\T}}}
\nc{\Ob}{\on{Ob}}
\nc{\bwr}{\mbox{\large$\wr$}}
\nc{\Img}{\on{Im}}
\nc{\Ab}{\mathcal{A}^{\mathrm{big}}}
\nc{\Sb}{\mathcal{S}^{\mathrm{big}}}
\nc{\As}{\mathcal{A}}
\nc{\Ss}{\mathcal{S}}
\nc{\ntens}{\widetilde{\otimes}}
\nc{\hR}{\widehat{R}}
\nc{\nconv}{\mathop{\mbox{\large $\odot$}}}
\nc{\snconv}{\mbox{\scriptsize$\odot$}}
\nc{\ts}{\tilde{s}}
\nc{\sho}{\mathcal{O}}
\nc{\bc}{\begin{cases}}
\nc{\ec}{\end{cases}}
\nc{\slnh}{{\widehat{\mathfrak{sl}}_N}}
\nc{\UA}{U_q'(\slnh)}
\nc{\cQ}{\mathcal{Q}}
\nc{\Irr}{\mathcal{I}rr}
\nc{\Irrw}[1]{{\Irr(\cC_{#1})}}
\nc{\Irrtw}[1]{{\Irr(\tcC_{#1})}}

\nc{\tQ}{\widetilde{\cQ}}
\nc{\bs}{\mathbf{s}}
\nc{\bL}{\mathbb{L}}
\nc{\tg}{\tilde{g}}

\nc{\conv}{\mathbin{\mbox{\large $\circ$}}}
\nc{\shconv}{\mathbin{\large\diamond}}
\nc{\sconv}{\mathbin{\large\Delta}}
\nc{\stens}{\mathbin{\large\Delta}}
\nc{\hconv}{\mathbin{\nabla}}
\nc{\htens}{\mathbin{\nabla}}

\nc{\Rm}{R^{\mathrm{ren}}}

\nc{\bQ}{\ol{Q}}

\nc{\de}{\on{\textfrak{d}\ms{1mu}}}

\nc{\xmono}{\ar@{>->}}
\nc{\xepi}{\ar@{->>}}
\nc{\db}[1]{\raisebox{-.5ex}[2ex][1.8ex]{$#1$}}
\nc{\wb}[1]{\mbox{$\rule[-1.1ex]{0ex}{2ex}#1$}}
\nc{\univ}{\mathrm{univ}}
\nc{\rM}{{}^*\mspace{-2mu}M}
\nc{\lM}{M^*}
\nc{\uqm}{\uq\smod}
\nc{\tR}{\widetilde{R}_{\gamma,\beta}}
\nc{\tx}{\tilde{x}}
\nc{\bi}{\mathbf{i}}
\nc{\ttau}{\widetilde{\tau}}

\nc{\tEnd}{\on{\widetilde{E}nd}}
\nc{\tHom}{\on{\widetilde{H}om}}

\nc{\K}{{J}}
\nc{\Kex}{{\K}_{\mathrm{ex}}}
\nc{\Kfr}{{\K}_{\mathrm{f\mspace{.01mu}r}}}
\nc{\coro}{\cor}
\nc{\tB}{\widetilde{B}}
\nc{\seed}{\mathscr{S}}

\nc{\up}{\mathrm{up}}
\nc{\bfa}{\mathbf{a}}
\nc{\bfb}{\mathbf{b}}
\nc{\bfc}{\mathbf{c}}
\nc{\bfm}{\mathbf{m}}
\nc{\hbfm}{ \widehat{\mathbf{m}}}

\newcommand{\Aqn}{A_q(\mathfrak n)}

\newcommand{\hA}{\widehat{\mathcal{A}}}

\newlength{\mylength}
\nc{\ov}[1]{\overline{#1}}
\nc{\Wlmj}[3]{\W_{#2,#3}^{(#1)}}
\nc{\Mkl}[2]{\M_\ttww(#1,#2)}
\nc{\mqs}{(-q^2)}
\nc{\Cquiver}{\upsigma}

\nc{\Kt}{\mathcal K_t}
\nc{\KT}{\mathbb{K}_t}
\nc{\yim}{y_{i,m}}
\nc{\yjm}{y_{j,m}}
\nc{\yjp}{y_{j,p}}
\nc{\yimp}{y_{i,m+1}}
\nc{\yjmp}{y_{j,m+1}}

\nc{\Refl}{\mathscr{S}}
\nc{\Reflinv}{{\Refl}^{-1}}

\nc{\Refn}{\mathsf{S}}

\nc{\catC}{\mathscr C}
\nc{\catA}{\mathcal A}
\nc{\shift}{{\mathrm T}}
\nc{\rE}{ \mathsf{E} }
\nc{\rW}{ \mathcal{W} }
\nc{\rES}{ \mathcal{E} }

\nc{\brd}{\sigma} 
\nc{\into}{\xymatrix@C=3ex{{}\ar@{^{(}->}[r]&{}}}
\nc{\dual}{\D\ms{0.2mu}}
\nc{\cdual}{\mathcal{D}}
\nc{\cat}[1][{\g}]{\catC_{#1}^0}

\nc{\catCO}{{\catC_\g^0}}
\nc{\catCOD}{{\catC_\g^{0, \ddD}}}
\nc{\catCQ}{{\catC_{\qQ}}}
\nc{\catCQd}{{\catC_{\widetilde{\qQ}}}}
\nc{\catCD}{{\catC_{\ddD}}}
\nc{\catCDK}{{\catC_{\ddD, \iK}}}
\nc{\Li}{{\La^\infty}}
\nc{\sig}{{\sigma(\g)}}
\nc{\sigZ}{{\sigma_0(\g)}}
\nc{\sigQ}{{\sigma_\qQ(\g)}}
\nc{\sigQd}{{\sigma_{\widetilde{\qQ}}(\g)}}
\nc{\phiQd}{\phi_{\widetilde{\qQ}}}
\nc{\sigD}{{\sigma_\ddD(\g)}}

\nc{\ZZ}{{\mathbf{Z}}}
\nc{\sP}{{\mathsf{P}}}
\nc{\sV}{{\mathsf{V}}}
\nc{\rxw}{{\underline{w}}}
\nc{\rxwz}{{\underline{w_0}}}
\nc{\boten}[1]{\overrightarrow{\bigotimes_{#1}}}
\nc{\cmA}{{\mathsf{A}}}
\nc{\cmC}{{\mathsf{C}}}
\nc{\ddD}{{\mathcal{D}}}
\nc{\ddDQ}{{\ddD_Q}}
\nc{\ddDQd}{{\ddD_{\widetilde{Q}}}}
\nc{\qQ}{{\mathcal{Q}}}
\nc{\gf}{{\g_{\mathrm{fin}}}}
\nc{\Df}{{\Delta_{\mathrm{fin}}}}
\nc{\If}{{I_{\mathrm{fin}}}}
\nc{\cmAf}{{\cmA_{\mathrm{fin}}}}
\nc{\weyl}{{\mathsf{W}}}
\nc{\weylf}{{\mathsf{W}_{\mathrm{fin}}}}
\nc{\sg}{{\mathfrak{S}}}
\nc{\weylA}{{\mathsf{W}_\cmA}}
\nc{\weylC}{{\mathsf{W}_\cmC}}
\nc{\Deg}{\mathrm{Deg}}
\nc{\Di}{\Deg^\infty}
\nc{\KRc}{{K_{q=1}(R_\cmC\gmod)}}
\nc{\prD}{{\Delta^+}}
\nc{\prDf}{{\Delta^+_{\mathrm{fin}}}}
\nc{\nrD}{{\Delta^-}}
\nc{\prDA}{{\Delta^+_\cmA}}
\nc{\prDC}{{\Delta^+_\cmC}}
\nc{\nrDC}{{\Delta^-_\cmC}}
\nc{\n}{{\mathfrak{n}}}
\nc{\Rt}{\mathsf{L}} 
\nc{\Cp}{\mathsf{C}} 
\nc{\cuspS}{{\mathscr{L}}}
\nc{\st}[1]{\{{#1}\}}
\nc{\bst}[1]{\bigl\{{#1}\bigr\}}
\nc{\WS}{Quantum affine Schur--Weyl duality\xspace}
\nc{\zz}{{{\mathbf{z}}}}
\nc{\wlP}{\mathsf{P}}
\nc{\twlP}{\widetilde{\wlP}}
\nc{\wl}{\wlP}
\nc{\clp}{{\mathrm{cl}}}
\nc{\wlPc}{{\wlP_\clp}}
\nc{\awlP}{\widehat{\mathsf{P}}}
\nc{\dM}{\mathsf{M}}
\nc{\dC}{\mathsf{C}}
\nc{\cC}{\mathcal{C}}
\nc{\tcC}{\widetilde{\cC}}

\nc{\lR}{\widetilde{{R}}}
\nc{\zero}{\mathrm{zero}}
\nc{\PD}{principal }

\nc{\prtl}[1][J]{\rootl_{#1}^+}
\nc{\hL}{\widehat{\Rt}}
\nc{\hF}{\widehat{\F}}

\nc{\Aff}{\mathrm{Aff}}
\nc{\rT}{\mathcal{T}}		
\nc{\rr}{rationally renormalizable\xspace}
\nc{\RA}{{R_\cmA}}		
\nc{\RC}{{R_\cmC}}		
\nc{\proolim}[1][]{\mathop{``{\varprojlim}{\mbox{''}}}\limits_{#1}}
\nc{\qtq}[1][\text{and}]{\quad\text{#1}\quad}
\newcommand{\gW}{\mathrm{W}}
\newcommand{\sgW}{\mathrm{W}^*}
\newcommand{\Sp}{\mathrm{span}_{\mathbb{R}_{\ge0}}}  	

\newcommand{\lT}{\mathcal{T}} 		
\newcommand{\lTs}{\mathcal{T}^{\ms{1mu}*}}
\nc{\corh}{\widehat{\cor}}
\nc{\ang}[1]{\langle{#1}\rangle}
\nc{\rc}{renormalizing coefficient\xspace}
\nc{\cz}{{\cor[z^{\pm1}]}}
\nc{\tp}{\ms{1.5mu}{\widetilde{p}}\ms{2mu}}
\nc{\G}{\mathcal{G}}
\nc{\cc}{\mathfrak{c}}
\nc{\rsP}{{\Phi_\g}}
\nc{\rsX}{{X_\g}}
\nc{\rs}{ \mathsf{s} }
\nc{\Dynkin}{\mathsf{D}}
\nc{\Dat}{\sigma}
\nc{\hf}{\xi}

\nc{\cd}{\mathrm{D}}
\nc{\cm}{\mathbf{m}}
\nc{\MS}{\mathsf{MS}}
\nc{\hMS}{\widehat{\mathsf{MS}}}
\nc{\rS}{\mathbf{S}}
\nc{\rSs}{\rS^*}
\nc{\brS}{\overline{\rS}}

\nc{\crI}{ \mathsf{I}}
\nc{\crhI}{\mathscr{S} }
\nc{\crD}{\mathsf{D}}
\nc{\crc}{\mathsf{c}}
\nc{\crB}{\mathscr{P}_n}
\nc{\cru}{\mathsf{u}}
\nc{\crsh}{\mathsf{sh}}
\nc{\bfs}{\mathbf{s}}
\nc{\crBB}[1]{\mathscr{P}_#1}

\nc{\gc}{{\g_{\cmC}}}
\nc{\cL}{\mathcal{L}}
\nc{\cLD}{{\cL_\ddD}}
\nc{\ccL}{\mathscr{L}}
\nc{\ccLD}{{\ccL_\ddD}}
\nc{\clen}{\mathsf{len}}

\nc{\hv}{\mathsf{1}}
\nc{\qt}[1]{\quad\text{#1}}
\nc{\snoi}{\smallskip\noindent}
\nc{\mnoi}{\medskip\noindent}

\nc{\ul}[1]{\underline{#1}}
\nc{\dul}[1]{\underline{\underline{#1}}}
\nc{\qh}{\qedhere}
\nc{\ca}{\mathsf{v}}
\nc{\sck}[1][k]{\ms{8mu}{\raisebox{-1.3ex}{$\scriptstyle{#1}$}\hs{-1.4ex}{\succ}}\ms{4mu}}
\nc{\scke}[1][k]{\ms{8mu}{\raisebox{-1.3ex}{$\scriptstyle{#1}$}\hs{-1.4ex}%
{\succcurlyeq}}\ms{4mu}}

\nc{\edot}{\emptyset}
\nc{\Nf}{N_{\gf}}
\nc{\fin}{\mathrm{fin}}
\nc{\trg}{\scalebox{.7}{$\triangle$}}
\nc{\ake}[1][1ex]{\rule[-#1]{0ex}{1ex}}
\nc{\akew}[1][1ex]{\rule[-1ex]{#1}{0ex}}
\nc{\akeu}[1][1ex]{\rule[#1]{0ex}{1ex}}

\nc{\bg}{\mathscr{B}}
\newcommand{\B}{B(\infty)}
\nc{\ipi}{{{_i}\pi}}
\nc{\pii}{{\pi_i}}
\nc{\Ld}{\mathcal{P}}
\nc{\Dc}{\Upsilon}
\nc{\bR}{\mathsf{k}}
\nc{\bRs}{\bR^*}
\nc{\rR}{\mathrm{R}}
\nc{\bfi}{\mathbf{i}}
\nc{\bfj}{\mathbf{j}}
\nc{\bfk}{\mathbf{k}}
\nc{\vt}{\vartheta}
\nc{\nP}{\Upsilon}
\nc{\hP}{\widehat{\nP}}
\nc{\fF}{ \widetilde{\mathrm{F}}}
\nc{\fE}{ \widetilde{\mathrm{E}}}
\nc{\fR}{\mathscr{R}}
\nc{\fT}{\mathscr{T}}
\nc{\orb}{\mathrm{orb}}
\nc{\fdr}{\mathrm{r}}

\nc{\Gu}{\mathrm{G}^{\rm up}}
\nc{\ve}{\varepsilon}
\nc{\sS}{\mathcal{S}}
\nc{\fr}{\mathrm{fr}}
\nc{\uf}{\mathrm{ex}}

\nc{\gR}{\mathrm{g}^{\mathrm{R}}}
\nc{\gL}{\mathrm{g}^{\mathrm{L}}}
\nc{\PBW}{\mathrm{PBW}}
\nc{\STR}{\mathrm{STR}}
\nc{\eqfr}{\equiv_{\mathrm{fr}}}
\nc{\Irrp}{\Irr_{\rm sd}}
\nc{\cI}{\mathcal{I}}
\nc{\Br}{\mathsf{Br}}
\nc{\Seq}{\mathrm{Seq}}
\nc{\mD}{\mathsf{D}}
\nc{\LB}{\widetilde{B}}
\nc{\fz}{\mathsf{c}}
\nc{\qmC}{\mathrm{C}} 
\nc{\eqq}{{\, \equiv_q \,}}
\nc{\inq}{{\, \in_q \,}}
\nc{\Zq}{{\Z[q^{\pm 1}]}}
\nc{\tA}{\widetilde{A}}
\nc{\fup}{f^{\up}}
\nc{\fupi}[1]{{\fup_{#1}}}
\nc{\inv}{{\mathsf{inv}}}
\nc{\cP}{\mathsf{P}}
\nc{\cS}{\mathsf{S}}
\nc{\tcP}{\widetilde{\cP}}
\nc{\tcS}{\widetilde{\cS}}
\nc{\strS}{\mathfrak{S}}
\nc{\pbwP}{\mathfrak{P}}
\nc{\triA}{\mathscr{A}}
\nc{\triAs}{\mathscr{A}_{\clubsuit}}
\nc{\triB}{\mathscr{B}}
\nc{\triBs}{\mathscr{B}_{\clubsuit}}
\nc{\triC}{\mathscr{C}}
\nc{\triCs}{\mathscr{C}_{\clubsuit}}
\nc{\triD}{\mathscr{D}}
\nc{\triDs}{\mathscr{D}_{\clubsuit}}
\nc{\triE}{\mathscr{E}}
\nc{\triF}{\mathscr{F}}
\nc{\triG}{\mathscr{G}}
\nc{\triX}{\mathscr{X}}
\nc{\emb}{\mathcal{E}}
\nc{\rv}{\mathsf{F}}
\nc{\im}{\imath}
\nc{\jm}{\jmath}
\nc{\qm}{{q_{\mathrm{sh}}}}
\nc{\ord}{\mathrm{ord}}
\nc{\cl}{\mathrm{cl}}
\nc{\hIQ}{\widehat{I}_\qQ}
\nc{\oIQ}{\overline{I}_\qQ}
\nc{\qmD}{\mathsf{D}}
\nc{\rev}{\mathrm{rev}}
\nc{\KR}{{Kirillov-Reshetikhin}}
\nc{\stb}{{stable} }

\newcommand{\bnom}[1]{\begin{bmatrix}#1\end{bmatrix}}
\nc{\Qh}{\Q(q^{1/2})}
\nc{\bD}{\overline{\mathcal{D}}}
\newcommand{\trmG}{\overline{\mathsf{G}}}

\makeatletter
\newcommand{\namedlabel}[2]{%
  #2%
  \def\@currentlabel{#2}%
  \label{#1}%
}
\makeatother

\newcommand{\cnum}[1]{%
\tikz[baseline=(char.base)]{
\node[
draw,
circle,
inner sep=0.2pt,
minimum size=1.1em,
font=\scriptsize
] (char) {$#1$};
}%
}

\title[Imaginary elements and crystals]
{Imaginary elements and crystals for quantum groups}

\author[K.-H.\ Lee]{Kyu-Hwan Lee$^{\star}$}
\thanks{$^{\star}$K.-H.\ Lee was supported by NSF grant DMS-2617278.}
\address[K.-H.\ Lee]{Department of Mathematics, University of Connecticut, Storrs, CT 06269, USA; \& Korea Institute for Advanced Study, Seoul     02455, Republic of Korea}
\email[K.-H.\ Lee]{khlee@math.uconn.edu}

\author[E.\ Park]{Euiyong Park$^{\dagger}$}
\thanks{$^{\dagger}$The research of E.\ Park was supported by the National Research Foundation of Korea (NRF) Grant funded by the Korea Government(MSIT)(RS-2023-00273425 and NRF-2020R1A5A1016126).}
\address[E.\ Park]{Department of Mathematics, University of Seoul, Seoul 02504, Republic of Korea;
\&
Department of Mathematics, University of Connecticut, Storrs, CT 06269, USA
}
\email[E.\ Park]{epark@uos.ac.kr}

\makeatletter
\@namedef{subjclassname@2020}{\textup{2020} Mathematics Subject Classification}
\makeatother
\date{November 19, 2026}

\begin{document}
\maketitle

\begin{abstract}
In this paper, we introduce the notion of \stb imaginary elements in quantum groups and use crystal basis theory to construct explicit one-parameter families of such elements for all finite types. We show that all of Leclerc's examples of imaginary elements are \stb imaginary. Using the generalized quantum Schur--Weyl duality, we also construct explicit examples of corresponding imaginary modules over quantum affine algebras of untwisted affine types.
\end{abstract}

\tableofcontents

\section{Introduction}

Let $U_q(\g)$ be the quantum group over the base field $\Q(q)$ associated with a generalized symmetrizable Cartan matrix $\cmC$, where $q$ is an indeterminate.  
For each element $w$ of the Weyl group $\weyl$ associated with $\cmC$, we denote by $A_q(\n(w))$ the \emph{quantum unipotent coordinate ring}, which is the subalgebra of the dual of the negative half $U_q^-(\g)$ generated by \emph{dual PBW vectors} associated with $w$. 
This algebra $A_q(\n(w))$ can be understood as a quantum analogue of the coordinate ring  $\C[N(w)]$ of the {unipotent subgroup}. The algebra $A_q(\n(w))$ is deeply related to several research areas including geometry, cluster algebras and categorification (see \cite{GLS13, Kimura12, KKKO18} and references therein).

The quantum unipotent coordinate ring $A_q(\n(w))$ has a quantum cluster algebra structure (\cite{GLS13, GY17}), in which each reduced expression of $w$ gives an initial seed consisting of the \emph{quantum unipotent minors}. 
It was shown in \cite{KKKO18} that the cluster algebra structure of $A_q(\n(w))$ can be interpreted from the monoidal categorical viewpoint via the categorification of $A_q(\n(w))$ using \emph{quiver Hecke algebras} (\cite{KL09,KL11,R08}). This is called the \emph{monoidal categorification} of $A_q(\n(w))$. Via the categorification, the initial seeds arising from quantum unipotent minors correspond to the monoidal seeds arising from \emph{determinantial modules} (see Section \ref{Sec: QHA} and \ref{Sec: Cw}). 
It was shown in \cite{KKKO18, Qin17, Qin20} that every cluster monomial of $A_q(\n(w))$ is contained in the distinguished basis $\Gup(w)$ of $A_q(\n(w))$, called the \emph{dual canonical basis} (\cite{L90, L91, LusztigBook}) or the \emph{upper global basis} (\cite{K91, K93}).
It is an interesting and important question to ask how cluster monomials can be described in the basis $\Gup(w)$ as a subset, which yields the \emph{reachability conjecture}, i.e.,  the set of cluster monomials of $A_q(\n(w))$ coincides with the set of \emph{real} elements of $\Gup(w)$ (see \cite[Conjecture 5.2]{HL16} and \cite[Conjecture 5.7]{Qin20A}). Here, an element $x\in \Gup(w)$ is \emph{real} if $q^A x^2$ is contained in $\Gup(w)$ for some integer $A\in \Z$, and $x$ is \emph{imaginary} otherwise. Since all cluster monomials are real, the point of the conjecture is to show that real elements are cluster monomials. 
Thus it is important to understand which elements of $\Gup(w)$ are real and which are imaginary.
The conjecture has been partially confirmed for cluster algebras of finite or affine type (\cite{Saka25}), but it remains widely open in full generality as far as the authors know.

In \cite{Lec03}, Leclerc gave examples of imaginary elements for types $A_5$, $B_3$, $C_3$, $D_4$ and $G_2$ to disprove a conjecture of Berenstein and Zelevinsky on $q$-commuting product of $\Gup(w)$. In the examples in \cite[Section 2]{Lec03}, the imaginary elements and their squares are written explicitly in terms of \emph{PBW data}. 
On the categorical side, \emph{imaginary modules}, which correspond to imaginary elements via monoidal categorification, have also been studied in the module categories of quantum affine algebras and quiver Hecke algebras.
For quantum affine algebras of type $A_n^{(1)}$, Lapid and Mínguez \cite{LM18} provided a criterion for determining whether a \emph{regular} simple module is real or imaginary, while Brito and Chari \cite{BC23} constructed imaginary modules by studying the tensor product of a \emph{higher-order Kirillov–Reshetikhin module} with its dual.
This construction was later extended to the setting of \emph{prime snake modules} by Brito and Moura (\cite{BM26}).   
In \cite{EL24}, Early and Li studied limit g-vectors of cluster algebras and constructed explicit examples of imaginary modules for quantum affine algebras of untwisted types $A$ and $D$. Further examples of imaginary modules in this direction were later given by Li in \cite[Section 6]{Li26}. 
In the case of monoidal module subcategories $\catC$ over quantum affine algebras and quiver Hecke algebras, under the condition that the Grothendieck ring $K(\catC)$ is a cluster algebra of \emph{affine type}, Sakamoto gave a complete classification of real and imaginary modules in $\catC$ using the combinatorics of affine type cluster algebras (\cite{Saka25}).  

In this paper, we introduce the notion of \emph{\stb imaginary elements} and construct one-parameter families of \stb imaginary elements in the basis $\Gup(w)$ for each \emph{finite type} (see Proposition \ref{Prop: triX in ft}), which yields one-parameter families of imaginary modules over quantum affine algebras via the \emph{generalized quantum Schur--Weyl duality} introduced in \cite{KKK18} (see Proposition \ref{Prop: imag for Cg}).
In the case of quantum affine algebras of untwisted affine types $A^{(1)}_{6}$, $B^{(1)}_{3}$, $C^{(1)}_{4}$, $D^{(1)}_{5}$, $E^{(1)}_{k}$ ($k=6,7,8$), $F^{(1)}_{4}$ and $G^{(1)}_{2}$,
explicit examples of such imaginary families are given in Section \ref{Sec: im for qaa}. It turns out that all of Leclerc's imaginary elements are \stb imaginary (see Section \ref{Sec: LE}), which suggests that \stb imaginary elements can be viewed as a generalization of Leclerc's examples. Our one-parameter families are generally distinct from previously known families, although some known imaginary modules appear in our families in low ranks (see, for example, Remarks \ref{rmk: LM} and \ref{rmk: Sakamoto}).

One of the key ingredients in the proof of our main result is the \emph{quantum twist automorphism} introduced by Kimura and Oya \cite{KiOy21}. This is an automorphism $\eta_w$ of the localized algebra $\tA_q(\n(w))$ (see Section \ref{Sec: Loc}) that quantizes the \emph{twist automorphism} introduced by Berenstein, Fomin, and Zelevinsky \cite{BFZ96, BZ97}.
The automorphism $\eta_w$ admits a categorification via the \emph{right dual functor} $\dual_w$ on the localized category $\tcC_w$. Here $\tcC_w$ is obtained by localizing the category $\cC_w$ of finite-dimensional modules over \emph{quiver Hecke algebras}, which categorifies $A_q(\n(w))$.
In \cite{JP25}, Jung and the second named author provided a crystal-theoretic description $\cdual_w$ of $\eta_w^{-1}$, which enables us to compute its action in a combinatorial way via crystal realizations. 
Let $B(w)$ be the subset of the \emph{infinite crystal} $B(\infty)$ corresponding to $\Gup(w)$ (Section \ref{Sec: B(w)}) and let $\Irr(\cC_w)$ be the set of the isomorphism classes of simple modules in $\cC_w$ (Section \ref{Sec: Cw}).
We define a \emph{$\cdual_w$-\stb imaginary} element $b \in B(w)$ to be 
a \emph{non-frozen} element that is fixed by $\cdual_w$ up to a multiplication by frozen variables (see Definition \ref{Def: frozen} and \ref{def: Dw PI}). Thanks to Lemma \ref{Lem: fixed by D}, all $\cdual_w$-\stb imaginary elements are imaginary, but the converse does not hold in general (see Remark \ref{rmk: pi and i}). 

Moreover, when $U_q(\g)$ is of finite type, Proposition \ref{Prop: bk pi} shows that if $b$ is $\cdual_w$-\stb imaginary, then so is every \emph{multiple} $b^{[k]}$ of $b$ for $k\in\Z_{>0}$ (see Definition \ref{def: k PBW}). This property is a key advantage of considering \stb imaginary elements, as it allows us to construct entire families from a single example.
By the equivalence in Lemma \ref{lem: equiv b}, 
whenever a \stb imaginary element $b \in B(w)$ is given, we have the one-parameter family $\{ \Gup(b^{[k]}) \}_{k\in \Z}$ (resp.\ 
$\{ \cL(b^{[k]}) \}_{k\in \Z}$) of imaginary elements (resp.\ imaginary modules), where $\Gup(b^{[k]})$ (resp.\ $\cL(b^{[k]})$) is the element of $\Gup(w)$ (resp.\ $\Irr(\cC_w)$) corresponding to $b^{[k]}$.

The \emph{crystal basis theory} (\cite{K95, KashBook02}) plays a crucial role to find \stb imaginary elements. Using module-theoretic arguments based on the \emph{integer-valued invariants} $\La$ and $\de$ for quiver Hecke algebras (Section \ref{Sec: iv La de}), we give a condition to find a \stb imaginary  element $b$ in terms of the crystal theory (see Theorem \ref{Thm: main} and Corollary \ref{cor: condition for PI}). In particular, the conditions in Corollary \ref{cor: condition for PI} is written solely in terms of weights and Weyl group elements, which is compatible with the \emph{folding} and \emph{embedding} of crystals (Section \ref{Sec: similarity}). 

By applying Corollary \ref{cor: condition for PI} with the aid of SageMath \cite{SageMath}, we construct \stb imaginary elements $b_\triX$ for all finite types (Proposition \ref{Prop: triX in ft}), where $\triX$ is a triple satisfying the conditions of Corollary \ref{cor: condition for PI}. As noted above, for each such element $b_\triX$, we obtain a one-parameter family of \stb imaginary elements $\{b_\triX^{[k]}\mid k\in\Z_{>0}\}$.
We also show that all of Leclerc's examples arise as instances of these \stb imaginary elements (see Section \ref{Sec: LE}).
We remark that there exists a \stb imaginary element which does not satisfy the conditions of Theorem \ref{Thm: main} and Corollary \ref{cor: condition for PI} (see Remark \ref{rmk: stb im G2}).
In addition, since the generalized quantum Schur--Weyl duality $\F_\ddD$ arising from a \emph{strong duality datum} $\ddD$ is an exact monoidal functor preserving simple modules, the images of $\cL(b_\triX^{[k]})$ ($k\in \Z_{>0}$) under the functor $\F_\ddD$ are imaginary modules over quantum affine algebras (see Proposition \ref{Prop: imag for Cg}). Using the duality datum arising from a $Q$-datum, together with computations in SageMath \cite{SageMath}, we construct explicit examples of such imaginary families in terms of dominant monomials for affine types $A^{(1)}_{6}$, $B^{(1)}_{3}$, $C^{(1)}_{4}$, $D^{(1)}_{5}$, $E^{(1)}_{k}$ ($k=6,7,8$), $F^{(1)}_{4}$ and $G^{(1)}_{2}$ (Section \ref{Sec: ex for qaa}).

This paper is organized as follows. 
In Section \ref{Sec: prelim}, we recall briefly necessary background on quantum groups, crystals and quiver Hecke algebras.
In Section \ref{Sec: QTA}, we review the crystal-theoretic description $\cdual_w$ of the quantum twist automorphism given in \cite{JP25}. 
In Section \ref{Sec: PI}, we introduce the notion of \stb imaginary modules and show their properties, and in Section \ref{Sec: crystal and pi} we give a condition to find \stb imaginary elements in terms of the crystal theory. Section \ref{Sec: PI for finite types} constructs explicit \stb imaginary elements $b_\triX$ for all finite types, and Section \ref{Sec: LE} shows that all of Leclerc's imaginary elements appear as examples of $b_\triX$. In Section \ref{Sec: im for qaa}, we construct imaginary modules arising from $b_\triX$ via the generalized quantum Schur--Weyl duality and gives examples of such imaginary modules for untwist affine types.  
In Appendix \ref{Sec: appA}, we establish a necessary property of the product involving $\eta_w$ that is needed for the proof of Lemma \ref{Lem: fixed by D} in the non-symmetric case.


\vskip 1em

{\bf Acknowledgments}.\ 
The authors would like to thank M. Kashiwara for valuable comments and suggestions on \stb imaginary modules and their properties (Remark \ref{rmk: new def of pi}). The authors also thank M. Kim and S.-j. Oh for fruitful discussions and helpful computations, and are grateful to H.~Oya for many helpful discussions and valuable comments.
E.P.\ would like to thank the Department of Mathematics at the University of Connecticut for its hospitality and excellent research environment during his visit.

\vskip 1em 

\textbf{Convention}. The following convention will be used in the paper.
\bna
\item  For a set $A$, we write $|A|$ for the cardinality of $A$.

\item For $a,b \in \mathbb{Z} \sqcup \{\pm\infty\}$ with $a \le b$, we define the \emph{interval} $[a,b]$ by
\[
[a,b] := \{ t \in \mathbb{Z} \mid a \le t \le b \}.
\]
If $a>b$, we set $[a,b]=\emptyset$. We sometimes write $[a]$ for $[a,a]$ if no confusion arises.
 \item For a monoid $A$ and an index set $K$, we denote by $A^{\oplus K}$ the direct sum of copies of $A$ indexed by $K$.
If no confusion arises, then we simply write $A^K$ for $A^{\oplus K}$.

\ee

\vskip 2em

\section{Preliminaries} \label{Sec: prelim}

\subsection{Quantum groups}
In this subsection, we briefly recall the definition of quantum groups and their properties. We refer the reader to \cite{LusztigBook} for details. 

\subsubsection{Cartan data}
A {\it Cartan datum}  $ ( \cmC, \wlP,\Pi,\wlP^\vee,\Pi^\vee) $ consists of a generalized (symmetrizable) \emph{Cartan matrix} $ \cmC =(c_{ij})_{i,j\in I}$,
 a free abelian group $\wlP$,  $\wlP^{\vee} := \Hom_{\Z}( \wlP, \Z )$, the set $\Pi := \{ \alpha_i \in \wlP \mid i\in I \}$ of \emph{simple roots} and the set $\Pi^{\vee} := \{ h_i \in \wlP^\vee \mid i\in I\}$ of \emph{simple coroots} that satisfy the following conditions: 
 \bna
\item $\lan h_i, \alpha_j \ran = c_{ij}$ for $i,j \in I$,
\item $\Pi$ is linearly independent over $\Q$,
\item for each $i\in I$, there is $\Lambda_i \in \wlP$ such that $\lan h_j,\Lambda_i \ran =\delta_{j,i}$ for all $j \in I$.
\item there is a symmetric bilinear 
form $( \cdot \, , \cdot )$ on $\wlP$ such that 
$(\al_i,\al_i)\in 2\Z_{>0}$ and 
$ \lan h_i,  \lambda\ran = {2 (\alpha_i,\lambda)}/{(\alpha_i,\alpha_i)}$. 
\ee
The abelian group $\wlP$ (resp.\  $\wlP^\vee$) is called the \emph{weight lattice} (resp.\ \emph{dual weight lattice}), and the element $\La_i \in \wlP$ is called the $i$th fundamental weight. We write $\wlP_+ := \{ \La \in \wlP \mid \ang{h_i, \La} \in \Z_{\ge 0} \text{ for all } i\in I  \}$, and call elements of $\wlP_+$ \emph{dominant integral weights}. We write $\rlQ:= \sum_{i\in I} \Z \al_i$, $\rlQ_+:= \sum_{i\in I} \Z_{\ge0} \al_i$ and denote by
$\prD$ (resp.\ $\nrD$) the set of all positive (resp.\ negative) roots of $\cmC$.  
For $\beta = \sum_{i\in I} c_i \al_i \in \rlQ_+$, we write $\Ht(\beta) := \sum_{i\in I} c_i \in \Z_{\ge0}$.
Throughout the paper, we take a bilinear form $( \cdot \, , \cdot )$ such that $ \min \{ (\al_i, \al_i)  \mid i\in I\}=2$.

Let $\weyl = \ang{s_i \mid i\in I}$ be the \emph{Weyl group} associated with $\cmC$,  where $s_i$ are the \emph{Coxeter generators} for $i\in I$. An expression $\underline{w} = s_{i_1} s_{i_2} \cdots s_{i_r} $ of an element $w \in \weyl$ is \emph{reduced} if $\underline{w}$ is minimal among all expressions of $w$. In this case, we write $\ell(w) := r$ and call $\bfi = (i_1, i_2, \ldots, i_r)$ a \emph{reduced sequence} of $w$. We set 
$$
\supp(w):= \{ i_1, i_2, \ldots, i_r \},
$$ 
where $(i_1, i_2, \ldots, i_r) $ is a reduced sequence of $w$. 
We define $\Seq(w)$ to be the set of all reduced sequences of $w$, i.e.,  
\begin{align}
\Seq(w) := \{ \bfi = (i_1, i_2, \ldots, i_{\ell(w)}) \in I^{\ell(w)} \mid  w = s_{i_1} s_{i_2} \cdots s_{i_{\ell(w)}} \} 
\end{align}
We write $v < w$ in the \emph{Bruhat order} if there exists a reduced sequence of $w$ containing a reduced sequence of $v$ as a subsequence.

\subsubsection{Quantum groups}
Let $q$ be an indeterminate and set $q_i := q^{(\al_i, \al_i)/2}$ for $i\in I$. For $i\in I$ and $n  \in \Z_{\ge0}$, define
$$ 
[n]_{q_i}= \dfrac{q_i^n-q_i^{-n}}{q_i-q_i^{-1}} \quad \text{ and }\quad [n]_{q_i}! = \prod_{k=1}^n [k]_{q_i}.
$$
The quantum group $U_q(\g)$ associated with $ (\cmC,\wlP,\Pi,\wlP^\vee,\Pi^\vee) $ is the associative $\Q(q)$-algebra generated by $f_i$, $e_i$ and $q^h$ for $i\in I$ and $ h \in \wlP^\vee$ satisfying certain defining relations (see \cite{LusztigBook} for details).
We write $U_{\Zq}^-(\g)$ for the $\Zq$-subalgebra of
$U_q(\g)$ generated by $f_i^{(n)} \seteq f_i^n/[n]_{q_i}!$ for $i \in I$ and $n\in \Z_{>0}$.  
For an element $x$ of the weight space $U_q^-(\g)_\beta$, 
we write $\wt(x):=\beta \in \rlQ_- := - \rlQ_+$. 
The \emph{quantum unipotent coordinate ring} $\Aqn$ is defined by 
$$
A_q(\n) \seteq \bigoplus_{\beta \in \rlQ_-} A_q(\n)_\beta, 
$$
where $A_q(\n)_\beta \seteq \Hom_{\Q(q)}(U_q^-(\g)_\beta,\Q(q))$ (see \cite[Section 8]{KKKO18} for details).
For $x,y \in A_q(\n)$ and  a subset $G \subset A_q(\n)$, we write 
$$
x \eqq y \qtq x \inq G
$$  
if  $x = q^k y$ and $q^t x \in G$ for some  $k,t\in \Z$.
We denote by $\Gup(\infty)$ the \emph{dual canonical basis} or the \emph{upper global basis} of $A_q(\n)$ (see \cite{K91,K93, L90, LusztigBook} and references therein for precise definitions, and see also \cite[Section 8]{KKKO18}). 
An element $x \in \Gup(\infty)$ is called \emph{real} if $x^2 \inq \Gup(\infty) $. Otherwise, $x$ is called \emph{imaginary.}

\subsubsection{PBW bases} \label{Sec: PBW}
For $i\in I$, Lusztig's \emph{braid symmetries} are written as 
\begin{align} \label{Eq: braid}
\lT_i := T_{i,-1}' \quad \text{ and } \quad \lTs_i := T_{i,1}'',
\end{align}
where $T_{i,-1}'$ and $ T_{i,1}''$ are defined in \cite[Chapter 37.1]{LusztigBook}. 
For  $\bfi = (i_1, i_2, \ldots, i_r) \in \Seq(w)$,  $k \in [1, r]$ and $c \in \Z_{ \ge 0}$,  we set  
\begin{equation} \label{Eq: betak}
\begin{aligned} 
\beta_k &:= s_{i_1} \cdots s_{i_{k-1}}(\al_{i_k}) \in \prD, \\
 \rv_{\bfi} (\beta_k) ^{(c)} &:= \lT_{i_1}\lT_{i_2} \cdots \lT_{i_{k-1}}(f_{i_k}^{(c)})  \in U_{\Zq}^-(\g),
\end{aligned}
\end{equation}
and, for any $\bfa=(a_k)_{k\in [1,r]} \in \Z_{\ge 0}^{\oplus r}$, define 
\begin{align} \label{Eq: PBW ep}
\rv_{\bfi} (\bfa) := 
\rv_{\bfi} (\beta_r) ^{(a_r)} \rv_{\bfi} (\beta_{r-1}) ^{(a_{r-1})} \cdots 	\rv_{\bfi} (\beta_{1}) ^{(a_{1})}.
\end{align}
The $\rv_{\bfi} (\beta_k)$ and  $ \rv_{\bfi} (\bfa)  $ are called the \emph{PBW vectors} and  \emph{PBW monomials} respectively with respect to $\bfi$ and $\bfa$.
We denote by $\rv^{\up}_{\bfi} (\beta_k)$ and $ \rv^{\up}_{\bfi} (\bfa) $ the dual of $\rv_{\bfi} (\beta_k)$ and $ \rv_{\bfi} (\bfa)  $ respectively with respect to the Kashiwara bilinear form (\cite[Proposition 3.4.4]{K91}). We denote by $A_q(\n(w))$ the subalgebra of $ A_q(\n)$ generated by all dual PBW monomials $ \rv^{\up}_{\bfi} (\bfa) $. 
These dual PBW monomials form a linear basis of $A_q(\n(w))$, which is called the \emph{dual PBW basis}. 
It was shown in \cite{Kimura12} that 
the algebra $A_q(\n(w))$ respects the basis $\Gup(\infty)$, i.e., 
$$
\Gup(w) := A_q(\n(w)) \cap \Gup(\infty) 
$$
forms a $\Q(q)$-linear basis of $A_q(\n(w))$.

\subsection{Crystals}
In this subsection, we briefly recall the notion of crystals and their properties. 
A \emph{crystal} is a set $B$ together with 
$\te_i$, $\tf_i \col  B \rightarrow B\,\sqcup\{0\}$,  
$\wt\col  B \rightarrow \wlP$, and $\varphi_i$,  $\ep_i \col  B \rightarrow \Z \,\sqcup \{ -\infty \}$	 
for all $i\in I$ satisfying certain axiomatic conditions \cite[Section 7.2]{K95}. 
We refer the reader to \cite{K95,KashBook02} for details.

\subsubsection{Infinite crystal $B(w)$} \label{Sec: B(w)}
We denote by $B(\infty)$ the \emph{infinite crystal} of $A_q(\n)$ (or $U_{q}^-(\g)$) with the highest weight vector $\one$. 
The crystal $B(\infty)$ is obtained from $\Gup(\infty)$ by specializing  at $q=0$, and this specialization gives a 1-1 correspondence between $B(\infty)$ and $\Gup(\infty)$. 
 For each $b \in B(\infty)$, we denote by $\Gup(b)$ the corresponding element of $\Gup(\infty)$ under the specializing $q=0$. Then we have 
$ \Gup(\infty) = \{ \Gup(b) \mid b\in B(\infty) \}.$
For $w\in \weyl$, we define 
\begin{align}\label{Eq: Bw}
B(w) := \{ b \in B(\infty) \mid \Gup(b) \in A_q(\n(w)) \}
\end{align}
so that 
\begin{align} \label{Eq: Gup B}
\Gup(w) = \{ \Gup(b) \mid b \in B(w) \}.
\end{align}

\subsubsection{Parameterizations of $B(w)$}
Let $\bfi = (i_1,i_2, \ldots, i_r) \in \Seq(w)$ and $b\in B(w)$. There exists $\bfa \in \Z_{\ge 0}^{\oplus r}$ such that  
\begin{align} \label{Eq: Gup PBW}
\Gup(b) \equiv \rv^{\up}_{\bfi}  (\bfa) \mod qL(\infty), 
\end{align}
where $L(\infty)$ is the $\Z[q]$-lattice of $A_q(\n(w))$ generated by $\Gup(w)$.
In this case we write  
\begin{align} \label{Eq: PBWi}
\PBW_{\bfi}(b) := \bfa 
\end{align}
(see \cite{LusztigBook} for details), which is called the \emph{PBW datum} of $b$ with respect to $\bfi$.
Note that this notation for PBW data is consistent with that used in \cite{KKOP18, P23, JP25}. 
On the other hand, we define 
\begin{align} \label{Eq: STRi}
\STR_{\bfi} (b) := (c_k)_{k\in [1,r]},
\end{align}
where $b_1 = b$,  $c_k = \ve_{i_k}(b_{k})$ and $b_{k+1} = \te^{c_k}_{i_k}(b_{k}) $ for $k\in [1,r]$. This is called the \emph{string datum} of $b$ with respect to $\bfi$.

\subsubsection{Saito crystal reflections}
Let $*: U_q(\g) \rightarrow U_q(\g)$ be the anti-automorphism over $\Q(q)$ given by 
$ e_i \mapsto e_i$, $ f_i \mapsto f_i$, and $q^h \mapsto q^{-h}$ for $i\in I$ and $h \in \wlP^\vee$. We define 
$$
\tes_i := * \circ \te_i \circ * \qtq \tfs_i := * \circ \tf_i \circ * \qquad \text{ for $i\in I$,}
$$
which gives another $*$-twisted crystal structure to $B(\infty)$, in which 
 the notations $\eps (b)$ and $\varphi_i^*(b)$ are used. 
For $b\in B(\infty)$,  write 
$$
\tem_i(b) := \te_i^{\ep_i(b)}(b)\qtq \tesm_i(b) := {\te_i}^{* \eps_i(b)}(b).
$$

We denote by $T_i$ the \emph{Saito crystal reflection} (\cite{Saito94}), i.e., 
for any $b\in \Binf$ with $\ep_i(b)=0$,
$$
T_i(b) = \tf_i^{ \varphi_i^*(b) } \te_i^{* \eps_i(b) }(b).
$$
The reflection $T_i$ is a crystal-theoretic counterpart of the braid symmetry $\T_i$. 
Note that $ \eps_i(T_i(b))=0$ and $T_i^{-1} = * \circ T_i \circ *$. 

\subsubsection{Highest weight crystals $B(\La)$} \label{Sec: HWC}
For a dominant integral weight $\La \in \wlP_+$,  $B(\La)$ is the crystals of the highest weight module $V(\La)$ of highest weight $\La$. We denote by $b_\La$ (resp.\ $u_\la$) the highest weight vector of $B(\La)$ (resp.\ $V(\La)$).
Using the Weyl group action on $B(\La)$ (\cite[Section 11]{K95}), for any $v = s_{j_1} s_{j_2} \cdots  s_{j_m} \in \weyl$, we define $ b_{\La, v\La} := v b_{\La} =  s_{j_1} s_{j_2} \cdots s_{j_m} (b_\La) $, which is called the \emph{extremal vector} of $B(\La)$ with extremal weight $v\La$. Note that $B(\La)_{v\La} = \{ b_{\La, v\La } \}$. If no confusion arises, we simply write $b_{v\La} = b_{ \La, v\La} $. We write $u_{v\La} \in V(\La) $ for the element corresponding to $b_{v\La}$, which is also called the extremal vector of $V(\La)$.

\subsection{Quiver Hecke algebras} \label{Sec: QHA}
In this subsection, we recall the definition of quiver Hecke algebras (see \cite{KL09, KL11, R08} and see also \cite{KKKO18, KKOP18, JP25}). 

\subsubsection{Quiver Hecke algebras}
For $i,j\in I$, let $\qQ_{i,j}(u,v) $ be polynomials over a field $\bR$ such that
$\qQ_{i,j}(u,v) = \qQ_{j,i}(v,u)$ and 
\begin{align*}
\qQ_{i,j}(u,v) =\bc
                   \sum\limits
_{p(\alpha_i , \alpha_i) + q(\alpha_j , \alpha_j)=-2(\alpha_i , \alpha_j)} t_{i,j;p,q} u^pv^q &
\text{if $i \ne j$,}\\[1.5ex]
0 & \text{if $i=j$,}
\ec
\end{align*}
where $t_{i,j; -c_{ij},0} \in  \bR^{\times}$.
For $\beta  \in \rlQ_+$ with $n = \Ht(\beta)$,
the {\em quiver Hecke algebra} $R(\beta)$ associated with  $(\qQ_{i,j}(u,v))_{i,j\in I}$
is the $\bR$-algebra generated by
$$
\{e(\nu) \mid \nu \in I^\beta \}, \; \{x_k \mid 1 \le k \le n \},
\; \{\tau_l \mid 1 \le l \le n-1 \}
$$
satisfying certain defining relations (see \cite{KL09, KL11, R08} and see also \cite[Section 2]{KKKO18}), where 
$$
I^\beta=  \left\{\nu=(\nu_1, \ldots, \nu_n ) \in I^n \mid \sum_{k=1}^n\alpha_{\nu_k} = \beta \right\}.
$$
The $\Z$-grading of $R(\beta)$ is defined as 
\begin{align*}
\deg(e(\nu))=0, \quad \deg(x_k e(\nu))= ( \alpha_{\nu_k} ,\alpha_{\nu_k}), \quad  \deg(\tau_l e(\nu))= -(\alpha_{\nu_{l}} , \alpha_{\nu_{l+1}}).
\end{align*}
We write $\wt(L):= - \beta \in -\rlQ_+$ for an $R(\beta)$-module $L$, and set $R := \bigoplus_{\beta\in \rlQ_+} R(\beta)$.
When $\cmC$ is symmetric and $\qQ_{i,j}(u,v)$ is a polynomial in $u-v$ for all $i,j\in I$, 
the algebra $R(\beta)$ is called \emph{symmetric}.

Let $R(\beta)\gmod$ be the category of finite-dimensional graded $R(\beta)$-modules. We write $R\gmod := \bigoplus_{\beta\in \rlQ_+} R(\beta)\gmod$. 
For graded $R(\beta)$-modules $M$ and $N $, $Hom_{R}(M,N)$ denotes the space of degree preserving homomorphisms, and define
\[
\Hom_{R}( M,N ) \seteq \bigoplus_{k \in \Z} Hom_{R}( q^kM,N ),
\]
where $q$ is the {\em grading shift functor}.
If no confusion arises, we write $M\simeq N$ if $M$ and $N$ are isomorphic to each other up to a grading shift.

We set $e(\beta, \gamma) := \sum_{\nu \in I^\beta, \mu \in I^{\gamma}} e(\nu \cdot \mu)  $, where $\nu \cdot \mu$ is the concatenation of $\nu $ and $ \mu$. For $L\in R(\beta+\gamma)\gmod$, $M\in R(\beta)\gmod$ and $N\in R(\gamma)\gmod$, the \emph{restriction} $\Res_{\beta, \gamma} (L)$ is defined as 
$$
\Res_{\beta, \gamma} (L) := e(\beta, \gamma ) L,
$$
and the \emph{convolution product} $M\conv N \in R(\beta+\gamma)\gmod$ is defined as
$$
M \circ N := R(\beta+\gamma) e(\beta, \gamma)  \tens_{R(\beta)\tens R(\gamma)} \left( M \tens N  \right), 
$$
where $M \tens N$ is the $R(\beta)\tens R(\gamma)$-module whose underlying vector space is the tensor product of $M$ and $N$ over the base field $\bR$.
Note that there is an isomorphism   
$$
\Hom_{R(\beta+\gamma)}(X\conv Y, Z) \simeq \Hom_{R(\beta)\tens R(\gamma)}(X\tens Y, \Res_{\beta, \gamma}Z)
$$
for $X \in R(\beta)\gmod$, $Y \in R(\gamma)\gmod$ and $Z \in R(\beta + \gamma)\gmod$ (see \cite{LV11}).

Let $\phi$ be the involutive algebra automorphism of $R(\beta)$ defined by 
\begin{align*}
e(\nu_1,\ldots ,\nu_n)\mapsto e(\nu_n,\ldots ,\nu_1), \quad  x_k\mapsto x_{n+1-k}, \quad 
\tau_k\mapsto  -  \tau_{n-k}.
\end{align*}
We denote by $\phi_*(X)$ the pullback of a $R$-module $X$ via $\phi$ (\cite[Section 1.2.2]{KKOP23A}). 
We then have
\begin{align} \label{Eq: phi *}
\phi_*(M\conv N) \simeq \phi_*(N) \conv \phi_*(M) \qquad \text{for $R$-modules $M$ and $N$.}
\end{align}

We simply write $ M^{\circ k} := \underbrace{M\conv \cdots \conv M}_{k \text{ factors}}$ for $k \in \Z_{>0}$, and understand $M^{\circ 0}$ as the trivial module $\one$. 
If no confusion arises, then we simply write $M^k$ for $M^{\circ k}$. 
We set $M \hconv N := \hd( M \conv N  )$ and $M \sconv N := \soc( M \conv N  )$, where $\hd (X)$ (resp.\ $\soc(X)$) is the head (resp.\ socle) of $X$ for a module $X$. We say that simple modules $M$ and $N$ \emph{strongly commute} if $M \conv N $ is simple. 
A simple module $L$ is \emph{real} if $L \conv L$ is simple. If $L$ is not real, then $L$ is called \emph{imaginary}.

\subsubsection{Categorification}
It was proved in \cite{KL09, KL11, R08} that the Grothendieck ring $K(R\gmod)$ is isomorphic to the quantum unipotent coordinate ring $A_q(\n)$, i.e, 
$$
K(R\gmod) \simeq A_\Zq(\n).
$$
Let $[X]$ denote the element in $K(R\gmod)$ corresponding to $X \in R\gmod$. If no confusion arises, we simply write $X$ for $[X]$.
We denote by $\Irr(R\gmod)$ the set of the isomorphism classes of simple $R$-modules (up to grading shifts). For each $i\in I$, we denote by $\ang{i}$ the 1-dimensional $R(\al_i)$-module of degree 0, and write 
$\ang{i^m} := q^{m(m-1)/2} \ang{i}^{\circ m}$ for $m \in \Z_{\ge0}$. Note that $\ang{i^m}$ is simple for $m >0$ since $\ang{i}$ is real.
We set 
\begin{align*}
\tF_i(M) := \ang{i}\hconv M \qtq  \tE_i(M) := \soc(\Res_{\al_i, \beta-\al_i}(M)), \\
\tFs_i(M) :=  M\hconv \ang{i} \qtq  \tEs_i(M) := \soc(\Res_{\beta-\al_i, \al_i}(M))
\end{align*}
for any $M \in \Irr(R\gmod)$. Then we have the following.
\bna
\item $\Irr(R\gmod)$ together with $\tF_i, \tE_i, \tFs_i$ and $\tEs_i$ has a crystal structure isomorphic to $B(\infty)$ (see \cite{LV11}), i.e., 
\begin{align} \label{Eq: Irr B infty}
\Irr(R\gmod) \simeq B(\infty).
\end{align}
Under this isomorphism, the pullback $\phi_*$ given in \eqref{Eq: phi *} corresponds to the $*$-involution of $B(\infty)$.
\item When $R$ is symmetric and $\bR$ is of characteristic $0$, $\Irr(R\gmod)$ corresponds to  $\Gup(\infty)$ under the categorification $K(R\gmod)\simeq A_q(\n)$ (\cite{VV11, R12}).
\ee
Note that the trivial module $\one$ is the highest weight vector of $\Irr(R\gmod)$ as a crystal. 
For any $b\in B(\infty)$, we denote by $\cL(b)$ the self-dual simple $R$-module corresponding to $b$ under the crystal isomorphism \eqref{Eq: Irr B infty}.

\subsubsection{Integer-valued invariants} \label{Sec: iv La de}
A pair $(M,N)$ is \emph{$\La$-definable} if there is a non-zero homomorphism $\rmat{}$ such that $\Hom_R (M\conv N,N\conv M)=\bR \rmat{}$. We call $\rmat{}$ the \emph{R-matrix}, and define 
\begin{align*}
\La(M,N)\seteq  \deg(\rmat{}) \qtq \de(M,N) &\seteq \dfrac{1}{2} \left(\La(M,N)+\La(N,M) \right).
\end{align*}
If one of $M$ and $N$ is \emph{affreal}, that is real and admits an \emph{affinization}, then $(M,N)$ and $(N,M)$ are $\La$-definable (\cite{KP18}). Note that every simple module over a \emph{symmetric} quiver Hecke algebra admits an affinization (\cite[Section 2.2]{KKKO18}). 
For $i\in I$ and a simple module $M$, we define  
\begin{align} \label{Eq: def of di}
\de_i(M) \seteq   \ep_i(M) + \ep^*_i(M) + \ang{h_i, \wt(M)}.
\end{align}
Note that $ \frac{(\al_i, \al_i)}{2} \de_i(M) = \de(\ang{i}, M)$ (\cite[Lemma 2.15]{KKOP23A}).

A sequence $ (L_1,\ldots,L_r)$ of simple modules is \emph{almost affreal} if all $L_i$ $(1 \le i \le r)$ are  affreal  except for at most one.
 An almost affreal sequence $ \underline{L} = (L_1,\ldots,L_r)$
is  a \emph{normal sequence} if the composition
\begin{align*}
\rmat{\underline{L} }  \seteq \prod_{1 \le i < k \le r} \rmat{L_i,L_k }   
& \colon \; q^{\sum_{1 \le i < k \le r} \La(L_i,L_k)} L_1 \conv \cdots \conv L_r \To L_r \conv \cdots L_1 
\end{align*}
does not vanish. Note that, if an almost affreal sequence $ (L_1,\ldots,L_r)$ is normal, then $\hd (L_1 \conv \cdots \conv  L_r)$ is simple and isomorphic to the image of $\rmat{\underline{L} }$, and 
\begin{align*}
\La(L_1, \hd(L_2 \conv \cdots \conv L_r)) &= \sum_{k\in [2,r]} \La(L_1, L_k), \\
\La(\hd(L_1 \conv \cdots \conv L_{r-1}), L_r) &= \sum_{k\in [1,r-1]} \La(L_k, L_r).
\end{align*}
(see \cite[Section 2.3]{KK19} and \cite[Section 2.2]{KKOP23A} for details).

\subsection{Categories $\cC_w$} \label{Sec: Cw}
In this subsection, we recall the definition of the subcategory $\cC_w$ of $R\gmod$ and its localization (\cite{KKOP18, KKOP21A, KKOP23A}). 

\subsubsection{Categories $\cC_w$} 
For $M \in R(\beta)\gmod$, we define 
\begin{align*}
\gW(M) & := \{  \gamma \in  \rlQ_+ \cap (\beta - \rlQ_+)  \mid  \Res_{\gamma, \beta-\gamma}  M \ne 0  \}, \\ 
\gW^*(M) & := \{  \gamma \in  \rlQ_+ \cap (\beta - \rlQ_+)  \mid  \Res_{\beta-\gamma, \gamma } M \ne 0  \}. 
\end{align*}
For $M, N \in \Irr(R\gmod)$, a pair $(M,N)$ is \emph{unmixed} if 
$$
\sgW(M) \cap \gW(N) = \emptyset.
$$ 
\begin{lemma}[{{\cite{KK19},\cite[Lemma 2.8, Corollary 2.13]{KKOP23A}}}] \label{Lem: normal}
 Let $(L,M,N)$ be an almost affreal sequence of simple $R$-modules. 
Then the triple $(L,M,N)$ is normal if one of the following conditions holds:
\bna
\item $L$ is real and $\de(L,M)=0$,
\item $N$ is real and $\de(M,N)=0$,
\item $L$ or $N$ is real and $(L,N)$ is unmixed.
\ee
\end{lemma}

For any $w\in \weyl$, let $\cC_{w}$ be the  full subcategory of $R\gmod$ consisting of objects $M$ satisfying $\gW(M) \subset \Sp( \prD \cap w \nrD ).$
Then the category $\cC_w$ gives the categorification of $A_q(\n(w))$, i.e., 
\begin{align} \label{Eq: categorification Aqnw}
K(\cC_w) \simeq A_\Zq(\n(w))
\end{align}
(see \cite{KL09, KL11, R08} and see also \cite{KKKO18, KKOP18}).
By the definition of $B(w)$ (see \eqref{Eq: Bw}) and \eqref{Eq: Irr B infty}, 
the set $\Irrw{w}$ of the isomorphism classes of simple modules in $\cC_w$ is in 1-1 correspondence with $B(w)$, i.e.,    
\begin{align} \label{Eq: Irr B}
\Irrw{w} = \{ \cL(b) \mid b \in B(w)\}.
\end{align}

\subsubsection{Determinantial modules} \label{Sec: DM} 
For $\La \in \wlP_+$ and  $ \bfi = (i_1, i_2, \ldots, i_r) \in \Seq(w)$ of $ w \in \weyl$, we define
$$
\dM(w\La, \La) := \tF_{i_1}^{c_1} \tF_{i_2} ^{c_2} \cdots \tF_{i_r} ^{c_r} \one,
$$
where $c_k = \ang{h_{i_k}, s_{i_{k+1}} \cdots s_{i_{r}}\La }$ for $k \in [1,r]$. 
For $v\in W$ with $v\le w$ and $\bfj = (j_1, j_2, \ldots, j_t) \in \Seq(v)$, we define 
$$
\dM(w\La, v\La) := (\tE_{j_1}^*)^{  d_1} (\tE_{j_2}^*)^{d_2} \cdots (\tE_{j_t}^*)^{d_t} (\dM(w\La, \La)),
$$
where $d_k = \ang{h_{j_k}, s_{j_{k+1}} \cdots s_{j_{t}}\La } $ for $k\in [1,t]$.
This simple module $\dM(w\La, v\La)$ is called the \emph{determinantial module}.
Note that the simple module $\dM(w\La, v\La)$ does not depend on the choice of $\bfi$ and $\bfj$. We have the following.
\bna 
\item Under the categorification \eqref{Eq: categorification Aqnw}, $\dM(w\La, v\La)$ corresponds to the \emph{quantum unipotent minor} $\qmD (w \La, v\La)$, i.e., $[\dM(w\La, v\La)] =  \qmD (w \La, v\La)$  (see \cite{KKKO18} and see also \cite[Section 4]{KKOP18} for details). Note that $\qmD (w \La, v\La)$ is contained in $\Gup(w)$.
\item $\dM(w\La, v\La)$ is affreal (\cite[Proposition 4.2]{KKOP18} and \cite[Theorem 3.26]{KKOP21A}). 
\item $\dM(w\La, \La)$ commutes strongly with any simple in $\cC_w$ (\cite[ Proposition 5.1]{KKOP21A}).
\ee

For $ \bfi = (i_1, i_2, \ldots, i_r) \in \Seq(w)$, we write $w_{\le 0} := \id$, $w_{\le k} := s_{i_1} s_{i_2} \cdots s_{i_k} $ and $w_{< k} := w_{\le k-1}$ for $k\in [1,r]$. For each $k \in [1,r]$, we define
\begin{align} \label{Eq: cuspS}
\cuspS_{\bfi, k} := \dM(w_{\le k}\La_{i_k}, w_{<k} \La_{i_k}),
\end{align}
which is called the $k$th \emph{cuspidal module} with respect to $\bfi$.
Under the categorification \eqref{Eq: categorification Aqnw}, the cuspidal module $\cuspS_{\bfi, k}$ corresponds to the dual PBW vector $\rv^{\up}_{\bfi} (\beta_k)$, i.e., $[\cuspS_{\bfi, k}] =  \rv^{\up}_{\bfi} (\beta_k)$ for $k\in [1,r]$ (see \cite[Corollary 4.7]{KKOP18}). For a simple module $X \in \Irrw{w}$ such that  $X\simeq \cL(b)$ for some $b\in B(w)$ (see \eqref{Eq: Irr B}), we set 
\begin{align} \label{Eq: PBW for QHA}
\PBW_\bfi(X) := \PBW_\bfi(b)
\end{align}
(see \eqref{Eq: PBWi}). Note that, if $\PBW_\bfi(X) = (a_1, \ldots, a_r)$, then 
$X \simeq \hd (\cuspS_{\bfi, r}^{\circ a_r} \conv \cdots \conv \cuspS_{\bfi, 2}^{\circ a_2} \conv \cuspS_{\bfi, 1}^{\circ a_1})$ (see \cite{McNa15} and \cite{KKOP18}).

A weight $\la \in \wlP$ is \emph{$w$-dominant} if $(\beta, \la) \ge0 $ for any $\beta \in \prD \cap w^{-1} \nrD$. Note that a dominant weight $\La \in \wlP_+$ is  $w$-dominant for any $w \in \weyl$.

\begin{theorem} [{\cite[Theorem 2.18]{KKOP23A}}] \label{Thm: gdM}
Let $w\in W$ and let $\la\in\wlP$.
Assume that $\la$ is $w$-dominant.
Then there exists a self-dual simple $R(\la-w\la)$-module
$\dM_w(w\la,\la)$ which satisfies the following conditions.
\bna
\item 
If $i\in I$ satisfies $s_iw<w$, then
$\dM_w(w\la,\la)\simeq  \ang{i^{m}} \hconv \dM_{s_iw}(s_iw\la,\la)$
where $m=\ang{h_i,s_iw\la}\in\Z_{\ge0}$.
\item
If $i\in I$ satisfies $ws_i<w$, then
$\dM_w(w\la,\la)\simeq \dM_{ws_i}(w\la,s_i\la)\hconv \ang{i^{m}} $
where $m=\ang{h_i,\la}\in\Z_{\ge0}$.
\item
If $i\in I$ satisfies $\ang{h_i,w\la}\ge0$,
then $\ep_i\bl\dM_w(w\la,\la)\br=0$.

\item
If $i\in I$ satisfies $\ang{h_i,\la}\le0$,
then $\ep^*_i \bl \dM_w(w\la,\la) \br=0$.
\item
For any $\mu\in\wlP_+$ such that
$\la+\mu\in\wlP_+$,  we have
$$\dM(w\mu,\mu)\hconv \dM_w(w\la,\la)\simeq\dM\bl w(\la+\mu),\la+\mu\br$$
up to a grading shift.
\ee
Moreover such an $\dM_w(w\la,\la)$ is unique up to an isomorphism.
\end{theorem}
We call $\dM_w(w\la,\la)$ the \emph{generalized determinantial module} associated with $w$ and $\la$. The following will be used crucially in Section \ref{Sec: ct app}.
\bna
\item If $\la \in \weyl \La$ for $\La \in \wlP_+$, the generalized determinantial module $\dM_w(w\la,\la)$ coincides with  $\dM(w\la,\la)$ (\cite[Remark 2.22 (i)]{KKOP23A}).
\item For a $w$-dominant $\la \in \wlP$, 
there is  an isomorphism 
\begin{align} \label{Eq: Phi * dM}
\phi_*(\dM_w(w\la,\la))
\simeq\dM_{w^{-1}}(-\la,-w\la),
\end{align}
where $\phi_*$ is the pullback given in \eqref{Eq: phi *} (\cite[Lemma 2.23]{KKOP23A}).

\item 
Let $w \in \weyl$, $\la \in \wlP_+$, and $\la' \in \sum_{j\in I \setminus \supp(w)} \Z \La_j$ such that $\la+ \la'$ is $w$-dominant. 
Then, by the definition of $\dM_w( w \la,  \la )$ (see Theorem \ref{Thm: gdM}), we have 
\begin{align} \label{Eq: Mw la la'}
\dM_w( w \la,  \la ) = \dM_w( w (\la + \la'), \la+\la' ). 
\end{align}

\ee

\vskip 2em 

\section{Quantum twist automorphisms} \label{Sec: QTA}

In this section, we review the crystal-theoretic counterpart $\cdual_w$, introduced in \cite{JP25}, of Kimura-Oya's quantum twist automorphism $\eta_w$ (\cite{KiOy21}).
Let $ (\cmC,\wlP,\Pi,\wlP^\vee,\Pi^\vee) $ be a Cartan datum and let $A_q(\n)$ and $R$ be the corresponding quantum unipotent coordinate ring and quiver Hecke algebra respectively.

\subsection{Localization} \label{Sec: Loc}
Let $w \in \weyl$ with $\supp(w) = I$. Recall the determinantial modules $\dM(w\La, \La) \in \cC_w$ and quantum unipotent minors $\qmD(w\La, \La)$ (see Section \ref{Sec: DM}).
Since the quantum unipotent minors $\qmD (w \La_i, \La_i)$ are central in $A(\n(w))$ up to powers of $q$, one can localizes $A_q(\n(w))$ by $\qmD (w \La_i, \La_i)$ for $i\in I$. We denote by $\tA_q(\n(w))$ the resulting localization. 
Then we have the basis $\tGup (w)$ of $\tA_q(\n(w))$ described as 
\begin{align*} 
\tGup (w) = \{ q^{a(x,\bfa)} x \cdot  \qmD(\bfa) \mid x \in \Gup(w),\ \bfa = (a_i)_{i\in I} \in \Z^{\oplus I}  \}, 
\end{align*}
where 
$
\qmD(\bfa) := \prod_{i\in I} ( \qmD ({w \La_{i}, \La_{i}} ) )^{a_i}
$
and $a(x,\bfa)$ is chosen so that $q^{a(x,\bfa)} x \cdot  \qmD(\bfa)$ is \emph{dual bar-invariant} (\cite[Proposition 4.5.]{KiOy21}).

We now define the crystal-theoretic counterpart of $\tGup (w)$.
The set $B(w) \times \Z^{\oplus I}$ has the equivalent relation $\sim$ defined as follows:
For any $ (b, \bfa), (b', \bfa')  \in B(w) \times \Z^{\oplus I}$, we write $ (b, \bfa) \sim (b', \bfa') $ if $ \Gup(b) \qmD(\bfa) \eqq  \Gup(b') \qmD(\bfa')$. 
Define   
\begin{align*} 
\LB(w) := B(w) \times \Z^{\oplus I} / \sim, 
\end{align*}
which is called the \emph{localized crystal}. For any $x = (b, \bfa) \in \LB(w) $, we denote by $\Gup(x) $ the element of $\tGup (w)$ such that  $ \Gup(x) \eqq\Gup(b) \qmD(\bfa)$. 
By construction, $\LB(w)$ is in 1-1 correspondence with $\tGup(w)$ (see \eqref{Eq: Gup B}), i.e., 
\begin{align} \label{Eq: tGup tB}
\tGup (w) = \{ \Gup(x) \mid x \in \LB(w) \}.
\end{align}

The categorical counterpart of $\tA_q(\n(w))$ is the localized category $\tcC_w$, which is the \emph{localization} of the category $\cC_w$ by the determinantial modules $\dM(w\La_i, \La_i)$ for $i\in I$ (see \cite[Section 5]{KKOP21A} for precise definition). By construction, we have  
\begin{align} \label{Eq: categorification localized}
K(\tcC_w) \simeq \tA_\Zq(\n(w)).
\end{align}
We denote by $\Irrtw{w}$ the set of the isomorphism classes of simple objects in $\tcC_w$. Since the simplicity is preserved under the natural functor $\cC_w \rightarrow \tcC_w$, we regard $\Irrw{w}$ as a subset of $\Irrtw{w}$. 
Then we have 
\begin{align*} 
\Irrtw{w} = \{ L \conv  \dM(\bfa) \mid L \in \Irrw{w},\ \bfa = (a_i)_{i\in I} \in \Z^{\oplus I}  \}, 
\end{align*}
where $\dM(\bfa) := \conv_{i\in I} ( \dM ({w \La_{i}, \La_{i}} ) )^{a_i}$ (see \cite[Section 5.1]{KKOP21A}). For $x = (b, \bfa) \in \LB(w)$, we write $\cL (x) = \cL(b) \conv  \dM(\bfa) \in \Irrtw{w}$. 
Then $\LB(w)$ is in 1-1 correspondence with $\Irrtw{w}$ (see \eqref{Eq: Irr B}), i.e.,
$$
\Irrtw{w} = \{ \cL(x) \mid x \in \LB(w) \}. 
$$

\subsection{Quantum twist automorphism} \label{Sec: qta}

Let $\eta_w$ be the \emph{quantum twist automorphism} on $\tA_q(\n(w))$ introduced by Kimura--Oya (\cite{KiOy21}).
The $\eta_w$ is the $\Q(q)$-automorphism of $\tA_q(w)$ 
defined by 
\begin{align} \label{Eq: eta_w}
\eta_w(\qmD (v \La, \La))  \eqq
\qmD (w \La, \La)^{-1} \qmD(w\La, v\La)
\end{align}
for any $v \in \weyl$ with $v \le w$ and $\La \in \wlP_+$ (see \cite[Theorem 6.1]{KiOy21} for precise definition). Since $\eta_w$ permutes the dual canonical basis $\tGup(w)$ (\cite[Theorem 6.1]{KiOy21}), the map $\eta_w$ induces a permutation of $\LB(w)$ via \eqref{Eq: tGup tB}. We denote by  $\cdual_w$ the permutation of $\LB(w)$ induced by the inverse $\eta_w^{-1}$ of $\eta_w$.

Since the category $ \tcC_w $ is \emph{rigid} (\cite{KKOP23A}), 
there exists the right dual functor $\dual_w$ in $\tcC_w$. 
It was shown  that the right dual functor $\dual_w$ categorifies the quantum twist automorphism $\eta_w^{-1}$  on $\tGup(w)$  (see \cite[Theorem B.22]{IO21_arXiv2} and see also \cite[Proposition 1.4]{JP25}). 
Thus we have the following commutative diagram:
\begin{equation} \label{Eq: DDeta} 
\begin{aligned} 
\xymatrix{
 \Irrtw{w} \ar[d]_{ \dual_w }  &&  \ar[ll]_\sim^{\cL}  \LB(w)  \ar[rr]^{\sim}_{\Gup} \ar[d]_{\cdual_w} && \tGup(w) \ar[d]_{\eta_w^{-1}} \\
\Irrtw{w} && \ar[ll]_\sim^{\cL} \LB(w) \ar[rr]^\sim_{\Gup} && \tGup(w)
}
\end{aligned}
\end{equation}
(see \cite[Section 2]{JP25} for details).

\subsection{Crystal-theoretic description of $\cdual_w$} \label{Sec: crystal para}
In this subsection, we recall the crystal-theoretic description of $\cdual_w$ given in \cite[Theorem 3.3]{JP25} in terms of PBW and string data.

Recall $\PBW_\bfi (b)$ and $\STR_\bfi (b)$ defined in \eqref{Eq: PBWi} and \eqref{Eq: STRi}. We choose and fix a reduced sequence $ \bfi = (i_1, i_2, \ldots, i_r) \in \Seq(w)$.
\bna 
\item  \textbf{(PBW parameterization)} \ 
Let $\tcP_\bfi(w) := \Z^r$ and set $P_i := \PBW_\bfi(\qmD(w\La_i, \La_i))$ for $i\in I$. For any $x = (b, (a_i)_{i\in I}) \in \LB(w)$, we set
$$
\PBW_\bfi(x) := \PBW_\bfi(b) + \sum_{i\in I} a_i P_i \in \Z^{r}.
$$
Then the map $\pbwP_\bfi :  \LB(w) \buildrel \sim \over \longrightarrow \tcP_\bfi(w), \  x \mapsto \PBW_\bfi(x)$
is bijective (\cite[Section 2.1]{JP25}), and the map $\pbwP_\bfi$ 
is called the \emph{PBW parametrization} of $\LB(w)$ with respect to $\bfi$.  

\item  \textbf{(String parameterization)} \ 
Let $\tcS_\bfi(w) := \Z^r$ and set $S_i := \STR_\bfi(\qmD(w\La_i, \La_i))$ for $i\in I$. For any $x = (b, (a_i)_{i\in I}) \in \LB(w)$, we set
$$
\STR_\bfi(x) := \STR_\bfi(b) + \sum_{i\in I} a_i S_i \in \Z^{r}.
$$
Then the map $\strS_\bfi :  \LB(w) \buildrel \sim \over \longrightarrow \tcS_\bfi(w), \ x \mapsto \STR_\bfi(x)$
is bijective (\cite[Section 2.2]{JP25}), and the map $\strS_\bfi$ is called the \emph{string parametrization} of $\LB(w)$ with respect to $\bfi$.  

\ee

We define a bijection $\psi_{\bfi}: \tcP_\bfi(w) \buildrel \sim \over\longrightarrow \tcS_\bfi(w)  $ so that the following diagram commutes:
\begin{equation} \label{Eq: STR PBW isom}
\begin{aligned} 
\xymatrix@R0.7em{
& \LB(w) \ar[ld]_{\pbwP_\bfi} \ar[rd]^{\strS_\bfi} &
\\
  \tcP_\bfi(w) \ar[rr]^ {\psi_\bfi} && \tcS_\bfi(w), 
}
\end{aligned}
\end{equation}
and define square matrices  
$N_\bfi = (n_{p,q})_{p,q \in [1,r]}$ and  
$M_\bfi = (m_{p,q})_{p,q \in [1,r]}$
by
\begin{align*}
n_{p,q} = 
\begin{cases}
1 & \text{ if } p=q,  \\
-1 & \text{ if } p^+ = q, \\
0  & \text{ otherwise,}
\end{cases}
\qquad
m_{p,q} = 
\begin{cases}
 \langle h_{i_{ p}}, s_{i_{ p+1}} s_{i_{p+2}} \cdots s_{i_{q}} \La_{i_q}   \rangle  & \text{ if } p\le q,   \\
0  & \text{ otherwise,}
\end{cases}
\end{align*}
where $p^+ := \min \left( \{ t\in [1,r] \mid t > p, \ i_t = i_p \} \cup \{ +\infty \} \right)$. 

\begin{theorem} [{\cite[Theorem 3.3]{JP25}}] \label{Thm: ta}
Let $w\in \weyl$ and $\bfi \in \Seq(w)$.
For any $b \in \LB(w)$,
\bni 
\item $ \PBW_{\bfi}( \cdual_w(b) )  =  \psi_\bfi^{-1} \circ M_\bfi \circ N_\bfi (-\PBW_{\bfi}(b )) $,
\item $ \STR_{\bfi}( \cdual_w(b) )  = -   M_\bfi \circ N_\bfi \circ \psi_\bfi^{-1} (\STR_{\bfi}(b )) $.
\ee
\end{theorem}

\vskip 2em 

\section{Stable imaginary elements} \label{Sec: PI}

In this section, we introduce the notion of \emph{\stb imaginary} elements using the quantum twist automorphism and investigate their properties. 

Let $A_q(\n)$ be the unipotent coordinate ring associated with a symmetrizable generalized Cartan matrix $\cmC$, and let $R$ be the corresponding quiver Hecke algebra.
Recall $\eta_w$, $ \cdual_w$ and  $ \dual_w$ given in Section \ref{Sec: qta}.

\begin{definition} \label{Def: frozen}
Let $w \in W$.
\bna
\item An element $x \in \Gup(w)$ is called \emph{frozen} if $x \eqq \qmD(w\La, \La)$ for some $\La \in \wlP_+$.
\item An element $b \in B(w)$ is called \emph{frozen} if $\Gup(b)$ is frozen. 
\item A simple module $X \in \Irrw{w}$ is called \emph{frozen} if $X \simeq \dM(w\La, \La)$ for some $\La \in \wlP_+$.
\ee
\end{definition}
From the perspective of cluster algebras, $\qmD (w \La, \La)$ corresponds to a monomial in frozen variables of $A_q(\n(w))$.
Note that, for $\La$ and $ \La' \in \wlP_+$, 
\begin{equation} \label{Eq: additivity for D M}
\begin{aligned}
\qmD (w \La, \La) \cdot \qmD (w \La', \La') & \eqq \qmD (w \La+\La', \La+\La'),\\ 
\dM (w \La, \La) \conv \dM (w \La', \La') &\simeq \dM (w \La+\La', \La+\La')
\end{aligned}
\end{equation}
(see \cite[Lemma 1.6, Proposition 4.2]{KKOP18}).

For any $x,y \in \tA_q(\n(w))$, we define
\begin{align} \label{Eq: eqfr for G}
x \eqfr y \quad \Longleftrightarrow \quad \qmD (w \La, \La)   x \eqq  y \quad \text{for some  $\La  \in \wlP_+$,}
\end{align}
and, for any $b,b' \in \tB(w)$, define
\begin{align} \label{Eq: eqfr for B}
b \eqfr b' \quad \Longleftrightarrow \quad  \Gup(b)  \eqfr  \Gup(b').
\end{align} 
Using the parameterizations given in Section \ref{Sec: crystal para}, we also write 
$\PBW_{\bfi} (b) \eqfr \PBW_{\bfi} (b')  $ and $\STR_{\bfi} (b) \eqfr \STR_{\bfi} (b')  $ if $b \eqfr b'$.
Similarly,  for any objects $M$ and $N$ in $\tcC_w$, we write 
\begin{align} \label{Eq: eqfr for Irr}
M \eqfr N \quad \Longleftrightarrow \quad \dM(w\La, \La)  \conv M \simeq  N \quad \text{for some  $\La  \in \wlP_+$.}
\end{align}

\begin{lemma} \label{Lem: fixed by D}
Let $x \in \Gup(w)$  and $X \in \Irrw{w}$.
\bni
\item If $ \eta_w(x) \eqfr x$, then $x$ is either frozen or imaginary.
\item If $\dual_w(X) \eqfr X $, then $X$ is either frozen or imaginary. 
\ee
\end{lemma}
\begin{proof}
 We first prove (ii).
We assume that $X$ is real since otherwise it is trivial. 
Since $\dual_w (X)$ is a right dual, there is a nonzero map $X \conv \dual_w(X) \twoheadrightarrow \one$ in the localized category $\tcC_w$, where $\one$ is a unit object in $\tcC_w$ (e.g., \cite[Section 1.5]{KKOP21A}). Note that every frozen modules are invertible and $q$-central in $\tcC_w$. 
By the assumption $\dual_w(X) \eqfr X $, there are frozen modules $F_1$ and $F_2$ such that  $ \dual_w(X) \conv F_2 \simeq  X \conv F_1  $ (see \cite[Section 5.1]{KKOP21A}). 
Then we have $X \conv \dual_w(X) \conv F_2 \simeq X \conv X \conv F_1$, and $X \conv X \conv F_1$ is simple in $\cC_w$ since $X$ is real. Thus the map $X \conv \dual_w(X) \twoheadrightarrow \one$ implies  $ X \circ X \circ F_1  \simeq F_2$.
Since $[\dM (w \La_i, \La_i) ]$ are prime elements in $K(\cC_w) \simeq A_q(\n(w))$, the simple module $X \circ X$ is a frozen module by \eqref{Eq: additivity for D M}. As $ \PBW_\bfi(X\conv X) = 2\, \PBW_\bfi(X) $ for  $\bfi \in \Seq(w)$, we conclude that $X$ itself is frozen.

We now prove (i). If $U_q(\g)$ is of symmetric type, then it follows from (ii) and the fact that $\Gup(w)= [\Irrw{w}]$ under the categorification (\cite{VV11, R12}).
Suppose that $U_q(\g)$ is of non-symmetric type, and assume that $x$ is real.
By Proposition \ref{prop: eta x x} in Appendix \ref{Sec: appA}, the expansion of $\eta_w(x)x$ with respect to $\tGup(w)$ has a nonzero constant term. Since $x^2 \eqfr \eta_w(x)x$ belongs to $\tGup(w)$ up to a power of $q$, we conclude that $x^2$ is frozen. Then, by the same argument as in (ii), we conclude that $x$ is frozen.
\end{proof}

\begin{definition} \label{def: Dw PI}
Let $w \in \weyl$. 
\bnum
\item \label{def: Dw PI qha c}
An element $b \in B(w)$ is called a \emph{$\cdual_w$-\stb imaginary} element if $\cdual_w(b) \eqfr b $ and $b$ is not frozen. 

\item \label{def: Dw PI Gup}
An element $x \in \Gup(w)$ is called \emph{$\eta_w$-\stb imaginary} if $\eta_w(x) \eqfr x $ and  $x$ is not frozen.

\item \label{def: Dw PI qha}
A simple module $L \in \Irrw{w}$ is called \emph{$\dual_w$-\stb imaginary} if $ \dual_w(L) \eqfr L$ and $L$ is not frozen.
\ee
\end{definition}

By the commuting diagram \eqref{Eq: DDeta}, we have the following.
\begin{lemma} \label{lem: equiv b}
For any $b \in B(w)$, the following are equivalent:
\bna
\item $b$ is $\cdual_w$-\stb imaginary, 
\item $\Gup(b)$ is $\eta_w$-\stb imaginary, 
\item $ \cL(b)$ is $\dual_w$-\stb imaginary. 
\ee
\end{lemma}

\begin{definition} \label{def: k PBW}
Let $b \in B(w)$, and choose  $\bfi \in \Seq(w)$.
For any $m\in \Z_{>0}$, we define $b^{[m]} \in B(w)$ to be the element of $B(w)$ such that
\begin{align*} 
\PBW_{\bfi} (b^{[m]}) =  m \, \PBW_{\bfi} (b).
\end{align*}
We call $ b^{[m]}$ the \emph{$m$-th multiple} of $b$.
\end{definition}
Note that, by \cite[Section 3]{BZ97} (and also \cite[Section 2.3]{SST18}),  we have 
\begin{align} \label{Eq: k mul of b}
\PBW_{\bfj} (b^{[m]}) =  m \, \PBW_{\bfj} (b) \qquad \text{ for any choice of $\bfj \in \Seq(w)$,}
\end{align}
which says that Definition \ref{def: k PBW} is well-defined.

\begin{prop} \label{Prop: bk pi}
Let  $U_q(\g)$ be of finite type and $w \in \weyl$. If $b \in B(w)$ is $\cdual_w$-\stb imaginary, then so is $b^{[k]}$ for any $k \in \Z_{>0}$.
\end{prop}
\begin{proof}
Let $w_0$ be the longest element. 
Let $i\in I$ and take $\bfi = (i_1, i_2, \ldots, i_\ell) \in \Seq(w_0)$ such that $i_1 = i^*$, where $\ell = \ell(w_0)$  and  $\al_{i^*} = -w_0(\al_i)$.  Since $\beta_\ell = s_{i_1} \cdots s_{i_{\ell-1}}(\alpha_{i_\ell}) = \al_i$,  we have the PBW vector $\rv^{\up}_{\bfi} (\beta_\ell)=  f_i $. Letting $ \PBW_{\bfi} (b)  = (a_{k})_{k\in [1,\ell]} $, we have
$$
\ep_i(b) = a_\ell \qtq \PBW_{\bfi} (\tem_i(b)) = (a_1, \ldots, a_{\ell-1},  0)
$$
(see \cite[Proposition 2.12]{SST18} for example). Thanks to \eqref{Eq: k mul of b}, for any $k\in \Z_{>0}$, we have 
\begin{align} \label{Eq: epi bk}
\ep_i(b^{[k]}) = k \ep_i(b) \qtq (\tem_i(b))^{[k]} = \tem_i(b^{[k]}).
\end{align}

Let $b$ be $\cdual_w$-\stb imaginary. It follows from \eqref{Eq: epi bk} that the isomorphism $\psi_{\bfi}: \tcP_\bfi(w) \buildrel \sim \over\longrightarrow \tcS_\bfi(w)  $ given in \eqref{Eq: STR PBW isom} is compatible with multiplying by $k > 0$, i.e., 
$$
\psi_{\bfi} (k \cdot \bfa ) = k \cdot \psi_{\bfi} ( \bfa ).
$$
Since $b$ is \stb imaginary, 
 Theorem \ref{Thm: ta} tells us that 
\begin{align*}
\PBW_{\bfi} (\cdual_w(b^{[k]})) & =  \psi_\bfi^{-1} \circ M_\bfi \circ N_\bfi (-\PBW_{\bfi}(b^{[k]} )) =  \psi_\bfi^{-1} \circ M_\bfi \circ N_\bfi ( k (-\PBW_{\bfi}(b) ) ) \\
&= k \left( \psi_\bfi^{-1} \circ M_\bfi \circ N_\bfi (  -\PBW_{\bfi}(b) ) \right) = k \left( \PBW_{\bfi} (\cdual_w(b)) \right)\\
&\eqfr k \left( \PBW_{\bfi} (b) \right)  =  \PBW_{\bfi} (b^{[k]}),
\end{align*} 
which completes the proof.
\end{proof}

\begin{remark} \label{rmk: pi and i}
In general, imaginary modules in $ \Irrw{w}$ are not necessarily $\dual_w$-\stb imaginary. Let $L := \cL( b_{\triC_{3}})$ be the imaginary $R_{C_3}$-module corresponding to the $\cdual_{w_0}$-\stb imaginary element $b_{\triC_{3}}$ given in Proposition \ref{Prop: triX in ft} \ref{P3} \ref{P3: b} (see also Section \ref{sec: Lec B3C3}).  Since $\de_1(L)=0$,  the modules $L$ and $\ang{1}$ commute strongly. Thus, $L\conv \ang{1}$ is an imaginary simple module, which is not fixed by $\dual_{w_0}$ (up to a multiple of frozens) since $\ang{1}$ is not fixed by $\dual_{w_0}$.
\end{remark}

\begin{remark} \label{rmk: new def of pi} 
The notion of \stb imaginary modules and multiples of crystal elements can be understood in a more abstract manner.
We are grateful to M. Kashiwara for explaining the following perspective to us.

\bni 
\item 
In the quiver Hecke algebra setting, we can consider a more abstract definition as follows: a simple module $L \in R\gmod$ is called \emph{\stb imaginary} if there exist an affreal simple module $C \in R\gmod$ and a non-zero homomorphism $f: L \conv L \rightarrow C$ such that  (a) $f$ is not an isomorphism,  and (b) $\de(L,C)=0$. If a simple module $L$ is $\dual_w$-\stb imaginary, then there is an epimorphism  $ L \conv L \twoheadrightarrow F $ with nontrivial kernel for some frozen module $F$. Thus every $\dual_w$-\stb imaginary module satisfies the conditions (a) and (b) for any $w \in \weyl$. In this sense, this new definition is more abstract and intrinsic, and it generalizes Definition \ref{def: Dw PI} \ref{def: Dw PI qha} in the quiver Hecke algebra case.
However, in the paper, we will focus on $\dual_w$-\stb imaginary modules since most of the results in the paper rely on the properties of $\dual_w$.
\item 
The notion of $b^{[m]}$ can be extended to elements  of $ B(\infty)$ using the similarity of crystals (\cite{Kash95}). Let $S_\infty: B(\infty) \rightarrow B(\infty)$ be the similarity map corresponding to $m\in \Z_{>0}$ introduced in \cite[Theorem 3.2]{Kash95}, and define $ b^{[m]} := S_\infty(b)$. 
Then this definition agrees with Definition \ref{def: k PBW} whenever $b\in B(w) \subset B(\infty)$, which means that $b^{[m]}$ does not depend on the choice of $w$. 
\ee
\end{remark}

\vskip 2em 

\section{Crystals and imaginary elements} \label{Sec: crystal and pi}

In this section, we present a condition to find $\cdual_w$-\stb imaginary elements in terms of crystals.  In what follows, we consider $B(\La)$ as a subset of $B(\infty)$ via the embedding \begin{align} \label{Eq: iota BLa Binf}
\iota_\La: B(\La) \hookrightarrow B(\infty)
\end{align}
defined by  
$\iota_\La (b_\La) = \one$ and $\iota_\La (\te_i(b)) = \te_i (\iota_\La (b))$ for any $b\in B(\La)$ and $i\in I$.

\subsection{Crystal-theoretic approach} \label{Sec: ct app}
We first prove several lemmas needed for the main theorem.
For $X \in \Irr(R\gmod)$, we simply write $X^k$ for $X^{\circ k}$ for $k\in \Z_{>0}$ if no confusion arises.

\begin{lemma} \label{Lem: criterion for imaginary}
Let $x, m \in B(\infty)$ and set $\beta := - \wt(x) \in \rlQ_+$. 
 Assume that 
\bna
\item $\Hom_{R(\beta) \tens R(\beta)}( \cL(x) \tens \cL(x), \Res_{\beta, \beta} (\cL(m))  ) \ne 0$, 
\item there exists $i\in I$ such that either $2 \ep_i(x) \ne \ep_i(m)$ or $2 \eps_i(x) \ne  \eps_i(m)$. 
\ee 
Then $x$ is imaginary. 
\end{lemma}
\begin{proof}
Let $X = \cL(x)$ and $M = \cL(m)$. We may assume that 
\begin{align} \label{Eq: 2ep ne}
2 \ep_i(x) \ne \ep_i(m)
\end{align} 
 since the other case for $2 \eps_j(x) \ne \eps_j(m)$ can be proved in the same manner.
By (a), we have a nonzero surjection $X\circ X \twoheadrightarrow M$.
If $X$ is real, then we have $ X \circ X \simeq M$. 
This tells us that 
$$
2 \ep_i(x) = 2\ep_i(X) = \ep_i(X \conv X) = \ep_i(M) = \ep_i(m), 
$$ 
which is a contradiction to \eqref{Eq: 2ep ne}.
\end{proof}

\begin{lemma} \label{Lemma: unique head}
Let $M \in \Irr(R\gmod)$ and   
$J \subset I$ such that 
$$
(\al_j, \al_{j'}) =0 \qquad \text{ for any $j,j' \in J$ with $j \ne j'$.} 
$$
Let $k_j \in \Z_{\ge0}$ for $j\in J$ and define $Z := \conv_{j\in J} \ang{j}^{  k_j}$.
Then, for any $a,b \in \Z_{\ge0}$ such that
$$
\de_j(M)  \ge k_j (a+b)  \qquad \text{ for all $j\in J$},
$$ 
we have 
$$
(Z^a \hconv M) \hconv Z^b \simeq Z^a \hconv (M \hconv Z^b).
$$ 
\end{lemma}
\begin{proof}
We shall use induction on $n := |J|$. 
If $n=1$, then the statement follows from \cite[Theorem 2.17]{KKOP23A}. 

Suppose that $n>1$. We choose $i \in J$ and set $J' := J \setminus \{i\}$, $Y := \conv_{j\in J'} \ang{j}^{ k_j} $. 
Note that $Y$ and $Z$ are affreal (\cite[Theorem 3.26]{KKOP21A}).
As $ (\al_i, \al_{j}) =0$ for any $j\in J'$, we have
\begin{align} \label{Eq: iY Yi}
\text{$(\ang{i}, Y)$ and $(Y,\ang{i})$ are unmixed, and $\de(\ang{i}, Y)=0$.}
\end{align}
By the induction hypothesis, we have 
\begin{align} \label{Eq: YaM}
(Y^a \hconv M) \hconv Y^b \simeq Y^a \hconv (M \hconv Y^b).
\end{align}
Combining Lemma \ref{Lem: normal} together with \eqref{Eq: iY Yi},
we have  
\begin{equation} \label{Eq: ZaMYb}
\begin{aligned}
(Z^a \hconv M )\hconv Y^b &\simeq \left( (\ang{i}^{k_ia} \conv Y^a) \hconv M \right)\hconv Y^b \simeq \left(  \ang{i}^{k_ia} \hconv (Y^a \hconv M) \right)\hconv Y^b \\
&\simeq  \ang{i}^{k_ia} \hconv \left(  (Y^a \hconv M) \hconv Y^b \right) \simeq  \ang{i}^{k_ia} \hconv \left(  Y^a \hconv (M \hconv Y^b) \right) \\
&\simeq  (\ang{i}^{k_ia} \hconv   Y^a) \hconv (M \hconv Y^b)  \\
&\simeq  Z^a \hconv (M \hconv Y^b),
\end{aligned}
\end{equation}
where the fourth isomorphism follows from \eqref{Eq: YaM}. By the same argument, we have 
$$
Y^a \hconv (M \hconv Z^b) \simeq (Y^a \hconv M) \hconv Z^b.
$$

We now set $N :=(Y^a \hconv M) \hconv Y^b$. By \cite[Lemma 2.24]{KKOP23B}, we have 
\begin{align*}
\de_i( N) & = \frac{2}{(\al_i, \al_i)} \de(\ang{i}, N) = 
\frac{2}{(\al_i, \al_i)} \left( \de(\ang{i}, Y^a \hconv M) + \de (\ang{i}, Y^b) \right) \\
& = \frac{2}{(\al_i, \al_i)} \left( \de(\ang{i}, Y^a ) + \de (\ang{i}, M) \right) \\
& = \frac{2}{(\al_i, \al_i)} \de (\ang{i}, M) = \de_i( M) \ge k_i (a+b).
\end{align*}
Applying \cite[Theorem 2.17]{KKOP23A} to $N$, we have that $ \hd ( \ang{i}^{k_ia} \conv N \conv \ang{i}^{k_ib} )$ is simple. By Lemma \ref{Lem: normal} and \eqref{Eq: iY Yi}, we have 
\begin{align*}
\hd ( \ang{i}^{k_ia} \conv N \conv \ang{i}^{k_ib} ) &\simeq  
\left( \ang{i}^{k_ia} \hconv \left( Y^a \hconv ( M \hconv Y^b) \right) \right) \hconv \ang{i}^{k_ib}  \\
&\simeq  \left(  \left( \ang{i}^{k_ia} \hconv  Y^a \right) \hconv ( M \hconv Y^b)  \right) \hconv \ang{i}^{k_ib} \\
&\simeq   (  Z^a  \hconv ( M \hconv Y^b)   ) \hconv \ang{i}^{k_ib} \\
&\simeq   ( ( Z^a  \hconv  M) \hconv Y^b  )  \hconv \ang{i}^{k_ib} \\
&\simeq    ( Z^a  \hconv  M) \hconv ( Y^b    \hconv \ang{i}^{k_ib} ) \\
&\simeq    ( Z^a  \hconv  M) \hconv Z^b,   
\end{align*}
where the fourth isomorphism follows from \eqref{Eq: ZaMYb}. In the same manner, we have 
$$
\hd ( \ang{i}^{k_ia} \conv N \conv \ang{i}^{k_ib} ) \simeq 
  Z^a  \hconv ( M \hconv Z^b), 
$$
which completes the proof.
\end{proof}

\begin{lemma} \label{Lem: key homo}
Let  $L,M,N \in \Irr(R\gmod)$ and let  
$J \subset I$ such that 
$$
(\al_j, \al_{j'}) =0 \qquad \text{ for any $j,j' \in J$ with $j \ne j'$.} 
$$
We define $Z := \conv_{j\in J} \ang{j}^{ k_j}$ for some $k_j \in \Z_{\ge0}$, and set
$$
\zeta := -\wt (Z) ,  \qquad \beta := -\wt(L), \qquad \gamma := -\wt(M).
$$
Suppose that 
\bna
\item $(L,Z)$ and $(Z, M)$ are unmixed,
\item $\Res_{\beta, \gamma + 2\zeta} N \simeq L\tens \left( Z^{ 2} \hconv M  \right)$,  
\item $\Res_{\beta+ 2\zeta, \gamma } N \simeq \left(  L \hconv Z^{ 2}   \right) \tens M $.  
\ee
Then there is a nonzero homomorphism 
$$
\left(  L \hconv Z    \right) \tens \left(   Z \hconv M  \right) \longrightarrow \Res_{\beta+ \zeta, \gamma+ \zeta} N.
$$
\end{lemma}
\begin{proof} As $Z$ is affreal,  the heads of $L \conv Z^{  k} $ and $Z^{ k} \conv M$ are simple for any $k \in \Z_{>0}$. 
By (a) and (b), we have a nonzero homomorphism 
$$
L \tens Z  \tens Z  \tens M \longrightarrow \Res_{\beta,  \zeta,  \zeta, \gamma}N,
$$
which gives a nonzero homomorphism 
\begin{align*}
f: \left(L \conv Z  \right) \tens  \left( Z  \conv M \right) \longrightarrow \Res_{\beta +  \zeta,  \zeta + \gamma}N.
\end{align*}

Let $K_L$ and $K_M$ be the maximal submodules of $L\conv Z $ and $Z  \conv M$. 
We now  show that $ f( \left(L \conv Z  \right) \tens  K_M) =0$.
Since $Z^{  2}\conv M$ has a simple head, there is a surjection 
$$
Z  \conv \left( Z  \hconv M \right) \twoheadrightarrow Z^{ 2}\hconv M,
$$
which yields that there is a nonzero homomorphism 
$$
Z  \tens \left( Z  \hconv M \right) \longrightarrow 
\Res_{\zeta, \gamma+ \zeta} \left( Z^{ 2}\hconv M \right).
$$
Applying $L \tens -$ on both sides, we have 
\begin{align*}
L \tens Z  \tens \left( Z  \hconv M \right) \longrightarrow 
 L \tens \Res_{ \zeta, \gamma+ \zeta} \left( Z^{ 2 }\hconv M \right) \simeq \Res_{\beta,  \zeta, \gamma+ \zeta} \left( N \right),
\end{align*}
where the isomorphism on the right follows from (b). This gives a nonzero homomorphism 
$$
\left( L \conv Z  \right) \tens \left( Z  \hconv M \right) \longrightarrow 
  \Res_{\beta+  \zeta, \gamma+ \zeta} \left( N \right),
$$
which tells us that 
$$ 
f \left( \left(L \conv Z  \right) \tens  K_M \right) =0.
$$
By the same argument with (c), we have  
$f \left( K_L \tens    \left( Z  \conv M \right) \right) =0.$
Therefore, the homomorphism $f$ factors through $\left( L\hconv Z  \right) \tens \left( Z  \hconv M \right)$, i.e., 
$
 \left(L \hconv Z  \right) \tens  \left( Z  \hconv M \right) \buildrel f \over \longrightarrow \Res_{\beta +  \zeta, \gamma + \zeta}N.
$
\end{proof}

Recall the extremal vector $b_{u\La}$ of $B(\La)$ in Section \ref{Sec: HWC}. 
As mentioned at the beginning of this section, 
we regard $B(\La)$ as a subset of $B(\infty)$ via the embedding \eqref{Eq: iota BLa Binf}.
Then we have the main theorem of this section.

\begin{theorem} \label{Thm: main}
Let $\La \in \wlP_+$, $w,u\in \weyl$ and $J  \subset I$.
Suppose that 
\bna
\item \label{Thm: main a} 
$ (\al_j, \al_{j'})=0$ for any $j,j'\in J$ with $j\ne j'$,
\item \label{Thm: main b} 
$\supp ( u ) \subset \supp (w)$ and 
$w\La  + \La \in \sum_{k\in I \setminus \supp ( w) } \Z \La_k $, 
\item \label{Thm: main c}
$uw \La = \left(\prod_{j\in J}s_j \right) u \La$, 
\item \label{Thm: main d}
$\left(\prod_{j\in J}s_j \right) u \La$ is $u^{-1}$-dominant, 
\item \label{Thm: main e}
$2\langle h_j, \La \rangle = \langle h_j, u\La \rangle $ for any $j\in J$,
\item \label{Thm: main f}
$us_j < u$ for any $j\in J$, 
\item \label{Thm: main g}
if we set $b := \left(\prod_{j\in J} \te^{* \langle h_j, \La \rangle}_j \right)  b_{u\La} $, then $*b = b$.
\ee
Then, for any $x \in \weyl$ with $x\La = w\La$, the element  
$$
\left(\prod_{j\in J} \tf^{ \langle h_j, \La \rangle}_j \right)  b_{u\La}
\in B(\La) \subset B(\infty)
$$ 
is  $\cdual_x$-\stb imaginary. 
\end{theorem}

\begin{proof}
Recall generalized determinantial modules $\dM_w(\la, \mu)$ (see Section \ref{Sec: DM}).
We set 
$$
z := \prod_{j\in J} s_j.
$$ 
Since $\La$ is a dominant integral weight, $\La$ is $u$-dominant. Thus we have
\begin{equation} \label{Eq:*uLL}
\begin{aligned}
\phi_*(\dM(u\La, \La)) &= \phi_*(\dM_u(u\La, \La)) \simeq \dM_{u^{-1}}(-\La, -u\La) =\dM_{u^{-1}}(w\La, uw\La)\\
& =\dM_{u^{-1}}(u^{-1}(uw\La), uw\La) =\dM_{u^{-1}}(u^{-1}(zu\La), zu\La) \\
&= \dM(u^{-1}(zu\La), zu\La) = \dM(w\La, zu\La),
\end{aligned}
\end{equation}
where 
the first and sixth equalities follow from \ref{Thm: main d} and \cite[Remark 2.22 (i)]{KKOP23A},
the second follows from \eqref{Eq: Phi * dM},
the third follows from \ref{Thm: main b} and \eqref{Eq: Mw la la'}, 
the fifth and seventh follows from \ref{Thm: main c}.

We set $X := \dM(u\La, z\La)$, and 
define $ k_j := \langle h_j, \La \rangle \in \Z_{\ge0}$ for $j\in J$ and 
$$
Z :=  \dM(z\La, \La) = \conv_{j\in J} \ang{j}^{k_i}. 
$$
Since \ref{Thm: main g} says $\phi_*(\dM(u\La, z\La)) \simeq \dM(u\La, z\La)$,
 we have 
\begin{equation} \label{Eq: ZX MwLzuL}
\begin{aligned}
Z\hconv X & \simeq  \dM(z\La, \La) \hconv \dM(u\La, z \La) \simeq 
 (\phi_* \dM(z\La, \La) ) \hconv   ( \phi_*\dM(u\La, z \La ) ) \\ 
 &\simeq 
\phi_* \left( \dM(u\La, z \La )  \hconv \dM(z\La, \La) \right) \simeq \phi_* ( \dM(u\La,  \La ) ) \\
&\simeq \dM(w\La,  zu\La ),
\end{aligned}
\end{equation}
where 
the second follows from $ \phi_* (Z) = Z  $ by \ref{Thm: main a},
the third follows from  \eqref{Eq: phi *},
the fourth follows from \ref{Thm: main f} and \cite[Proposition 4.6]{KKOP18}, and
the fifth follows from \eqref{Eq:*uLL}.
Since $\langle h_j, u\La \rangle = 2k_j $ by \ref{Thm: main e}, we have 
\begin{align*}
\de_j(X) -2k_j &= \langle h_j, u \La - z \La \rangle + \ep_j(X) + \eps_j(X) - 2k_j \\ 
&= \langle h_j, u \La  \rangle - \langle h_j,  z \La \rangle - 2k_j \\ 
& = \langle h_j,   \La \rangle = k_j \ge 0,
\end{align*}
 where the second equality follows from \cite[Lemma 1.7 and Proposition 4.1]{KKOP18}.
By Lemma \ref{Lemma: unique head} and \eqref{Eq: ZX MwLzuL}, we have 
\begin{equation} \label{Eq: MZ ZM}
\begin{aligned} 
\dM(w\La, zu \La) \hconv Z & \simeq (Z \hconv X )\hconv Z  \simeq
Z \hconv (X \hconv Z) \\
& \simeq Z \hconv \dM(u \La, \La). 
\end{aligned}
\end{equation}

We set 
$$
L := \dM(w\La, zu \La), \qquad  M :=\dM(u \La, \La)  \qtq N := \dM(w\La, \La),
$$ 
and define $T:= L\hconv Z$. 
By \eqref{Eq: MZ ZM}, we have $T \simeq Z \hconv M$,
\begin{align} \label{Eq: key info}
T \simeq \cL \left( \left(\prod_{j\in J} \tf^{ \langle h_j, \La \rangle}_j \right)  b_{u\La} \right),\quad 
k_t = \eps_t(T) = \eps_t(N) \quad \text{ for $t\in J$},
\end{align}
and, by \ref{Thm: main e},
$$
\dM(w\La, u \La) \simeq L \hconv Z^{2} \qquad \text{and}\qquad 
\dM(zu \La, \La) \simeq Z^{2} \hconv M .
$$

We write $\zeta := -\wt(Z)$, $\beta := -\wt(L)$ and $\gamma:= - \wt(M)$. Since 
$(\dM(x\La', y\La'),  \dM(y\La', \La'))$ is unmixed whenever $x\ge y$ (see \cite[proof of Proposition 4.6]{KKOP18}), we have
\begin{align*}
\Res_{ \beta+ 2\zeta, \gamma }N &\simeq  \dM(w\La, u \La) \tens \dM(u\La, \La)
 = (L \hconv Z^{2}) \tens M, \\
 \Res_{ \beta , \gamma + 2\zeta }N &\simeq  \dM(w\La, zu \La) \tens \dM(zu\La, \La)
 = L  \tens   ( Z^{2} \hconv M).
\end{align*}
Lemma \ref{Lem: key homo} tells us that there is a nonzero homomophism 
\begin{align} \label{Eq: TT N}
 T\tens T \simeq  (L \hconv Z) \tens ( Z \hconv M) \longrightarrow \Res_{\beta+ \zeta , \gamma + \zeta} N.
\end{align}
Combining Lemma \ref{Lem: criterion for imaginary} with \eqref{Eq: key info}, we conclude that 
$T$ is imaginary. 

Let  $x \in \weyl$ with $x\La = w\La$. 
Then $N = \dM(x \La, \La)$ is frozen in the category $\cC_{x}$.
Moreover, \eqref{Eq: TT N} says that there is a surjective homomorphism 
$$
T \conv T \twoheadrightarrow N.
$$ 
Therefore, by \cite[Lemma 5.6]{KKOP21A}, we have $T \eqfr \dual_x(T)$.
\end{proof}

The following corollary gives a more explicit condition to find \stb imaginary elements in terms of Weyl group elements and weights, which will play a crucial role in Section  \ref{Sec: PI for finite types}.

\begin{corollary}  \label{cor: condition for PI}
Let $v \in \weyl$, $J \subset I$, and $\La := \sum_{j\in J} k_j \La_j \in \wlP_+$, 
and set 
$$
z := \prod_{j\in J}s_j \in \weyl.
$$ 
Suppose that 
\bna
\item \label{cor: condition for PI a} 
$(\al_j, \al_{j'})=0  $ for $j,j' \in J$ with $j \ne j'$,
\item \label{cor: condition for PI b}
$vz \La - \sum_{j\in J} k_j \al_j \in \sum_{j \in I \setminus (\supp(v) \cup J) } \Z \La_j $, 
\item \label{cor: condition for PI c}
$v=v^{-1}$,
\item \label{cor: condition for PI d}
$\ell(z v z v z) = 2\ell(v) + 3|J|$, 
\item \label{cor: condition for PI e}
$\supp(v) \subset I \setminus J$. 
\ee
Then we have the following.
\bnum
\item   \label{cor: condition for PI 1}
For any $x \in \weyl$ with $ x\La = z v z v z  \La$, the element  
$$
b_{(v,J,\La)} := \left(\prod_{j\in J} \tf^{ \langle h_j, \La \rangle}_j \right)  b_{vz\La}  
$$ 
is $\cdual_x$-\stb imaginary. This element $b_{(v,J,\La)} $ is called the \stb imaginary element 
associated with the triple $(v,J,\La)$. 
\item   \label{cor: condition for PI 2}
Let $\bfi = (i_1, i_2, \ldots, i_l) \in \Seq(zvz)$ such that $J=\{ i_1, i_2, \ldots, i_{|J|} \}$. Then we have 
$$
\STR_\bfi(b_{(v,J,\La)}) = (c_k)_{k\in [1,l]},
$$
where 
$
c_k = 
\begin{cases}
\ang{h_k, \La} & \text{if } 1 \le k \le |J|,\\
\ang{ h_{i_k}, s_{i_{k+1}} s_{i_{k+2}} \cdots s_{i_{l}} \La  }  & \text{if } |J| < k \le l,
\end{cases}  \quad \text{for $k\in [1,l]$.}
$ 

\ee
\end{corollary}

\begin{proof}
Let $u := vz$ and $w := u^{-1} z u =  z v z v z$. 

(i)
It is easy to show that all conditions except \ref{Thm: main b}, \ref{Thm: main e}, \ref{Thm: main g} of Theorem \ref{Thm: main} hold.

We write  
\begin{align} \label{Eq: uLala}
u \La - \sum_{j\in J} k_j \al_j = \la
\end{align} 
for some $\la \in \sum_{j \in I \setminus (\supp(v) \cup J) } \Z \La_j$. Then we have 
$$
\ang{h_j, u\La} = \ang{h_j, \sum_{j\in J} k_j \al_j + \la} = 2k_j = 2\ang{h_j, \La}
\quad \text{ for $j\in J$,}
$$
and 
\begin{align*}
w\La  = u^{-1} z (u \La) = u^{-1} z \left( \sum_{j\in J} k_j \al_j  + \la\right)
= u^{-1} \left( -\sum_{j\in J} k_j \al_j  + \la \right) = -\La + 2\la,
\end{align*}
which implies \ref{Thm: main b} and \ref{Thm: main e} of Theorem \ref{Thm: main}.

Since $z\La$ is $uz$-dominant and $z u\La$ is $z u^{-1}$-dominant, we have
\begin{equation} \label{Eq: suLzLzu} 
\begin{aligned} 
\phi_*(\dM(u\La, z\La)) &= \phi_*(\dM(v (z\La), z\La)) = \phi_*( \dM_{v}( v (z\La), z\La)) \\
&\simeq  \dM_{v}(-z\La, -v (z\La)) = \dM_{v}(zw\La, vzw\La) \\
&=  \dM_{v}(v (zu\La), zu\La) =  \dM(v (zu\La), zu\La), 
\end{aligned}
\end{equation}
where the second and sixth equalities follow from \cite[Remark 2.22 (i)]{KKOP23A},
the third  follows from \eqref{Eq: Phi * dM} and \ref{cor: condition for PI c},
the fourth  follows from \eqref{Eq: Mw la la'} and \ref{Thm: main b} of Theorem \ref{Thm: main}.

We write a reduced expression $ v = s_{t_1} s_{t_2} \cdots s_{t_p} $.
 From \eqref{Eq: uLala}, we have 
$$
z\La = \La - \sum_{j\in J} k_j \al_j \qtq zu\La =  - \sum_{j\in J} k_j \al_j + \la
$$
and $ \langle h_{t_a}, \La \rangle = 0$ for $a\in [1,p]$ by \ref{cor: condition for PI e}. Thus we have  
\begin{align*}
\langle h_{t_a}, s_{t_{a+1}} \cdots s_{t_{p}}  z\La \rangle = 
 \langle h_{t_a}, s_{t_{a+1}} \cdots s_{t_{p}}  (- \sum_{j\in J} k_j \al_j ) \rangle = 
\langle h_{t_a}, s_{t_{a+1}} \cdots s_{t_{p}} zu\La \rangle  
\end{align*}
for any $a\in [1,p]$. Then \cite[Lemma 1.7]{KKOP18} and \eqref{Eq: suLzLzu} imply 
$$
\dM(u\La, z\La) \simeq \dM(v(z\La), z\La) \simeq \dM(v(zu\La), zu\La) \simeq \phi_*(\dM(u\La, z\La)),
$$
which gives \ref{Thm: main g} of Theorem \ref{Thm: main}.

(ii) The assertion follows from the fact that $u_{vz\La}$ is an extremal vector and $(\al_j, \al_{j'})=0 
$ for $j,j' \in J$ with $j \ne j'$.
\end{proof}

\begin{remark} \label{rmk: stb im G2}
There exists a \stb imaginary element which does not satisfy the conditions given in Theorem \ref{Thm: main} and Corollary \ref{cor: condition for PI}.
Let $A_q(\n(w))$ be of type $G_2$ (see Proposition \ref{Prop: triX in ft} \ref{P9}), and let $\bfi = ( 1, 2, 1, 2, 1,2) \in \Seq(w_0)$. We choose $b \in B(\infty)$ such that $\PBW_\bfi (b) = (1, 1, 0, 0, 1, 1)$. Note that $b = \tf_1 \tf_2^4 \tf_1^2 \tf_2 \one $. Then one can check that $b$ is a $\cdual_{w_0}$-\stb imaginary, but it does not satisfy the conditions of Theorem \ref{Thm: main}.
\end{remark}

\subsection{Similarity of crystals} \label{Sec: similarity}
Let $\cmC = (c_{i,j})_{i,j\in I}$ (resp.\ $\cmC' = (c_{i,j}')_{i,j\in I'}$) be a generalized Cartan matrix with an index set $I$ (resp.\ $I'$), and 
let $\weyl_\cmC$, $\wlP_\cmC$ and $B_\cmC$ (resp.\ $\weyl_{\cmC'}$ and $\wlP_{\cmC'}$ and $B_{\cmC'}$) denote the corresponding Weyl group, weight lattices and crystals. 
We choose 
\begin{align} \label{Eq: mi Ki}
 K_i \subset I' \quad  \text{for $i\in I$}, \qquad  m_{i,k} \in \Z_{\ge0}, \quad \text{ for $k\in K_i$}
\end{align}
and define 
$$
K := \bigcup_{i\in I} K_i \subset I', \qquad \wlP_{I'\setminus K} := \sum_{j\in I' \setminus K} \Z \La_j
$$ 
and  
$$
\theta : \wlP_\cmC \rightarrow \wlP_{\cmC'} \qquad  \text{ by $\La_i \mapsto  \sum_{k\in K_i} m_{i,k}\La_k $ \quad for $i\in I$.} 
$$
Suppose that 
\begin{enumerate}[label=(C\arabic*), ref=C\arabic*]
\item \label{cond:C1}
$K_i \cap K_j = \emptyset$ whenever $i, j\in I$ are distinct,
\item \label{cond:C2} for each $i\in I$, $(\al_k, \al_{k'})=0$ whenever $k, k'\in K_i$ are distinct, 
\item \label{cond:C3} for $i,j\in I$, $(\al_i, \al_j)=0$ if and only if $ (\al_k, \al_{k'})=0$ for any $k\in K_i$ and $k' \in K_j$,
\item \label{cond:C4} for any $i\in I$,
$\theta( \al_i ) \equiv \sum_{k\in K_i} m_{i,k}\al_k  \mod \wlP_{I'\setminus K} $,
\item \label{cond:C5} the map $\Theta: \weyl_\cmC \rightarrow \weyl_{\cmC'} $ defined by 
$$
\Theta (s_i) = \prod_{k\in K_i} s_k\qquad \text{ for $i\in I$}
$$
is an injective group homomorphism,
\item \label{cond:C6} for any $\La \in \wlP_\cmC^+$, there is an injective crystal morphism $ \Upsilon: B_{\cmC}(\La) \rightarrow B_{\cmC'}(\theta(\La))$ such that $ \Upsilon(b_{\La}) = b_{\theta(\La)}$ and 
$$
\Upsilon(\tf_i(b)) = \prod_{k\in K_i} \tf_k^{m_k} (\Upsilon(b)), \qquad 
\Upsilon(\te_i(b)) = \prod_{k\in K_i} \te_k^{m_k} (\Upsilon(b))
$$
for any $i\in I$ and $b \in B_{\cmC}(\La)$.
\end{enumerate}
It follows from \eqref{cond:C2}, \eqref{cond:C4} and \eqref{cond:C5} that, for any $i\in I$ and $\la\in \wlP_\cmC$, 
$$
\theta( s_i \la) \equiv \Theta(s_i) \theta(\la) 
\qquad  \mod \wlP_{I'\setminus K}. 
$$
Then one can prove the following easily.

\begin{corollary} \label{Cor: folding triple}
Let us consider $K_i$ and $m_{i,k}$ given in \eqref{Eq: mi Ki} satisfying the conditions \eqref{cond:C1}--\eqref{cond:C6}. Let $v \in \weyl$,  $J \subset I$ and $\La = \sum_{j\in J} k_j \La_j \in \wlP_+$. 
We set $K_J := \cup_{j\in J} K_j$.
Then the following are equivalent:
\bna
\item  the triple $(v,J, \La)$ satisfies all conditions \ref{cor: condition for PI a}--\ref{cor: condition for PI e} in Corollary \ref{cor: condition for PI},
\item  the triple $(\Theta(v), K_J, \theta(\La))$ satisfies all conditions \ref{cor: condition for PI a}--\ref{cor: condition for PI e} in Corollary \ref{cor: condition for PI}.
\ee
\end{corollary}

\begin{example} \label{Ex: folding 1} 

Let $U_q(\cmC)$ be the quantum group of a generalized Cartan matrix $\cmC = (c_{i,j})_{i,j\in I}$ with an index set $I$, and let $H $ be a subset of $I$. We write $\cmC_H := \cmC|_{H \times H}$ and denote by $U_H$ the subalgebra of $U_q(\cmC)$ generated by $f_i$, $e_i$, $q^h$ for $i\in H$ and  $h\in \wlP_{\cmC_H}$, where $ \wlP_{\cmC_H}$ is regarded as a subset of $ \wlP_{\cmC}$.
 In this case, $\theta: \wlP_{\cmC_H} \to \wlP_{\cmC}$ is given by $\La_i \mapsto \La_i$ for $i\in H$ and $U_q(\cmC_H)$ is isomorphic to $U_H$, i.e., there exists the natural embedding 
$ U_q(\cmC_H) \buildrel \sim \over \longrightarrow U_H \subset U_q(\cmC)$. This induces
 the crystal embedding
\begin{align} \label{Eq: embedding H I} 
\emb_{ I, H}: B_{H}(\infty) \rightarrowtail B(\infty),
\end{align} 
where $B_{H}(\infty)$ (resp.\ $B(\infty)$) is the infinite crystal of $U_q^-(\cmC_H)$ (resp.\ $U_q^-(\cmC)$), and $B(\infty)$ is regarded as a $U_q(\cmC_H)$-crystal.
If we set  $K_i := \{ i \}$ and  $m_{i,k} :=1$  for all $i\in I$ and $k\in K_i$,
then the conditions \eqref{cond:C1}--\eqref{cond:C6} hold in this setting.

\end{example}

\begin{example} \label{Ex: folding 2} 
The \emph{folding} of crystals introduced in \cite{Kash95} gives an example. 
Let $\sigma$ be one of the following Dynkin diagram automorphisms:  
\begin{equation} \label{Eq: Folding Dynkin}
\begin{aligned} 
	\scriptsize
& \xymatrix@R=0.1em@C=0.6em{
 &\circ \ar@{-}[r] \ar@{.}[dd] & \cdots \ar@{-}[r]  & \circ \ar@{-}[r] \ar@{.}[dd] & \circ \ar@{-}[dr]  \ar@{.}[dd] \\
 A_{2n-1} : & &&&& \circ  \ar@{.>}[ddd] \\
&  \circ \ar@{-}[r] \ar@{.>}[dd] & \cdots \ar@{-}[r] & \circ \ar@{-}[r] \ar@{.>}[dd] & \circ \ar@{-}[ru] \ar@{.>}[dd]  \\
&    &&&& \\
B_{n} :&  \circ \ar@{-}[r] & \cdots \ar@{-}[r] & \circ \ar@{-}[r] & \circ \ar@{=>}[r] & \circ 
}
\quad  \qquad
\xymatrix@R=0.1em@C=0.6em{
	&  \circ \ar@{-}[r]  \ar@{.}[dd]  & \circ  \ar@{-}[rd] \ar@{.}[dd] &&   \\
	E_6 : &   && \circ \ar@{-}[r] \ar@{.>}[ddd] & \circ \ar@{.>}[ddd]  \\
	& \circ \ar@{-}[r] \ar@{.>}[dd] & \circ \ar@{-}[ur] \ar@{.>}[dd] &&    \\
	&    &&&& \\
	F_4 :&  \circ \ar@{-}[r] & \circ \ar@{=>}[r] & \circ \ar@{-}[r] & \circ  
}
\\
& \xymatrix@R=0.1em@C=0.6em{
	& &&&& \circ  \ar@{.}[dd] \\
	D_{n+1} : & \circ \ar@{-}[r] \ar@{.>}[ddd] & \cdots \ar@{-}[r] & \circ  \ar@{-}[r] \ar@{.>}[ddd] & \circ  \ar@{-}[ru] \ar@{-}[rd]  \ar@{.>}[ddd] &  \\
	& &&&& \circ  \ar@{.>}[dd]   \\
	&    &&&& \\
	C_{n} :&  \circ \ar@{-}[r] & \cdots \ar@{-}[r] & \circ \ar@{-}[r] & \circ \ar@{<=}[r] & \circ  
}
\qquad \qquad 
\xymatrix@R=0.1em@C=0.6em{
	 &\circ  \ar@{-}[rd]  \ar@{.}[d]   \\
	D_{4} : &  \circ \ar@{-}[r] \ar@{.}[d] & \circ \ar@{.>}[ddd]  \\
	 & \circ  \ar@{-}[ru] \ar@{.>}[dd]    \\
	&    &&&& \\
	G_{2} :&  \circ    \ar@3{->}[r] & \circ 
}
\end{aligned}	 	
\end{equation}
and 
\begin{equation}  \label{Eq: Folding Dynkin 2}
\begin{aligned}  
	\small
& \xymatrix@R=0.1em@C=0.6em{
	&\circ \ar@{-}[r] \ar@{.}[dd] & \cdots \ar@{-}[r]  & \circ \ar@{-}[r] \ar@{.}[dd] & \circ \ar@{-}[dr]  \ar@{.}[dd] \\
	A_{2n-1} : & &&&& \bullet  \ar@{.>}[ddd] \\
	&  \circ \ar@{-}[r] \ar@{.>}[dd] & \cdots \ar@{-}[r] & \circ \ar@{-}[r] \ar@{.>}[dd] & \circ \ar@{-}[ru] \ar@{.>}[dd]  \\
	&    &&&& \\
	C_{n} :&  \circ \ar@{-}[r] & \cdots \ar@{-}[r] & \circ \ar@{-}[r] & \circ \ar@{<=}[r] & \circ 
}
\quad  \qquad
\xymatrix@R=0.1em@C=0.6em{
	&\circ  \ar@{-}[r]  \ar@{.}[d] & \circ \ar@{-}[ld] \ar@{.}[d]    \\
	A_{5} : &  \bullet \ar@{-}[r] \ar@{.}[d] & \circ \ar@{.>}[ddd]  \\
	& \circ  \ar@{-}[ru] \ar@{.>}[dd]    \\
	&    &&&& \\
	G_{2} :&  \circ    \ar@3{->}[r] & \circ 
}
\end{aligned}	
\end{equation}
Let $\cmC$ and $\cmC'$ be the corresponding folded and unfolded Cartan matrices respectively, and $\pi$ is the projection between the index sets induced from $\sigma$.
Then we set
$$
K_i :=   \pi^{-1}(i) \qtq 
m_{i,k} :=  
\begin{cases}
	1 & \text{ if $k$ is at $\circ$},\\
	2 & \text{ if $k$ is at $\bullet$}.
\end{cases}
$$ 
Then conditions \eqref{cond:C1}--\eqref{cond:C6} hold in this setting.
\end{example}

\vskip 2em

\section{Imaginary elements for finite types} \label{Sec: PI for finite types}

In this section, we construct explicit imaginary elements in finite type using Corollary \ref{cor: condition for PI} with the aid of SageMath (\cite{SageMath}). For simplicity, throughout this section, a \stb imaginary element means a $\cdual_w$-\stb imaginary element for some $w \in \weyl$.

For $\bfi = (i_1, \ldots, i_l)$, define 
$$
s_\bfi := s_{i_1} \cdots s_{i_l} \in \weyl. 
$$
Recall the crystal elements $b^{[k]}$ given in Definition \ref{def: k PBW} for $b\in B(\infty)$ and $k\in \Z_{>0}$.

\begin{prop} \label{Prop: triX in ft} 
Let $\cmC$ be a Cartan matrix of finite type, and  
let $\triX = (v, J, \La)$ be one of the triples listed in  \ref{P1}--\ref{P9} below.  
Then the triple $\triX$ satisfies all conditions \ref{cor: condition for PI a}--\ref{cor: condition for PI e} in Corollary \ref{cor: condition for PI}. Thus the crystal elements $b_\triX^{[k]} \in B(\infty)$ are \stb imaginary for all $k\in \Z_{>0}$, where $b_\triX$ is the \stb imaginary element associated with $\triX$ in Corollary \ref{cor: condition for PI}. 
\bnum 
\item[\namedlabel{P1}{{\rm(P1)}}]    
Let $\cmC$ be the Cartan matrix of the Dynkin diagram of type $A_n$ ($n\ge5$). 
\[
\xymatrix@C=1.5em@R=0.1em{
\circ \ar@{-}[r] &
\circ \ar@{-}[r] &
\cdots  \ar@{-}[r] &
\circ \ar@{-}[r] &
\circ \\
1 & 2 && n-1 & n
}
\]

\bna
\item \label{P1: a} 
Suppose that $n=5$. Define $ \triAs := (v, J, \La)$ by 
$$
v = s_2 s_4, \qquad J = \{1, 3,5\}, \qquad \La = \La_1 + 2\La_3 + \La_{5}.
$$
\item \label{P1: b}
Let $n \ge 5$. Define $ \triA_{n} := (v, J, \La)$ by 
$$
v = s_1 s_\bfi s_n, \qquad J = \{2, n-1\}, \qquad \La = \La_2 + \La_{n-1},
$$
where $\bfi = ( 3,4, \ldots, n-3, n-2 , n-3, \ldots, 4,3) \in I^{2n-9}$.
\ee

\item[\namedlabel{P2}{{\rm(P2)}}]  
Let $\cmC$ be the Cartan matrix of the Dynkin diagram of type $B_n$ ($n\ge3$).
\[
\xymatrix@C=1.5em@R=0.1em{
\circ \ar@{-}[r] &
\circ \ar@{-}[r] &
\cdots \ar@{-}[r] &
\circ \ar@{=>}[r] &
\circ \\
1 & 2 && n-1 & n
}
\]

\bna
\item \label{P2: a}  Suppose that $n=3$. Define $ \triBs := (v, J, \La)$ by 
$$
v = s_2, \qquad J = \{1, 3\}, \qquad \La = \La_1 + 2\La_3.
$$
\item \label{P2: b}  Let $n \ge 3$. Define $ \triB_{n} := (v, J, \La)$ by
$$
v = s_1 s_\bfi, \qquad J = \{2\}, \qquad \La = \La_2,
$$
where $\bfi = ( 3,4, \ldots, n-1, n , n-1, \ldots, 4,3) \in I^{2n-5}$.
\ee

\item[\namedlabel{P3}{{\rm(P3)}}]  
Let $\cmC$ be the Cartan matrix of the Dynkin diagram of type $C_n$ ($n\ge3$).
\[
\xymatrix@C=1.5em@R=0.1em{
\circ \ar@{-}[r] &
\circ \ar@{-}[r] &
\cdots \ar@{-}[r] &
\circ \ar@{<=}[r] &
\circ \\
1 & 2 && n-1 & n
}
\]

\bna
\item \label{P3: a}  Suppose that $n=3$. Define $ \triCs := (v, J, \La)$ by
$$
v = s_2, \qquad J = \{1, 3\}, \qquad \La = \La_1 +  \La_3.
$$
\item \label{P3: b}  Let $n \ge 3$. Define $ \triC_{n} := (v, J, \La)$ by
$$
v = s_1 s_\bfi, \qquad J = \{2\}, \qquad \La = \La_2,
$$
where $\bfi = ( 3,4, \ldots, n-1, n , n-1, \ldots, 4,3) \in I^{2n-5}$.
\ee

\item[\namedlabel{P4}{{\rm(P4)}}]  
Let $\cmC$ be the Cartan matrix of the Dynkin diagram of type $D_n$ ($n\ge4$).
\[
\xymatrix@C=1.5em@R=0.1em{
&&& \hskip 1em     \circ \hskip 0.5em n \ar@{-}[ddd]  & \\ \\ \\ 
\circ \ar@{-}[r] &
\circ \ar@{-}[r] &
\cdots \ar@{-}[r] &
\circ \ar@{-}[r] &
\circ \\
 \ \ 1 \ \ & 2 && n-2 & n-1
}
\]
\bna
\item \label{P4: a}  Suppose that $n=4$. Define $ \triDs := (v, J, \La)$ by 
$$
v = s_2, \qquad J = \{1, 3, 4\}, \qquad \La = \La_1 +  \La_3 +  \La_4.
$$
\item \label{P4: b}  Let $n \ge 4$. Define $ \triD_{n} := (v, J, \La)$ by
$$
v = s_1 s_\bfi, \qquad J = \{2\}, \qquad \La = \La_2,
$$
where $\bfi = (3,4)$ if $n=4$ and 
$\bfi = ( 3,4, \ldots, n-3, n-2, n-1, n , n-2, n-3, \ldots, 4,3) \in I^{2n-6}$ if $n>4$.
\ee

\item[\namedlabel{P5}{{\rm(P5)}}]  
Let $\cmC$ be the Cartan matrix of the Dynkin diagram of type $E_6$.
\[
\xymatrix@C=1.5em@R=0.1em{
&& \hskip 1em     \circ \hskip 0.5em 2 \ar@{-}[ddd]  & \\ \\ \\ 
\circ \ar@{-}[r] &
\circ \ar@{-}[r] &
\circ \ar@{-}[r] &
\circ \ar@{-}[r] &
\circ \\
 \ \ 1 \ \ & 3 & 4 & 5 &  \ \ 6 \ \ 
}
\]
\bna
\item \label{P5: a}   Define $ \triE_{6, \bigstar} := (v, J, \La)$ by
$$
v = s_5 s_4 s_3 s_2 s_4 s_5 s_2 s_4 s_3 , \qquad J = \{1,6\}, \qquad \La = \La_1+\La_6.
$$
\item \label{P5: b}  Define  $ \triE_{6, \clubsuit} := (v, J, \La)$ by 
$$
v = s_4 s_5 s_6 s_3 s_4 s_5 s_1 s_3 s_4, \qquad J = \{2\}, \qquad \La = \La_2.
$$
\ee

\item[\namedlabel{P6}{{\rm(P6)}}]  
Let $\cmC$ be the Cartan matrix of the Dynkin diagram of type $E_7$.
\[
\xymatrix@C=1.5em@R=0.1em{
&& \hskip 1em     \circ \hskip 0.5em 2 \ar@{-}[ddd]  & \\ \\ \\ 
\circ \ar@{-}[r] &
\circ \ar@{-}[r] &
\circ \ar@{-}[r] &
\circ \ar@{-}[r] &
\circ \ar@{-}[r] &
\circ \\
 \ \ 1 \ \ & 3 & 4 & 5 &  \ \ 6 \ \ & \  7  \
}
\]
Define $ \triE_{7, \bigstar} := (v, J, \La)$ by
$$
v = s_3 s_4 s_5 s_6 s_7 s_2 s_4 s_5 s_6 s_3 s_4 s_5 s_2 s_4 s_3 , \qquad J = \{1\}, \qquad \La = \La_1.
$$

\item[\namedlabel{P7}{{\rm(P7)}}] 
Let $\cmC$ be the Cartan matrix of the Dynkin diagram of type $E_8$.
\[
\xymatrix@C=1.5em@R=0.1em{
&& \hskip 1em     \circ \hskip 0.5em 2 \ar@{-}[ddd]  & \\ \\ \\ 
\circ \ar@{-}[r] &
\circ \ar@{-}[r] &
\circ \ar@{-}[r] &
\circ \ar@{-}[r] &
\circ \ar@{-}[r] &
\circ \ar@{-}[r] &
\circ \\
 \ \ 1 \ \ & 3 & 4 & 5 &  \ \ 6 \ \ & \  7   & \  8  
}
\]
Define $ \triE_{8, \bigstar} := (v, J, \La)$ by
\begin{align*}
v &= s_7 s_6 s_5 s_4 s_2 s_3 s_1 s_4 s_3 s_5 s_4 s_2 s_6 s_5 s_4 s_3 s_1 s_7 s_6 s_5 s_4 s_2 s_3 s_4 s_5 s_6 s_7, \\
\qquad J &= \{8\}, \qquad \La = \La_8.
\end{align*}

\item[\namedlabel{P8}{{\rm(P8)}}]  
Let $\cmC$ be the Cartan matrix of the Dynkin diagram of type $F_4$.
\[
\xymatrix@C=1.5em@R=0.1em{
\circ \ar@{-}[r] &
\circ \ar@{=>}[r] &
\circ \ar@{-}[r] &
\circ \\
1 & 2 &3& 4 
}
\]

\bna
\item \label{P8: a}   Define $ \triF_{ \bigstar} := (v, J, \La)$ by
$$
v = s_2 s_3 s_4 s_2 s_3 s_2 , \qquad J = \{1\}, \qquad \La = \La_1.
$$
\item \label{P8: b}  Define $ \triF_{ \clubsuit} := (v, J, \La)$ by
$$
v = s_3 s_2 s_3 s_1 s_2 s_3, \qquad J = \{4\}, \qquad \La = \La_4.
$$
\ee

\item[\namedlabel{P9}{{\rm(P9)}}]  
Let $\cmC$ be the Cartan matrix of the Dynkin diagram of type $G_2$.
\[
\xymatrix@C=1.5em@R=0.1em{
\circ \ar@3{->}[r] &
\circ \\
1 & 2  
}
\]

\bna
\item \label{P9: a}   Define $ \triG_{ \clubsuit} := (v, J, \La)$ by
$$
v = s_2, \qquad J = \{1\}, \qquad \La = \La_1.
$$

\item \label{P9: b}  Define $ \triG_{\bigstar } := (v, J, \La)$ by
$$
v = s_1, \qquad J = \{2\}, \qquad \La = \La_2.
$$
\ee
\ee
\end{prop}
\begin{proof}
By Corollary \ref{Cor: folding triple} together with the folding of crystals given in Example \ref{Ex: folding 2}, we have 
\bni
\item if $\triA_{\clubsuit}$ satisfies all the conditions in Corollary \ref{cor: condition for PI}, then so does $\triX_{\clubsuit}$ for each $\triX = \triB$, $\triC$, $\triD$ and $\triG$,
\item if $\triA_{n}$ satisfies all the conditions in Corollary \ref{cor: condition for PI} for any $n\in \Z_{\ge 5}$, 
then so do $\triG_{\bigstar}$, $\triD_{t}$  for $t\ge 4$, and $\triX_{m}$ for $\triX = \triB, \triC$ and $m\ge 3$,
\item if $\triE_{6, \bigstar}$ and $\triE_{6, \clubsuit}$ satisfy all the conditions in Corollary \ref{cor: condition for PI}, then so do $\triF_{ \bigstar}$ and $\triF_{\clubsuit}$.
\ee
Thus it suffices to check the cases for types $A$ and $E$. In this case, the conditions  \ref{cor: condition for PI a} and \ref{cor: condition for PI e} of Corollary \ref{cor: condition for PI}  are obvious.  

We consider the case $\triA_n$. It is straightforward to verify that $v=v^{-1}$ and that $zvzvz$ is reduced. Moreover,
\begin{align*}
vz (\La_2 + \La_{n-1}) &= s_1 s_{\bfi} s_n s_2 s_{n-1} (\La_2 + \La_{n-1}) \\
&= \La_2 + \La_{n-1} - \al_1 - \al_2 - \left( \sum_{k\in [3,n-2]} 2 \al_k \right) - \al_{n-1} - \al_n \\
&= \al_2 + \al_{n-1}.
\end{align*}
Thus, conditions \ref{cor: condition for PI b}, \ref{cor: condition for PI c}, and \ref{cor: condition for PI d} of Corollary \ref{cor: condition for PI} are satisfied.

The remaining cases $\triA_{\clubsuit}$ and $\triE$ can be checked directly in a similar manner. 
\end{proof}

Combining Proposition \ref{Prop: triX in ft} with the embedding given in Example \ref{Ex: folding 1}, we have the following theorem.

\begin{theorem}
Let $U_q(\cmC)$ be the quantum group associated with a finite-type Cartan matrix $\cmC$ with index set $I$.
For each subset $H\subseteq I$, let $U_q(\cmC_H)$ be the quantum group associated with the submatrix $\cmC_H=\cmC|_{H\times H}$, and let $\mathfrak{I}_H$ denote the set of elements $b_\triX^{[k]}$ arising from all triples $\triX$ given in Proposition \ref{Prop: triX in ft} for $U_q(\cmC_H)$ and all $k\in\Z_{>0}$.
Then every element of
$$
\mathfrak{I}:=\bigcup_{H\subseteq I}\emb_{I,H}(\mathfrak{I}_H)\subseteq B(\infty)
$$
is \stb imaginary, where $\emb_{I,H}:B_H(\infty)\rightarrowtail B(\infty)$ is the natural crystal embedding defined in \eqref{Eq: embedding H I}.
\end{theorem}

\vskip 2em 

\section{Leclerc's examples} \label{Sec: LE}

In \cite{Lec03}, Leclerc constructed examples of imaginary elements in finite types. In this section, we will see that Leclerc's examples are special cases of the \stb imaginary elements given in Proposition \ref{Prop: triX in ft}.

Recall Lusztig's braid symmetries  $\lT_i $ and $ \lTs_i$ in \eqref{Eq: braid} and the positive roots $\beta_k$ in \eqref{Eq: betak}.  Let $\bfi = (i_1, i_2, \ldots, i_\ell) \in \Seq(w)$ be a reduced expression of $w \in \weyl$ and let $\epsilon = \pm 1$. For any $c \in \Z_{ \ge 0}$ and $k\in [1,\ell]$,  define  
\begin{align*} 
\rv_{\bfi, \epsilon} (\beta_k) ^{(c)} := 
\begin{cases}
	\lT_{i_1}\lT_{i_2} \cdots \lT_{i_{k-1}}(f_{i_k}^{(c)})  	 & \text{ if } \epsilon = -1, \\
	\lTs_{i_1}\lTs_{i_2} \cdots \lTs_{i_{k-1}}(f_{i_k}^{(c)})   & \text{ if } \epsilon = 1.
\end{cases}
\end{align*} 
For any $\bfa=(a_k)_{k\in [1,\ell]} \in \Z_{\ge 0}^{\oplus \ell}$, we define 
\begin{align*} 
\rv_{\bfi, \epsilon} (\bfa) := 
\begin{cases}
	\rv_{\bfi, \epsilon} (\beta_\ell) ^{(a_\ell)} \rv_{\bfi, \epsilon} (\beta_{\ell-1}) ^{(a_{\ell-1})} \cdots 	\rv_{\bfi, \epsilon} (\beta_{1}) ^{(a_{1})}  & \text{ if } \epsilon = -1, \\
	\rv_{\bfi, \epsilon} (\beta_{1}) ^{(a_{1})}  \cdots \rv_{\bfi, \epsilon} (\beta_{\ell-1}) ^{(a_{\ell-1})}   \rv_{\bfi, \epsilon} (\beta_\ell) ^{(a_\ell)}   & \text{ if } \epsilon = 1.
\end{cases}
\end{align*}
The monomials $ \rv_{\bfi, \epsilon} (\bfa)  $ are called the \emph{PBW monomials} with respect to $\bfi$,  $\epsilon$ and $\bfa$.
We also have the dual $ \rv^{\up}_{\bfi, \epsilon} (\bfa) $ with respect to the Kashiwara bilinear form (\cite[Proposition 3.4.4]{K91}).
For any $b\in B(w)$, there exists $\bfa \in \Z_{\ge 0}^{\oplus \ell}$ such that  
\begin{align} \label{Eq: Gup PBW e}
\Gup(b) \equiv \rv^{\up}_{\bfi, \epsilon} (\bfa) \mod qL(\infty).  
\end{align}
In this case we write  $\PBW_{\bfi, \epsilon}(b) := \bfa $.
Note that  
$\rv^{\up}_{\bfi} (\bfa) = \rv^{\up}_{\bfi, -1} (\bfa) $ and $\PBW_{\bfi}(b) = \PBW_{\bfi, -1}(b)$,
where  $\rv^{\up}_{\bfi} (\bfa)$ and $\PBW_{\bfi}(b)$ are given in Section \ref{Sec: PBW}.
Then we have the following lemma (see \cite[2.11]{L90} and see also \cite[Lemma 2.2]{P23}).
\begin{lemma}    \label{Lem: two ld}	
Assume that $U_q(\g)$ is of finite type. Let $b \in B(\infty)$, $ \bfi = (i_1, i_2,  \ldots, i_\ell) \in \Seq(w_0) $, and set 
$\bfi^\vee := (i_\ell^*, \ldots, i_2^*,  i_1^*)$, where $\al_{i^*} = - w_0(\al_i)$.
If we write 
$$
\PBW_{\bfi, -1}(b) = (a_k)_{k\in [1, \ell]} \qtq \PBW_{\bfi^\vee, 1}(b) = (a_k')_{k\in [1, \ell]},
$$ 
then we have $ a_k = a'_{\ell-k+1}$ for $k\in [1,\ell]$, i.e., as a vector,
\begin{align} \label{eq: ep 1 = ep -1}
    (a_1, a_{2}, \ldots, a_\ell )	= (a'_{\ell}, a'_{\ell-1}, \ldots, a'_1).
\end{align}
\end{lemma}

Leclerc's examples in \cite{Lec03} use the parametrization $\PBW_{\bfi,1}(b)$. Lemma \ref{Lem: two ld} allows us to translate this parametrization into the convention $\PBW_{\bfi}=\PBW_{\bfi,-1}$ adopted in this paper. In what follows, we write
$$
\PBW_{\bfi}(b)=(a_1,a_2,\ldots,a_\ell)
$$
and denote by $\one$ the highest weight element of $B(\infty)$.

\subsection{Type $A_5$} 
Let $ \triA_{5} := (v, J, \La)$ be the triple given in Proposition \ref{Prop: triX in ft} \ref{P1} \ref{P1: b} for the case $n=5$, i.e., 
$$
v = s_1 s_3 s_5, \qquad J = \{2, 4\}, \qquad \La = \La_2 + \La_{4}.
$$
Proposition \ref{cor: condition for PI} \ref{cor: condition for PI 2} says 
$$
b_{\triA_{5}} = \tf_2 \tf_4 \tf_1 \tf_3^2 \tf_5 \tf_2 \tf_4 (\one) \in B(\infty).
$$
Let
$
\bfi = (5,4,5,3,4,5,2,3,4,5,1,2,3,4,5) \in \Seq(w_0).
$
Note that $\bfi^\vee$ coincides with the reduced expression of $w_0$ in \cite[Section 2.7]{Lec03} (see Lemma \ref{Lem: two ld}).
Then  the PBW datum of $b_{\triA_{5}} $ with respect to $\bfi$ is given by
$$
\PBW_\bfi (b_{\triA_{5}}) = (0, 1, 0, 0, 0, 1, 0, 1, 0, 0, 0, 0, 0, 1, 0),
$$
which coincides with that of $b$ in \cite[Section 2.7]{Lec03} (see \eqref{eq: ep 1 = ep -1}).

\subsection{Type $B_3$ and $C_3$} \label{sec: Lec B3C3}
Let $ \triX_{3} := (v, J, \La)$ ($\triX \in \{\triB, \triC \}$) be the triple given by 
$$
v = s_1 s_3, \qquad J = \{2\}, \qquad \La = \La_2
$$
(see Proposition \ref{Prop: triX in ft} \ref{P2} \ref{P2: b} and \ref{P3} \ref{P3: b}). 
Note that the index labeling of the Dynkin diagram in \cite[Sections 2.4 and 2.5]{Lec03} differs from that used in this paper (see Proposition \ref{Prop: triX in ft}).
Proposition \ref{cor: condition for PI} \ref{cor: condition for PI 2} says 
\begin{align*}
b_{\triB_{3}} = \tf_2 \tf_1 \tf_3^2 \tf_2 (\one) \qtq 
b_{\triC_{3}} = \tf_2 \tf_1 \tf_3 \tf_2 (\one).
\end{align*}
Let
$
\bfi = (1, 2, 1, 3, 2, 3, 1, 2, 3) \in \Seq(w_0).
$
Note that $\bfi^\vee$ coincides with the reduced expression of $w_0$ in \cite[Section 2.4 and 2.5]{Lec03} (see Lemma \ref{Lem: two ld}).
The PBW datum of $b_{\triX_3}$ is then given by
$$
\PBW_\bfi(b_{\triX_3})=(0,1,0,0,0,0,0,1,0),
$$
which coincides with that of $b$ in \cite[Sections 2.4 and 2.5]{Lec03} (see \eqref{eq: ep 1 = ep -1}).

\subsection{Type $D_4$} 
Let $ \triD_{4} := (v, J, \La)$ be the triple given by 
$$
v = s_1 s_3 s_4, \qquad J = \{2\}, \qquad \La = \La_2
$$
(see Proposition \ref{Prop: triX in ft} \ref{P4} \ref{P4: b}). 
Note that the index labeling of the Dynkin diagram in \cite[Sections 2.6]{Lec03} differs from that used in this paper (see Proposition \ref{Prop: triX in ft}).
Proposition \ref{cor: condition for PI} \ref{cor: condition for PI 2} says 
\begin{align*}
b_{\triD_{4}} = \tf_2 \tf_4 \tf_3 \tf_1 \tf_2 (\one) \in B(\infty).
\end{align*}
Let $\bfi = (4, 2, 4, 3, 2, 4, 1, 2, 4, 3, 2, 1) \in \Seq(w_0)$.
Note that $\bfi^\vee$ coincides with the reduced expression of $w_0$ in \cite[Section 2.4 and 2.5]{Lec03} (see Lemma \ref{Lem: two ld}).
Then  the PBW datum of $b_{\triD_{4}} $ is given by 
$$
\PBW_\bfi (b_{\triD_{4}}) = (0, 1, 0, 0, 0, 1, 0, 0, 0, 0, 1, 0),
$$
which coincides with that of $b$ in \cite[Section 2.6]{Lec03} (see \eqref{eq: ep 1 = ep -1}).

\subsection{Type $G_2$} 
Let  
$$
\triG_{ \clubsuit} = ( s_2, \{1\}, \La_1) \qtq 
 \triG_{\bigstar } = (s_1, \{2\}, \La_2 ) 
$$
 be the triple given in Proposition \ref{Prop: triX in ft} \ref{P9}. 
Proposition \ref{cor: condition for PI} \ref{cor: condition for PI 2} says 
\begin{align*}
b_{\triG_{ \clubsuit}} = \tf_1 \tf_2^3 \tf_1 (\one) \qtq 
b_{\triG_{ \bigstar}} =  \tf_2 \tf_1 \tf_2  (\one)
\end{align*}
Let
$\bfi = ( 1, 2, 1, 2, 1,2) \in \Seq(w_0)$.
Note that $\bfi^\vee$ coincides with the reduced expression of $w_0$ in \cite[Section 2.3]{Lec03} (see Lemma \ref{Lem: two ld}).
Then  the PBW datum of $b_{\triG_{ \clubsuit}}$ and  $b_{\triG_{ \bigstar}}$ are given by 
$$
\PBW_\bfi (b_{\triG_{ \clubsuit}}) = (1, 0, 0, 0, 1, 0) \qtq 
\PBW_\bfi (b_{\triG_{ \bigstar}}) =(0, 1, 0, 0, 0, 1).
$$
The  $\PBW_\bfi (b_{\triG_{ \bigstar}}) $ coincides with the PBW datum of $b$ in \cite[Section 2.3]{Lec03} and $\PBW_\bfi (b_{\triG_{ \clubsuit}})$ appears in \cite[Section 3.3]{Lec03} (see \eqref{eq: ep 1 = ep -1}).

\vskip 2em 

\section{Imaginary modules over quantum affine algebras} \label{Sec: im for qaa}

In this section, we discuss the imaginary modules obtained from the \stb imaginary elements in Proposition \ref{Prop: triX in ft} via the \emph{generalized quantum Schur--Weyl duality}, and present explicit examples in low rank cases.

Let $U_q'(\g)$ be a \emph{quantum affine algebra} over an algebraic closure of $\C(q)$. 
The \emph{Hernandez--Leclerc} category $\catCO$ is a distinguished monoidal subcategory of the category of finite-dimensional integrable $U_q'(\g)$-modules. We recall its definition later in this section; see also \cite{HL10, HL15} and \cite[Section 2]{KKOP24A}.

Let $\ddD=\{L_i\}_{i \in J}$ be a \emph{strong duality datum} of $\catCO$ of simply-laced finite type $X_n$ ($X \in \{ A,D,E \}$), which is the set of \emph{root modules} satisfying certain categorical conditions determined by the Cartan matrix of type $X_n$ (see \cite[Section 4.4]{KKOP23B}). 
The duality datum $\ddD$ gives rise to an exact monoidal functor, known as the generalized quantum Schur--Weyl duality,   
$$
\F_\ddD : R_{X_n}\gmod \rightarrow \catCO,
$$
where $R_{X_n}$ is the symmetric quiver Hecke algebra defined by the Cartan matrix of type $X_n$ (see \cite{KKK18} and see also \cite[Section 4]{KKOP23B}).
The following proposition follows directly from the fact that $\F_\ddD$ sends simple modules to simple modules (\cite[Theorem 4.10]{KKOP23B}).  

\begin{prop} \label{Prop: imag for Cg}
Let $\ddD$ be a strong duality datum of $\catCO$ of type $X_n$.
\bnum
\item For any  simple module $L$ in $R_{X_n}\gmod$,  $ L$ is imaginary if and only if  $\F_\ddD (L)$ is imaginary.
\item 
Let $\triX$ be the triple of type $X_n$ given in Proposition \ref{Prop: triX in ft}. Then $\F_\ddD ( \cL(b_\triX^{[t]}) )$ is imaginary   for any $t \in \Z_{>0}$.
\ee
\end{prop}

\subsection{$Q$-data for untwisted types}

Let $U_q'(\g)$ be of untwisted type $X_n^{(1)}$ 
and let $\catC_\g$ be the category of finite-dimensional integrable $U_q'(\g)$-modules. We adopt the notation for $U_q'(\g)$ from \cite[Section 2]{KKOP24A}.
We denote by $\triangle_0$ the subdiagram obtained from the affine Dynkin diagram of type $X_n^{(1)}$ by removing the index $0$. Let $I_0$ be the index set of $\triangle_0$.
For each $i \in I_0$, we denote by $\varpi_i$ the $i$th \emph{fundamental weight}. Then
there exists a unique simple module $V(\varpi_i)$,
called the \emph{fundamental module} of weight $\varpi_i$ satisfying  certain  conditions (see \cite[\S 5.2]{Kas02} and \cite[Section 2.2]{KKOP24A}). 
For a module $M \in \catC_\g$ and  $x \in \Q(q)^{\times}$, we set $M_x := M_z/ (z-x)M_z$, where $M_z$ is the \emph{affinization} of $M$ with indeterminate $z$.

We denote by $\gf$ the simply-laced finite simple Lie algebra determined by the category $\catCO$, which is listed in the below table
(see \cite{KKOP22} and see also \cite[Section 6]{KKOP24A}).  
\begin{align*} \label{Table: root system} \small
\begin{array}{|c||c|c|c|c|c|c|c|} 
\hline
\text{$U_q'(\g)$} & A_n^{(1)}  & B_n^{(1)} & C_n^{(1)} & D_n^{(1)} & E_{n}^{(1)} & F_{4}^{(1)} & G_{2}^{(1)}  \\
&(n\ge1)&(n\ge2)&(n\ge3)&(n\ge4)&(n=6,7,8) &&\\
\hline
\text{$\gf$} & A_n & A_{2n-1}    & D_{n+1}   &  D_n & E_{n} & E_{6} & D_{4}  \\
\hline
\end{array}
\end{align*} 
Let $\cmC_\fin$ and $I_\fin$ be the Cartan matrix of $\gf$ and its index set respectively. 
A $Q$-datum $\qQ = (\triangle_\fin, \sigma, \xi)$ is a triple of
\bna 
\item the Dynkin diagram $\triangle_\fin$ of $\gf$, 
\item the folding automorphism \(\sigma\) whose orbit set yields $\triangle_0$, i.e, $\sigma = \id$ if \(U_q'(\mathfrak g)\) is of simply laced affine type, and otherwise $\sigma$ is the one defined by \eqref{Eq: Folding Dynkin},
\item a \emph{height function} $\xi: I_\fin \rightarrow \Z$ defined in \cite[Definition 6.1]{KKOP24A}.
\ee
We write $d_\imath := | \{ \sigma^k(\imath) \mid k \in \Z_{\ge0} \}| $ for $\im \in I_\fin$.
For any sink $\imath \in I_\fin$ in $\qQ$, i.e., $\xi(\imath) \le \xi (\jmath)$ for any $\jmath$ with $(\al_\imath,  \al_\jmath)<0$, we define $s_\imath \qQ$ to be the $Q$-datum
$(\triangle_\fin, \sigma, s_\imath \xi)$ where $s_\imath \xi$ is the height function defined by $ s_\imath \xi (\jmath) := \xi (\jmath) + \delta_{\imath, \jmath} d_\imath $ for $\jmath \in I_\fin$. 
A reduced sequence $\bfi = (\imath_1, \imath_2, \ldots, \imath_\ell)$ of an element $w$ in the Weyl group $\weyl_\fin$ of $\gf$ is \emph{adapted to $\qQ$} if 
$$
\text{$\imath_k$ is a sink of $s_{\imath_{k-1}} s_{\imath_{k-2}} \cdots s_{\imath_{1}} \qQ $ for all $1 \le k \le \ell$.}  
$$

We denote by $\pi: I_\fin \twoheadrightarrow I_0$ the surjective map given by the folding automorphism $\sigma$.
We define $\hIQ  \seteq \{  (\imath, p) \in \If \times \Z \mid p-\xi(\imath) \in 2 d_\imath \Z  \}$ and 
\begin{align*}
\oIQ  \seteq \{  ( \pi(\imath), p) \in I_0 \times \Z  \mid  (\imath, p) \in \hIQ  \}. 
\end{align*}
Note that the map $ \hIQ  \buildrel \sim \over \longrightarrow \oIQ $ defined by $(\im, p) \mapsto (\pi(\im), p) $ is bijective.
We set  $\qm \seteq q^{1/\ord(\sigma)}$, where $\ord(\sigma)$ is the order of $\sigma$.  
For any $(i,p) \in \oIQ $, the associated fundamental module $V (Y_{i, p})\in \catC_\g$ is defined as
\begin{align*} 
V (Y_{i, p}) \seteq \begin{cases}
V(\varpi_{i})_{(-\qm)^p} & \text{ if } \g=A^{(1)}_n, \; C_n^{(1)}, \; D^{(1)}_n, \; E_{6,7,8}^{(1)}, \; G_2^{(1)}, \\
V(\varpi_{i})_{(-1)^{d(i, n)}(\qm)^p} & \text{ if } \g=B_n^{(1)}, \\
V(\varpi_{i})_{(-1)^{i}(\qm)^p} & \text{ if } \g=F_4^{(1)},
\end{cases}
\end{align*}
where $d(i, j)$ is the distance between $i$ and $j$ in $\triangle_0$.

We denote by $\catCO$ the \emph{Hernandez--Leclerc category}, that is, the smallest full subcategory of $\catC_\g$  that 
\bna 
\item  contains  $V (Y_{i, p})$ for all $ (i,p) \in \oIQ$, and
\item  is stable under taking subquotients, extensions and tensor products.
\ee
For any $ (i_1, p_1), (i_2, p_2), \ldots, (i_\ell, p_\ell) \in \oIQ$ with $p_1 \le p_2 \le \cdots \le p_{\ell}$, we define 
$$
V \left( \prod_{k} Y_{i_k, p_k } \right) := \hd \left( V (Y_{i_\ell, p_\ell}) \otimes  \cdots  V (Y_{i_2, p_2}) \otimes  V (Y_{i_1, p_1} ) \right),
$$
where $\hd (M)$ is the head of $M$.
We regard $\prod_{k} Y_{i_k, p_k } $ as a monomial in the polynomial ring in the variables $Y_{i,p}$.
Note that, for any monomial $m$ in $Y_{i,p}$, the module $V( m )$  is a simple module of $\catCO$ and every simple module of $\catCO$ has of this form (see \cite{CP91, CP95, Kas02} and see also \cite{HL10, HL16, KKOP23B}).

For each $Q$-datum $\qQ = (\triangle_\fin, \sigma, \xi)$, one can construct a \emph{complete duality datum} $\ddD_\qQ$ (see \cite[Section 6.2]{KKOP23B}), which gives a Schur--Weyl duality,   
$$
\F_\qQ : R_{\gf}\gmod \rightarrow \catCO.
$$

\begin{prop}[{\cite[Proposition 6.4]{KKOP23B}}] \label{Prop: cusp fund}
Let $\qQ = (\triangle_\fin, \sigma, \xi)$ be a $Q$-datum and $ \bfi \in \Seq(w_0) $ be a $\qQ$-adapted sequence of the longest element $w_0 \in \weyl_\fin$. Then the functor $\F_\qQ$ sends cuspidal modules $\cuspS_{\bfi, k}$ with respect to $\bfi$ to fundamental modules, where $\cuspS_{\bfi, k}$ is defined in \eqref{Eq: cuspS}. 
\end{prop}

\subsection{Examples arising from $Q$-data} \label{Sec: ex for qaa}

In this subsection, we provide examples of imaginary modules in $\catCO$ using 
Propositions \ref{Prop: imag for Cg} and \ref{Prop: cusp fund}.
In the following examples, we always denote by $\qQ = (\triangle, \sigma, \xi)$ a $Q$-datum and by $\bfi$ a $\qQ$-adapted reduced sequence of $w_0 \in \weyl_\fin$. Let $*$ be the involution of $\triangle_\fin$ defined by $\alpha_{i^*} = - w_0(\alpha_i)$ for $i\in \If$. For any sequence $\bfa = (a_1, \ldots, a_t)$ and $\bfb = (b_1, \ldots, b_k)$, we write $\bfa \cdot \bfb$ for the concatenation of $\bfa$ and $\bfb$. 
In the following subsections, we simply write $\bfa ^ k := \bfa \cdot \bfa^{k-1}$ for $k \in \Z_{>0}$, where we understand $\bfa^0 = \emptyset$.  Throughout this subsection, let $t \in \Z_{>0}$.

We present examples for untwisted affine types, omitting details for the exceptional types since they can be treated in the same way as the classical affine types.
Those examples are based on computations by SageMath (\cite{SageMath}).

\subsubsection{Type $A_6^{(1)}$}
In this case, $\gf$ is of type $A_6$. We choose
\begin{itemize}
\item $\qQ := (\triangle_{A_6}, \id, \xi)$, where  $\xi(i) = i-1 $ for $i\in I_0$, 
\item $\bfi := (1,2,1,3,2,4,1,3,5,2,4,6,1,3,5,2,4,1,3,2,1).$ 
\end{itemize}
Then $\F_\qQ$ sends the cuspidal module $\cuspS_{\bfi, k}$ to the fundamental modules $V(Y_{i,p})$, where $(i,p)$ is the position of the circled number $\cnum{k}$ in the following diagram:
\[ \small
\renewcommand{\arraystretch}{0.9} 
\setlength{\arraycolsep}{3pt}
\begin{array}{c|cccccccccccccccc}
 i\backslash p
   & &0&1&2&3&4&5&6&7&8&9&10\\
\hline 
1 && \cnum{1}& &\cnum{3}& & \cnum{7}& & \cnum{13}& & \cnum{18}& & \cnum{21}& &\\
2 &&       &\cnum{2}& &\cnum{5} & &\cnum{10} & &\cnum{16} & &\cnum{20} \\
3 &&        && \cnum{4}&& \cnum{8}&& \cnum{14}&& \cnum{19} \\
4 &&        &&&\cnum{6}&&\cnum{11}&&\cnum{17} \\
5 &&        &&&& \cnum{9} && \cnum{15} \\
6 &&        &&&&& \cnum{12}
\end{array}
\]

Let $ \triA_{6} $ be the triple given in Proposition \ref{Prop: triX in ft} \ref{P1} \ref{P1: b}.  Then we have 
$$
\PBW_{\bfi}(b_{\triA_{6}}^{[t]}) = (0, t, 0, 0, 0, 0, t, 0, 0, 0, t, 0, t, 0, 0, 0, 0, 0, 0, t, 0),
$$
which implies that $ \cL(b_{\triA_{6}}^{[t]}) \simeq \hd(\cuspS_{\bfi, 20}^{\circ t} \conv \cuspS_{\bfi, 13}^{\circ t} \conv \cuspS_{\bfi, 11}^{\circ t} \conv \cuspS_{\bfi, 7}^{\circ t} \conv \cuspS_{\bfi, 2}^{\circ t} ) $ is imaginary. Thus
\begin{align} \label{Ex: regular}
\F_\qQ(\cL(b_{\triA_{6}}^{[t]})) 
 \simeq V(Y_{2,9}^{t} Y_{1,6}^t Y_{4,5}^t Y_{1,4}^t Y_{2,1}^t )
\end{align}
is an imaginary module. 

\begin{remark} \label{rmk: LM}
When $t=1$, the module $V(Y_{2,9} Y_{1,6} Y_{4,5} Y_{1,4} Y_{2,1} )$ given in \eqref{Ex: regular} is regular, and hence its imaginarity follows from the result of \cite{LM18}. But, for $t > 1$,  the corresponding module is no longer regular.  
\end{remark}

\subsubsection{Type $B_3^{(1)}$}
In this case, $\g_\fin$ is of type $A_5$. We choose
\begin{itemize}
\item $\qQ := (\triangle_{A_5}, \sigma, \xi)$, where 
$\xi(1) = 5$,  $\xi(2) = 3$,   $\xi(3) = 0$, $\xi(4) = 1$, and $\xi(5) = 3$, 
\item $\bfi := (3,4,3,5,2,3,1,4,3,5,2,3,1,4,5)$. 
\end{itemize}
Then $\F_\qQ$ send $\cuspS_{\bfi, k}$  to $V(Y_{i,p})$, where $(i,p)$ is the position of the circled number $\cnum{k}$ in the following diagram:
\[ \small
\renewcommand{\arraystretch}{0.9} 
\setlength{\arraycolsep}{3pt}
\begin{array}{c|cccccccccccccccc}
 i\backslash p
   & &0&1&2&3&4&5&6&7&8&9&10 & 11\\
\hline 
1 && & &&  \cnum{4}& & \cnum{7}& & \cnum{10}& & \cnum{13}& & \cnum{15}\\
2 &&       &\cnum{2}& &\cnum{5} & &\cnum{8} & &\cnum{11} & &\cnum{14} \\
3 && \cnum{1}  && \cnum{3}&& \cnum{6}&& \cnum{9}&& \cnum{12} \\
\end{array}
\]
Note that this diagram coincides with that given in \cite[Section 6.2 Type $B_n^{(1)}$]{KKOP22} up to a spectral parameter shift.
\bna 
\item 
Let $ \triAs $ be the triple given in Proposition \ref{Prop: triX in ft} \ref{P1} \ref{P1: a}. Then we have
$$
\PBW_{\bfi}(b_{\triAs}^{[t]}) = (2t, 0, 2t, 0, 0, 0, 0, t, 0, t, 0, t, t, 0, t),
$$
which implies that $ \cL(b_{\triAs}^{[t]}) \simeq 
\hd(\cuspS_{\bfi, 15}^{\circ t} \conv \cuspS_{\bfi, 13}^{\circ t} \conv \cuspS_{\bfi, 12}^{\circ t} \conv \cuspS_{\bfi, 10}^{\circ t} \conv \cuspS_{\bfi, 8}^{\circ t} \conv \cuspS_{\bfi, 3}^{\circ 2t} \conv \cuspS_{\bfi, 1}^{\circ 2t} )  $ is  imaginary. Thus
\begin{align*}
\F_\qQ(\cL(b_{\triAs}^{[t]}))  \simeq V(Y_{1,11}^{t} Y_{1,9}^t Y_{3,8}^t  Y_{1,7}^t Y_{2,5}^t Y_{3,2}^{2t} Y_{3,0}^{2t} )
\end{align*}
is an imaginary module.

\item 
Let $ \triA_{5} $ be the triple given in Proposition \ref{Prop: triX in ft} \ref{P1} \ref{P1: b}. Then we have
$$
\PBW_{\bfi}(b_{\triA_{5}}^{[t]}) = (0, t, 0, 0, 0, t, 0, 0, t, 0, 0, 0, 0, t, 0),
$$
which implies that $ \cL(b_{\triA_{5}}^{[t]}) \simeq \hd(\cuspS_{\bfi, 14}^{\circ t} \conv \cuspS_{\bfi, 9}^{\circ t} \conv \cuspS_{\bfi, 6}^{\circ t} \conv \cuspS_{\bfi, 2}^{\circ t} ) $ is imaginary. Thus
\begin{align*}
\F_\qQ(\cL(b_{\triA_{5}}^{[t]}))\simeq V(Y_{2,9}^{t} Y_{3,6}^t Y_{3,4}^t  Y_{2,1}^t )
\end{align*}
is an imaginary module.

\ee

\subsubsection{Type $C_4^{(1)}$}
In this case, $\g_\fin$ is of type $D_5$. We choose
\begin{itemize}
\item $\qQ := (\triangle_{D_5}, \sigma, \xi)$, where 
$\xi(1) = 3$,  $\xi(2) = 2$,   $\xi(3) = 1$, $\xi(4) = 2$, and $\xi(5) = 0$, 
\item $\bfi := (5,3,2,4,1,3,2,5,1,3,2,4,1,3,2,5,1,3,2,1)$. 
\end{itemize}
Then $\F_\qQ$ send $\cuspS_{\bfi, k}$ to $V(Y_{i,p})$, where $(i,p)$ is the position of the circled number $\cnum{k}$ in the following diagram:
\[ \small
\renewcommand{\arraystretch}{0.9} 
\setlength{\arraycolsep}{3pt}
\begin{array}{c|cccccccccccccccc}
 i\backslash p
   & &0&1&2&3&4&5&6&7&8&9&10 & 11\\
\hline 
1 && & &&  \cnum{5}& & \cnum{9}& & \cnum{13}& & \cnum{17}& & \cnum{20}\\
2 &&      & &\cnum{3}& &\cnum{7} & &\cnum{11} & &\cnum{15} & &\cnum{19} \\
3 &&   & \cnum{2}&& \cnum{6}&& \cnum{10}&& \cnum{14} && \cnum{18} \\
4 && \cnum{1}  && \cnum{4}&& \cnum{8}&& \cnum{12}&& \cnum{16} \\
\end{array}
\]
Note that this diagram coincides with that given in \cite[Section 6.2 Type $C_n^{(1)}$]{KKOP22} up to a spectral parameter shift.

Let $ \triD_{5} $ be the triple given in Proposition \ref{Prop: triX in ft} \ref{P4} \ref{P4: b}. Then we have
$$
\PBW_{\bfi}(b_{\triD_{5}}^{[t]}) = (0, 0, t, 0, 0, 0, 0, 0, t, 0, 0, 0, t, 0, 0, 0, 0, 0, t, 0),
$$
which implies that $ \cL(b_{\triD_{5}}^{[t]}) \simeq \hd(\cuspS_{\bfi, 19}^{\circ t} \conv \cuspS_{\bfi, 13}^{\circ t} \conv \cuspS_{\bfi, 9}^{\circ t} \conv \cuspS_{\bfi, 3}^{\circ t} ) $ is imaginary. Thus
\begin{align*}
\F_\qQ(\cL(b_{\triD_{5}}^{[t]}))\simeq V(Y_{2,10}^{t} Y_{1,7}^t Y_{1,5}^t  Y_{2,2}^t )
\end{align*}
is an imaginary module.

\subsubsection{Type $D_5^{(1)}$}
In this case, $ \g_\fin$ is of type $D_5$. We choose
\begin{itemize}
\item $\qQ := (\triangle_{D_5}, \id, \xi)$, where 
$\xi(1) = 3$,  $\xi(2) = 2$,   $\xi(3) = 1$, $\xi(4) = 0$, and $\xi(5) = 0$, 
\item $\bfi := (4,5,3,2,4,5,1,3,2,4,5,1,3,2,4,5,1,3,2,1)$. 
\end{itemize}
Then $\F_\qQ$ send $\cuspS_{\bfi, k}$ to $V(Y_{i,p})$, where $(i,p)$ is the position of the circled number $\cnum{k}$ in the following diagram:
\[
\renewcommand{\arraystretch}{0.9} 
\setlength{\arraycolsep}{3pt}
\begin{array}{c|cccccccccccccccc}
 i\backslash p
   & &0&1&2&3&4&5&6&7&8&9\\
\hline 
1 && & &&  \cnum{7}& & \cnum{12}& & \cnum{17}& & \cnum{20}\\
2 &&      & &\cnum{4}& &\cnum{9} & &\cnum{14} & &\cnum{19}  \\
3 &&   & \cnum{3}&& \cnum{8}&& \cnum{13}&& \cnum{18}  \\
4 && \cnum{1}  && \cnum{5}&& \cnum{10}&& \cnum{15} \\
5 && \cnum{2}  && \cnum{6}&& \cnum{11}&& \cnum{16}\\
\end{array}
\]
Note that this diagram coincides with that given in \cite[Section 6.2 Type $D_n^{(1)}$]{KKOP22} up to a spectral parameter shift, except that the weights are twisted by \(*\).

Let $ \triD_{5} := (v, J, \La)$ be the triple given in Proposition \ref{Prop: triX in ft} \ref{P4} \ref{P4: b}. Then we have
$$
\PBW_{\bfi}(b_{\triD_{5}}^{[t]}) = (0, 0, 0, t, 0, 0, 0, 0, 0, 0, 0, t, 0, 0, 0, 0, 0, 0, t, 0),
$$
which implies that $ \cL(b_{\triD_{5}}^{[t]}) \simeq \hd(\cuspS_{\bfi, 19}^{\circ t} \conv \cuspS_{\bfi, 12}^{\circ t}  \conv \cuspS_{\bfi, 4}^{\circ t} ) $ is  imaginary. Thus 
\begin{align*}
\F_\qQ(\cL(b_{\triD_{5}}^{[t]}))\simeq V(Y_{2,8}^{t} Y_{1,5}^t Y_{2,2}^t )
\end{align*}
is an imaginary module.

\begin{remark} \label{rmk: Sakamoto}
Let us consider the case of type $D_4^{(1)}$.
Applying the same argument to the case of 
\begin{itemize}
\item $\qQ := (\triangle_{D_4}, \id, \xi)$, where 
$\xi(1) = 2$,  $\xi(2) = 1$,   $\xi(3) = 0$ and $\xi(4) = 0$,
\item $\bfi := (3,4,2,1,3,4,2,1,3,4,2,1)$, 
\end{itemize}
we obtain the following imaginary modules 
\begin{align*}
\F_\qQ(\cL(b_{\triDs}^{[t]}))\simeq V(Y_{1,6}^{t} Y_{1,4}^{t} Y_{2,3}^{t} Y_{3,0}^{t} Y_{4,0}^{t}      ), 
\qquad
\F_\qQ(\cL(b_{\triD_{4}}^{[t]}))\simeq V(Y_{2,5}^{t}  Y_{2,1}^t ).
\end{align*} 
For $t=1$, the above imaginary modules coincide with those given in \cite[Section 6.4]{Saka25}, and the second one also appears in \cite[Section 11]{EL24}.
\end{remark}

\subsubsection{Type $E_6^{(1)}$}
In this case, $ \g_\fin$ is of type $E_6$. We choose
\begin{itemize}
\item $\qQ := (\triangle_{E_6}, \id, \xi)$, where 
$\xi(1) = 0$,  $\xi(2) = 3$,   $\xi(3) = 1$, $\xi(4) = 2$, $\xi(5) = 3$, and $\xi(6) = 4$, 
\item $\bfi := (1,3,1,4) \cdot (3,2,5,1,4,6)^{4} \cdot (3,2,5,1,4,3,2,1) $. 
\end{itemize}

\bna 
\item 
Let $ \triE_{6, \bigstar} $ be the triple given in Proposition \ref{Prop: triX in ft} \ref{P5} \ref{P5: a}. 
Then  
\begin{align*}
\F_\qQ(\cL(b_{\triE_{6, \bigstar}}^{[t]}))  \simeq 
V(Y_{1,14}^{t}  Y_{1,12}^t  Y_{1,10}^t Y_{3,7}^t Y_{1,4}^{t} Y_{1,2}^{t} Y_{1,0}^{t} )
\end{align*}
is an imaginary module.

\item 
Let $ \triE_{6, \clubsuit} $ be the triple given in Proposition \ref{Prop: triX in ft} \ref{P5} \ref{P5: b}. 
Then 
\begin{align*}
\F_\qQ(\cL(b_{\triE_{6, \clubsuit}}^{[t]}))  \simeq 
V(Y_{2,13}^{t}  Y_{2,11}^t  Y_{2,5}^t Y_{2,3}^t )
\end{align*}
is an imaginary module.

\ee

\subsubsection{Type $E_7^{(1)}$}
In this case, $ \g_\fin$ is of type $E_7$. We choose
\begin{itemize}
\item $\qQ := (\triangle_{E_7}, \id, \xi)$, where 
$\xi(1) = 5$,  $\xi(2) = 4$,   $\xi(3) = 4$, $\xi(4) = 3$, $\xi(5) = 2$, $\xi(6) = 1$ and $\xi(7) = 0$, 
\item $\bfi := (7,6,5,7,4,6) \cdot (3,2,5,7,1,4,6)^{7} \cdot (3,2,5,1,4,3,2,1) $. 
\end{itemize}
Let $ \triE_{7, \bigstar} $ be the triple given in Proposition \ref{Prop: triX in ft} \ref{P6}. Then 
\begin{align*}
\F_\qQ(\cL(b_{\triE_{7, \bigstar}}^{[t]}))  \simeq 
V(Y_{1,21}^{t}  Y_{1,19}^t  Y_{1,17}^t Y_{1,9}^t Y_{1,7}^t Y_{1,5}^t )
\end{align*}
is an imaginary module.

\subsubsection{Type $E_8^{(1)}$}
In this case, $ \g_\fin$ is of type $E_8$. We choose
\begin{itemize}
\item $\qQ := (\triangle_{E_8}, \id, \xi)$, where 
$\xi(1) = 6$,  $\xi(2) = 5$,   $\xi(3) = 5$, $\xi(4) = 4$, $\xi(5) = 3$, $\xi(6) = 2$, $\xi(7) = 1$ and $\xi(8) = 0$, 
\item $\bfi := (8,7,6,8,5,7,4,6,8) \cdot (3,2,5,7,1,4,6,8)^{12} \cdot (3,2,5,7,1,4,6,3,2,5,1,4,3,2,1) $. 
\end{itemize}
Let $ \triE_{8, \bigstar} $ be the triple given in Proposition \ref{Prop: triX in ft} \ref{P7}. Then 
\begin{align*}
\F_\qQ(\cL(b_{\triE_{8, \bigstar}}^{[t]}))  \simeq 
V(Y_{8,28}^{t}  Y_{8,26}^t  Y_{8,24}^t Y_{8,22}^t Y_{8,20}^t 
Y_{8,8}^{t}  Y_{8,6}^t  Y_{8,4}^t Y_{8,2}^t Y_{8,0}^t
)
\end{align*}
is an imaginary module.

\subsubsection{Type $F_4^{(1)}$}
In this case, $ \g_\fin$ is of type $E_6$. We choose
\begin{itemize}
\item $\qQ := (\triangle_{E_6}, \sigma, \xi)$, where 
$\xi(1) = 5$,  $\xi(2) = 1$,   $\xi(3) = 3$, $\xi(4) = 0$, $\xi(5) = 1$, and $\xi(6) = 3$, 
\item $\bfi := (4,5,2) \cdot (4,6,3,2,4,1,5,2)^{4} \cdot (6) $. 
\end{itemize}

\bna 
\item 
Let $ \triE_{6, \bigstar} $ be the triple given in Proposition \ref{Prop: triX in ft} \ref{P5} \ref{P5: a}. 
Then
\begin{align*}
\F_\qQ(\cL(b_{\triE_{6, \bigstar}}^{[t]}))  \simeq 
V(Y_{1,19}^{t}  Y_{1,17}^t  Y_{1,15}^t Y_{1,7}^t Y_{1,5}^{t} Y_{1,3}^{t} )
\end{align*}
is an imaginary module.

\item 
Let $ \triE_{6, \clubsuit} $ be the triple given in Proposition \ref{Prop: triX in ft} \ref{P5} \ref{P5: b}. 
Then  
\begin{align*}
\F_\qQ(\cL(b_{\triE_{6, \clubsuit}}^{[t]}))  \simeq 
V(Y_{4,17}^{t}  Y_{4,15}^t  Y_{4,13}^t Y_{4,5}^t Y_{4,3}^t Y_{4,1}^t )
\end{align*}
is an imaginary module.

\ee

\subsubsection{Type $G_2^{(1)}$}
In this case, $ \g_\fin$ is of type $D_4$. We choose
\begin{itemize}
\item $\qQ := (\triangle_{D_4}, \sigma, \xi)$, where 
$\xi(1) = 4$,  $\xi(2) = 1$,   $\xi(3) = 2$, and $\xi(4) = 0$, 
\item $\bfi := (4,2,3,2,1,2,4,2,3,2,1,2)$. 
\end{itemize}

\bna 
\item 
Let $ \triDs$ be the triple given in Proposition \ref{Prop: triX in ft} \ref{P4} \ref{P4: a}. 
Then 
\begin{align*}
\F_\qQ(\cL(b_{\triDs}^{[t]}))  \simeq V(Y_{1,10}^{t}  Y_{2,7}^t  Y_{2,5}^t Y_{2,3}^t Y_{1,0}^{t} )
\end{align*}
is an imaginary module.

\item 
Let $ \triD_{4}$ be the triple given in Proposition \ref{Prop: triX in ft} \ref{P4} \ref{P4: b}. 
Then 
\begin{align*}
\F_\qQ(\cL(b_{\triD_{4}}^{[t]}))\simeq V(Y_{2,11}^{t} Y_{2,9}^t Y_{2,3}^t  Y_{2,1}^t )
\end{align*}
is an imaginary module.
\ee

\vskip 2em

\appendix
\section{Bosonic extensions and quantum twist automorphisms} \label{Sec: appA}

In this appendix, we shall show that the expansion of $\eta_w(x) x$ with respect to the basis $\tGup(w)$ has a non-zero constant term for any real element $x \in \tGup(w)$ (Proposition \ref{prop: eta x x}). This is used in the proof of Lemma \ref{Lem: fixed by D} for the case where $\Gup(w)$ is of non-symmetric type.
Since the proof requires additional background including \emph{bosonic extensions}, which is not used elsewhere in the paper, we present the proof here separately. We keep all the notation from Sections \ref{Sec: prelim} and \ref{Sec: QTA} and follow the notation used in \cite{Kimura12, KiOy21, KKOP25A} for quantum groups, quantum twist automorphisms and bosonic extensions.

\subsection{Unipotent quantum matrix coefficients}

Let $U_q(\g)$ be the quantum group associated with a generalized symmetrizable Cartan matrix $\cmC = (c_{i,j})_{i,j\in I}$ and let $A_q(\n)$ be the corresponding quantum unipotent coordinate ring. 
Throughout this appendix, we identify $A_q(\n)$ with the negative half $U_q^-(\g)$ via the Lusztig's symmetric bilinear form $(\ , \ )_L$ on $U_q^-(\g)$ (see \cite{LusztigBook} and see also \cite[Definition 3.7]{KiOy21}).

Let $ V(\la)$ be the irreducible highest wight $U_q(\g)$-module of highest weight $\la \in \wlP_+$ and $\Gup(\la)$ be the dual canonical basis (or upper global basis) of $V(\la)$. 
We write $\wt_\la (x)$ for the weight of $x$ in $V(\la)$, and recall  extremal vectors $u_{w\la}$ of $V(\la)$ given in Section \ref{Sec: HWC}. 
Via the dual of the projection $\pi_\la: U_q^-(\g) \rightarrow V(\la)$ mapping $x \mapsto xu_\la$, we regard $\Gup(\la)$ as a subset of $\Gup(\infty)$.  
For $w \in \weyl$, we write $B_w(\la)$ for the crystal of the \emph{Demazure module} $V_w(\la) := U_q^+ ( u_{w \la} )$, and write $\Gup_w(\la) := \{ \Gup(b) \mid b\in B_w(\la) \}$ (see \cite{K93A} and see also \cite[Proposition 3.34]{KiOy21}). We denote by $(\ , \ )_\la$  the non-degenerate bilinear form on $V(\la)$ such that $(u_\la, u_\la)_\la=1$ and $(f_ix,y)_\la = (x,e_iy)_\la$ for $i\in I$ and $x,y\in V(\la)$ (see \cite[Section 3.3]{KiOy21} for example).

For any $u,v \in V(\la)$, define the \emph{unipotent quantum matrix coefficient} $D_{u,v} \in U_q^-(\g)$ by 
$$
(D_{u,v}, f)_L = ( u, f v )_\la \qquad \text{ for any $f \in U_q^-(\g)$}
$$
(see \cite[Definition 3.40]{KiOy21}).

\begin{lemma} \label{Lem: DD=D}
Let $w\in \weyl$, $ \la, \mu \in \wlP_+$, $x \in \Gup_w(\la)$ and  $y \in \Gup_w(\mu)$. 
Then we have 
$$
D_{u_{w(\la+\mu)}, p(x \otimes y)} = q^{(w\la - \wt_\la(x), w \mu)} D_{u_{w\la}, x} D_{u_{w\mu}, y}, 
$$
where $p: V(\la) \otimes V(\mu) \rightarrow V(\la+\mu)$ is the dual map of the embedding $\iota: V(\la+\mu) \rightarrow  V(\la) \otimes V(\mu)$ mapping $ u_{\la+\mu}  \mapsto u_\la \otimes u_\mu$.
\end{lemma}
\begin{proof}
We set 
$(a\otimes b, c\otimes d )_{\la,\mu} := (a, b)_\la (c,d)_\mu  $ 
for $a\otimes b, c\otimes d \in V(\la) \otimes V(\mu) $
(see \cite[(2.5.2)]{K91}).
Let $\vee: U_q(\g)\buildrel \sim \over \rightarrow U_q(\g)$ be the $\Q(q)$-algebra involution defined by $ e_i^\vee = f_i$, $ f_i^\vee = e_i$, and $(q^h)^\vee = q^{-h}$. 
Let $\pi_{w,\la}^\vee : U_q^-(\g) \twoheadrightarrow V_w(\la)$ be the surjective $\Q(q)$-linear map given by $ \pi_{w,\la}^\vee (y) = y^\vee u_{w\la}$ for $y\in U_q^-(\g)$, and set $j^\vee_{w,\la}$ to be the dual of $\pi_{w,\la}^\vee$ (\cite[Definition 3.43]{KiOy21}). 

Let $\Delta : U_q(\g) \rightarrow U_q(\g) \otimes U_q(\g)$ be the coproduct of $U_q(\g)$ defined in \cite[Section 2.1.9]{Kimura12}, and $r: U_q^-(\g) \rightarrow U_q^-(\g) \otimes U_q^-(\g)$ be the twist coproduct of $U_q^-(\g)$ defined in \cite[Section 2.1.10]{Kimura12}.
Note that, for a homogeneous element $z \in U_q(\g)$ with $r(z) = \sum z_{(1)} \otimes z_{(2)}$, we have 
\begin{align*}
\Delta(z) = \sum z_{(1)} q^{-\xi (\wt (z_{(2)}))} \otimes z_{(2)},
\end{align*}
where $\xi : \wlP \buildrel \sim \over\rightarrow \wlP^\vee$ is the $\Q$-linear isomorphism determined by $\langle \xi(\nu), \nu' \rangle = (\nu, \nu')$ for $\nu, \nu' \in \wlP$ (see \cite[Lemma 2.5]{Kimura12} and see also \cite[Lemma 4.1]{KP20}). 
Since $\Delta \circ \vee = \text{flip} \circ (\vee \otimes \vee) \conv \Delta  $ where $\text{flip}(a\otimes b) = b \otimes a$, we have 
\begin{align} \label{Eq: Delta r}
\Delta(z^\vee) = \sum z_{(2)}^\vee  \otimes z_{(1)}^\vee q^{\xi (\wt (z_{(2)}))},
\end{align}

We set $\gamma := w\la - \wt_\la (x)$. For any homogeneous $z \in U_q^-(\g)$, we have 
\begin{align*}
( j_{w,(\la+\mu)}^\vee (p(x \otimes y)), z )_L 
&= ( p(x \otimes y), z^\vee u_{w(\la+\mu)} )_{\la+\mu} \\
&= (  x \otimes y , \iota ( z^\vee u_{w(\la+\mu)} ) )_{\la, \mu} \\
&= (  x \otimes y , \Delta( z^\vee) (u_{w\la} \otimes u_{w\mu}) )_{\la, \mu} \\
&=  \sum  q^{(\wt (z_{(2)}), w\mu)}  ( x, z_{(2)}^\vee u_{w\la} )_\la ( y, z_{(1)}^\vee u_{w\mu} )_\mu \\ 
&=  q^{(\gamma, w\mu)} \sum     ( y, z_{(1)}^\vee u_{w\mu} )_\mu  ( x, z_{(2)}^\vee u_{w\la} )_\la \\
&=  q^{(\gamma, w\mu)} \sum     ( j_{w,\mu }^\vee (y), z_{(1)})_L  ( j_{w, \la }^\vee(x), z_{(2)} )_L \\
&=  q^{(\gamma, w\mu)}     ( j_{w, \mu }^\vee (y) \otimes j_{w, \la }^\vee(x), \sum  z_{(1)} \otimes z_{(2)} )_L \\
&= ( q^{(\gamma, w\mu)} j_{w, \mu }^\vee (y)  j_{w, \la }^\vee(x), z )_L,
\end{align*}
where 
the first and sixth equalities follow from \cite[Definition 3.43]{KiOy21}, 
the second follows from \cite[(2.5.6)]{K91}, 
the third follows from $\iota(u_{w (\la +\mu) }) = u_{w \la} \otimes u_{w\mu} $, 
the forth follows from \eqref{Eq: Delta r},
the fifth follows from the fact that $ ( x, z_{(2)}^\vee u_{w\la} )_\la = 0 $ unless $ \wt(z_{(2)}) = \gamma $, 
the eighth follows from the compatibility of $r$ with $(\ ,\ )_L$ (\cite{LusztigBook}, \cite[Section 2.1.10]{Kimura12}). This implies that 
$$
j_{w,(\la+\mu)}^\vee (p(x \otimes y)) = q^{(\gamma, w\mu)} j_{w, \mu }^\vee (y)  j_{w, \la}^\vee(x).
$$
Applying the anti-involution $*$ to the above, we have 
\begin{align*}
D_{u_{w(\la+\mu)}, p(x \otimes y)} &= *j_{w,(\la+\mu)}^\vee (p(x \otimes y))
= q^{(\gamma, w\mu)} (* j_{w, \la }^\vee(x)) (*j_{w, \mu }^\vee (y) )\\
&= q^{(w\la - \wt_\la(x), w \mu)} D_{u_{w\la}, x} D_{u_{w\mu}, y}
\end{align*}
by \cite[Proposition 3.46 (ii)]{KiOy21}.
\end{proof}

\begin{lemma} \label{Lem: positive Gup}
Let $w\in \weyl$, $ \la, \mu \in \wlP_+$, $x \in \Gup_w(\la)$ and  $y \in \Gup_w(\mu)$. 
Then we have 
$$
q^{(w\la - \wt_\la(x), w \mu)} D_{u_{w\la}, x} D_{u_{w\mu}, y} 
\in \sum_{z \in \Gup_w(\la+\mu)} \Z[q]z.
$$
\end{lemma}
\begin{proof}
Recall the surjective homomorphism $p: V(\la) \otimes V(\mu) \rightarrow V(\la+\mu)$ given in Lemma \ref{Lem: DD=D}.
Let $\A = \Z[q^{\pm 1}]$ and write $V_\A (\lambda)$ for the $\A$-lattice of $V(\la)$ generated by $\Gup(\la)$. Since   $p( V_\A (\lambda) \otimes_\A V_\A (\mu) ) \subset V_\A (\lambda+\mu)$ (see \cite[Section 4.1]{KP20} for example), we have $p(x\otimes y) \in \sum_{z \in \Gup(\la+\mu) } \A z $.  
Since all coefficients of $p(x\otimes y)$ are regular at $q=0$ (\cite[Section 4.2]{K91}, \cite[Section 3.3.3]{Kimura12}), we conclude that $p(x\otimes y) \in \sum_{z \in \Gup(\la+\mu) } \Z[q] z $. Hence the assertion follows from Lemma \ref{Lem: DD=D} and \cite[Proposition 3.46]{KiOy21}.  
\end{proof}

\subsection{Bosonic extensions}

\begin{definition}
The \emph{bosonic extension} $\hA$ associated with $\cmC$ is the $\Q(q^{\pm1/2})$-algebra generated by 
$\{f_{i,p}\mid i\in I,p\in \Z\}$ with the following relations:
\bna
\item 
$\sum\limits_{k=0}^{1- \langle h_i,\al_j \rangle}(-1)^k\bnom{1-\ang{h_i,\al_j}\\k}_i
f_{i,p}^kf_{j,p}f_{i,p}^{1-\langle h_i,\al_j \rangle-k}=0$\quad for any $i,j\in I$ such that $i\not=j$ and $p\in \Z$,

\item 
$f_{i,m}f_{j,p}=q^{(-1)^{p-m+1}(\al_i,\al_j)}f_{j,p}f_{i,m}
+  \delta_{j,i} \delta_{p,m+1} (1-q_i^2)$\quad if $m<p$,
\ee
where $ \delta_{a,b} = 1$ if $a=b$ and $ \delta_{a,b} = 0$ otherwise.
\end{definition}
We set $\wt(f_{i,p}) := (-1)^{p+1} \al_i$ for $i\in I$, which give the $\rlQ$-graded algebra structure on $\hA$. 
We have the following (anti-)automorphisms of $\hA$:
\bna
\item the $\Qh$-algebra automorphism $ \bD $ given by
\begin{align*} 
  \bD(f_{i,p}) = f_{i, p+1} \qquad \text{ for any $i\in I$ and $p\in \Z$.}
\end{align*}
\item  the $\Q$-algebra anti-automorphism $\overline{\ \cdot \ }$ defined by
$$
\overline{q^{\pm1/2}}=q^{\mp1/2}\qtq \overline{f_{i,p}}=f_{i,p}.
$$ 
\ee
 We set $c(x) := q^{N(\wt(x))} \overline{x}$ for a homogeneous $x \in \hA$, where $N(\al) := (\al, \al)/2$ for $\al \in \rlQ$ (\cite[Section 5]{KKOP25A}).  

We now extend the scalar of $ A_q(\n)$ (resp.\ $ A_q(\n(w))$) to $\Q(q^{\pm 1/2})$, which is also denoted by $A_q(\n)$ (resp.\ $ A_q(\n(w))$). We denote by  
$$
Q_w:  A_q(\n) \rightarrow  \tA_q(\n(w)) 
$$
the quotient map given in \cite[Section 1.6]{KKOP23A} (see also \cite[Section 3.5, Theorem 4.13]{KiOy21}).
For an interval $[a,b]$, we define $ \hA[a,b]$ to be the subalgebra of $\hA$ generated by $f_{i,p}$ for $i\in I $ and $p\in [a,b]$. We simply write $\hA[a]:= \hA[a, a]$. For each $m\in \Z$, we define a $\Qh$-algebra isomorphism 
$$
\varphi_m : A_q(\n) \buildrel \sim \over \longrightarrow \hA[m] \qquad \text{ for $\varphi_m(\ang{i}) = q_i^{1/2}f_{i,m}$ ($i\in I$),}  
$$
where $\ang{i} \in \Gup(\infty)$ is the dual of $f_i$ (\cite[Section 7]{KKOP25A}).

Recall the quantum twist automorphism $\eta_w$ given in Section \ref{Sec: qta}.
\begin{lemma} \label{Lem: Phi_w Hw}
We define 
$$
H_w := \{ \varphi_0(a) \varphi_{-1}(b) \mid a \in A_q(\n(w)),  b \in A_q(\n)\} \subset \hA[-1,0].
$$
Then we have the following.
\bni
\item $H_w$ is a subalgebra of $\hA[-1,0]$.   
\item Define 
$$
\Psi_w: H_w \rightarrow \tA(\n(w)), \qquad \varphi_0(a) \varphi_{-1}(b) \mapsto a \cdot \eta_w(Q_w(b)).
$$  
Then $\Psi_w$ is an algebra homomorphism.
\ee
\end{lemma}
\begin{proof}
(i) It suffices to show that $ f_{i,-1} u $ is contained in $ H_w $ for any $i\in I$ and $u \in \varphi_0(A_q(\n(w)))$. We set $ E_{i,0}^\star (u) := f_{i,-1} u - q^{- (\al_i, \wt(u))} u f_{i, -1} $ (\cite[(6.2)]{KKOP25A}), that is, 
\begin{align} \label{Eq: fi-1 u}
f_{i,-1} u = E_{i,0}^\star (u) + q^{- (\al_i, \wt(u))} u f_{i, -1}.
\end{align}
Since $e_i^*(x) \in A_q(\n(w))$ for any $x\in A_q(\n(w))$, it follows from \cite[Lemma 6.1 (iii)]{KKOP25A} that $ f_{i,-1} u $ is contained in $H_w$.

(ii) It is enough to show that, for $x \in A_q(\n(w))$ and $y\in A_q(\n)$, 
\begin{align} \label{Eq: PsiPsi}
\Psi_w( \varphi_{-1}(y) \varphi_{0}(x)) = \Psi_w( \varphi_{-1}(y))  \Psi_w(\varphi_{0}(x)) = \eta_w (Q_w (y)) x.
\end{align}
Since $A_q(\n)$ is generated by $\ang{i}$ for $i\in I$, 
we may also assume that $y = q_i^{-1/2} \ang{i} $ so that $\varphi_{-1}(y) = f_{i, -1}$. We write $\beta := \wt(x)$ and $X := \varphi_0(x)$, and identify $\hA[0]$ with $U_q^-(\g)$ via the mapping $f_{i,0} \mapsto f_i$ for $i\in I$.

For $i\in I$, we set $ \mathsf{F}_i := D_{u_{w\La_i}, u_{\La_i}}$ and $\mathsf{E}_i := e_i^* D_{u_{w\La_i}, u_{\La_i}} = D_{u_{w\La_i}, u_{s_i\La_i}} $. By \cite[Theorem 6.1]{KiOy21}, we have 
\begin{align} \label{Eq: eta i}
\eta_w(\ang{i}) = \eta_w(D_{u_{s_i\La_i}, u_{\La_i}})
 = q_i \mathsf{F}_i^{-1} \mathsf{E}_i.
\end{align}
\bna 
\item Let $ \zeta_i = 1-q_i^2$ for $i\in I$. By \cite[Lemma 6.1 (iii)]{KKOP25A}, we have 
\begin{align*}
\varphi_{-1}(y) \varphi_{0}(x) &= f_{i,-1} X
= E_{i,0}^\star(X) + q^{-(\al_i, \beta)} X f_{i,-1} \\ 
&= q^{-(\al_i, \beta+\al_i)} \zeta_i e_i^*(X) + q^{-(\al_i, \beta)} X f_{i,-1},
\end{align*}
which implies 
\begin{align*} 
\Psi_w( \varphi_{-1}(y) \varphi_{0}(x)) &= 
q^{-(\al_i, \beta+\al_i)} \zeta_i  (q_i^{1/2} e_i^*(x)) + q^{-(\al_i, \beta)} x ( q_i^{-1/2} q_i \mathsf{F}_i^{-1} \mathsf{E}_i) \\
&= q_i^{-3/2}  q^{-(\al_i, \beta)}  \zeta_i  e_i^*(x) + q_i^{1/2}q^{-(\al_i, \beta)}  (x  \mathsf{F}_i^{-1} \mathsf{E}_i) , 
\end{align*}
where $q_i^{\pm 1/2}$ on the right hand side follows from the definition of $\varphi_m$ and the second term on the right hand side follows from \eqref{Eq: eta i}.

\item Note that $ \mathsf{F}_i x = q^{(w\La_i + \La_i, \beta)} x \mathsf{F}_i $ by \cite[Proposition 4.7 (3)]{KiOy21}. It follows from 
\begin{align*}
e_i^* (\mathsf{F}_i x) &= \mathsf{F}_i e_i^*(x) + q^{(\al_i, \beta)} e_i^*(\mathsf{F}_i) x = \mathsf{F}_i e_i^*(x) + q^{(\al_i, \beta)} \mathsf{E}_i x,  \\ 
e_i^* ( x \mathsf{F}_i ) &= 
x e_i^*(\mathsf{F}_i) + q^{(\al_i, w\La_i - \La_i)} e_i^*(x) \mathsf{F}_i = x  \mathsf{E}_i + q^{(\al_i, w\La_i - \La_i)} e_i^*(x) \mathsf{F}_i 
\end{align*}
that 
\begin{align*}
q^{(\al_i, \beta)} \mathsf{E_i} x &= -\mathsf{F_i} e_i^*(x) 
+ q^{(w\La_i+\La_i, \beta)} x \mathsf{E}_i 
+ q^{(w\La_i+\La_i, \beta)} q^{(\al_i, w\La_i-\La_i)} e_i^*(x) \mathsf{F}_i \\ 
&= -\mathsf{F_i} e_i^*(x) 
+ q^{(w\La_i+\La_i, \beta)} x \mathsf{E}_i 
+ q_i^{-2} \mathsf{F}_i e_i^*(x)  \\
&=  q_i^{-2}\zeta_i \mathsf{F_i} e_i^*(x) 
+ q^{(w\La_i+\La_i, \beta)} x \mathsf{E}_i.
 \end{align*}
Thus we have 
\begin{align*}
\eta_w(Q_w(y))x &= q_i^{-1/2}\eta_w(\ang{i})x = q_i^{1/2} \mathsf{F}_i^{-1} \mathsf{E}_i x \\ 
&= q_i^{-3/2} q^{-(\al_i, \beta)}  \zeta_i  e_i^*(x) 
+ q_i^{1/2} q^{-(\al_i, \beta)}  q^{(w\La_i+\La_i, \beta)} \mathsf{F}_i^{-1} x \mathsf{E}_i \\ 
&= q_i^{-3/2} q^{-(\al_i, \beta)}  \zeta_i  e_i^*(x) 
+ q_i^{1/2} q^{-(\al_i, \beta)}  (x \mathsf{F}_i^{-1} \mathsf{E}_i).
\end{align*}
\ee
By (a) and (b), we have \eqref{Eq: PsiPsi}.
\end{proof}

Let $\trmG$ be the \emph{normalized global basis} of $\hA$ given in \cite[Definition 7.7]{KKOP25A}. Combining \eqref{Eq: fi-1 u} with \cite[Lemma 7.4 and 7.5]{KKOP25A}, one can show that $H_w$ is compatible with the basis $\trmG$, i.e.,  $\trmG \cap H_w$ is a $\Qh$-basis of $H_w$. 

\begin{lemma} \label{Lem: L0 L0}
Let $\mathfrak{L}$ be the $\Z[q^{1/2}]$-lattice of $H_w$ generated by $\trmG \cap H_w$, and write $\mathfrak{L}_0$ for the subset of elements of weight $0$ in $\mathfrak{L}$. Let $\widetilde{L}$ be the $\Z[q^{1/2}]$-lattice of $\tA_q(\n(w))$ generated by $\tGup(\infty)$. Then we have 
$$
\Psi_w( \mathfrak{L}_0) \subset \widetilde{L}.
$$
\end{lemma}
\begin{proof}
Note that $\mathfrak{L}_0$ is generated by $\varphi_0(a) \varphi_{-1}(b)$ for $a \in \Gup(w)$ and $ b\in \Gup(\infty)$ with $\wt(a) = \wt(b)$ (\cite[(7.3)]{KKOP25A}). Let $a \in \Gup(w)$ and $ b\in \Gup(\infty)$ such that $\wt(a) = \wt(b)$ and $Q_{w}(b) \ne 0$. Then it suffices to show that $ \Psi_w(\varphi_0(a) \varphi_{-1}(b) ) \in \widetilde{L}$.

Note that the braid symmetry $T_i$ used in \cite[Definition 3.23]{KiOy21} coincides with $\lTs_i$ (see \eqref{Eq: braid}) in this paper.   
Taking this difference in notation into account, \cite[Proposition 3.46 and Theorem 3.48]{KiOy21} says that there exists $x \in \Gup(\la)$ and $y\in  \Gup(\mu)$ for some $\la, \mu \in \wlP_+$ such that
$ a = D_{u_{w\la}, x} $ and $ b = D_{y, u_\mu} $.

We set $\gamma:= \wt(a) = \wt(b) = \wt_{\mu}(y) -\mu = w\la - \wt_\la(x) $. By \cite[Proposition 4.7 (3) and Theorem 6.1]{KiOy21},  we have 
\begin{align*}
\Psi_w(\varphi_0(a) \varphi_{-1}(b)) &= a \eta_w(Q_w(b))
= q^{-(\mu, \gamma)} D_{u_{w\la}, x}  D_{u_{w\mu}, u_{\mu}}^{-1} D_{u_{w\mu}, y} \\
&=  D_{u_{w\mu}, u_{\mu}}^{-1} ( q^{(w\mu, \gamma)}  D_{u_{w\la}, x} D_{u_{w\mu}, y}).
\end{align*}
Since $ (\mu + w\mu, \wt(D_{u_{w\la}, x}) + \wt(D_{u_{w\mu},y}) ) = (\mu + w\mu,  w\mu - \mu) =0$, we have 
$$
\Psi_w(\varphi_0(a) \varphi_{-1}(b)) \in \widetilde{L}
$$
by \cite[Definition 4.6 and Proposition 4.7 (3)]{KiOy21} and Lemma \ref{Lem: positive Gup}.
\end{proof}

For any interval $[a,b]$, we set $\trmG[a,b] := \trmG \cap \hA[a,b]$. Note that $\trmG[a,b]$ is a $\Qh$-linear basis of $\hA[a,b]$. When $a=b$, we simply write $\trmG[a] = \trmG[a,a] $.

\begin{prop} \label{prop: eta x x}
Let $w\in \weyl$. 
\bni
\item If $a \in \trmG[0]$ is real, then we have 
$$ 
\bD^{-1}(a) a \in 1+ q^{1/2}  \mathfrak{L}_0 
$$

\item If  $x \in \Gup(w)$ is real, then 
$$
q^{-N(\wt(x))} \eta_w(x) x \in 1+ q^{1/2}  \widetilde{L}. 
$$
\ee
Here $ \mathfrak{L}_0 $ and $\widetilde{L}$ are the lattices given in Lemma \ref{Lem: L0 L0}.
\end{prop}
\begin{proof}
(i) Recall the symmetric bilinear form $(\ , \ )_{\hA}$ on $\hA$ introduced in \cite[(6.10)]{KKOP25A}, and the map $\mathbf{M}: \hA \rightarrow \Qh$ defined in \cite[(6.8)]{KKOP25A}. By the definition of $(\ , \ )_{\hA}$, we have $(x,y)_{\hA} = \mathbf{M} (x \bD(y))$ for $x,y\in \hA$.   

Let $z := \bD^{-1}(a) a$. Since $a$ is real and  its squares $a^2$ is bar-invariant, the element $a^2$ is contained in $\trmG$. Thus, by \cite[Lemma 6.2, Proposition 7.9]{KKOP25A},  we have 
\begin{equation}\label{Eq: zaa hA}
\begin{aligned}
(z,1)_{\hA} &= \mathbf{M} ( \bD^{-1}(a) a ) = \mathbf{M} ( a \bD(a)) = (a,a)_{\hA} \in 1 + q(\Z[[q]]), \\
(z,z)_{\hA} &= \mathbf{M} ( (a \bD(a^2) ) \bD^2(a) )
= \mathbf{M} ( a^2 \bD(a^2)   ) \\
&= (a^2,   a^2)_{\hA} \in 1 + q(\Z[[q]]).
\end{aligned}
\end{equation}
We now write $z = \sum_{k=1}^N c_k q^{t_k} y_k $ for some $c_k \in \Z \setminus\{0\}$, $t_k \in \Z$, $y_k \in \trmG[-1, 0]$ such that $ t_k \le t_s $ for $ k,s \in [1,N]$ with $k<s$. Since $ (y_k,  y_k)_{\hA} \in 1 + q(\Z[[q]]) $ (\cite[Proposition 7.9]{KKOP25A}), \eqref{Eq: zaa hA} says that $c_1=1$, $t_1=0$ and $t_k>0$ for $k\in [2,N]$, i.e., $ z \in 1+ q^{1/2}  \mathfrak{L}_0$.

(ii) By \cite[Lemma 7.1 (b)]{KKOP25A}, $\varphi_m(x)$ is $c$-invariant. Thus, 
$ q^{-N(\wt(x))/2} \varphi_0(x)$ is bar-invariant by \cite[(5.10)]{KKOP25A}.
Applying (i) to the setting $a := q^{-N(\wt(x))/2} \varphi_0(x)$ and then taking its image under $\Psi_w$, we have the assertion by Lemmas \ref{Lem: Phi_w Hw} and \ref{Lem: L0 L0}. 
\end{proof}

\vskip 1em 
{\bf Use of AI in Appendix \ref{Sec: appA}}.\ 
During the preparation of Appendix \ref{Sec: appA}, the authors used OpenAI GPT-6 Astra to assist in exploring mathematical ideas and developing arguments. The text of the appendix was written by the authors. Apart from this, AI was used only to check grammar and improve wording occasionally throughout the paper.

\vskip 2em


\end{document}